\def\article{1}
\ifnum\article=1
\documentclass{amsart}
\else
\documentclass[a4paper,oneside,10pt]{amsart}
\fi
\title{A Power Function with a Fixed Finite Gap Everywhere}
\author{Carmi Merimovich}
\ifnum\article=1
\address{
School of Mathematical Sciences
\\
Tel-Aviv University
\\
Tel-Aviv 69978
\\
ISRAEL
}
\email{carmi\_m@cet.ac.il}
\thanks{
This work is the final part of a research done towards a Ph.D. degree
at Tel-Aviv university under the supervision
of Moti Gitik. The author thanks Moti Gitik for his help with this work.
}
\fi
\date{May 17, 2000}
\subjclass{Primary 03E35, 03E55, 04A30}
\keywords{Forcing, Modified Radin forcing, Extender, Extender based forcing,
        Generalized continuum hypothesis, Singular cardinal hypothesis}

\allowdisplaybreaks

\ifnum\article=0
\usepackage{setspace}
\fi

\usepackage{amssymb}

\usepackage{pb-diagram}
\usepackage{lamsarrow}
\usepackage{pb-lams}

\theoremstyle{plain}
\newtheorem{theorem}{Theorem}[section]
\newtheorem*{theorem*}{Theorem}
\newtheorem{lemma}[theorem]{Lemma}
\newtheorem{proposition}[theorem]{Proposition}
\newtheorem{claim}[theorem]{Claim}
\newtheorem{corollary}[theorem]{Corollary}

\theoremstyle{definition}
\newtheorem{definition}[theorem]{Definition}
\newtheorem*{definition*}{Definition}

\theoremstyle{remark}

\newtheorem{note}[theorem]{Note}

\numberwithin{equation}{theorem}

\newcommand{\cA}{{\mathcal{A}}}

\newcommand{\ha}{\aleph}

\newcommand{\ga}{\alpha}
\newcommand{\gb}{\beta}
\newcommand{\gc}{\chi}

\newcommand{\gD}{\Delta}
\newcommand{\gee}{\epsilon}

\newcommand{\gga}{\gamma}
\newcommand{\gh}{\eta}

\newcommand{\gj}{\varphi}
\newcommand{\gk}{\kappa}
\newcommand{\gl}{\lambda}
\newcommand{\gm}{\mu}
\newcommand{\gn}{\nu}

\newcommand{\gp}{\pi}

\newcommand{\gs}{\sigma}
\newcommand{\gt}{\tau}

\newcommand{\gw}{\omega}
\newcommand{\gx}{\xi}
\newcommand{\gy}{\psi}
\newcommand{\gz}{\zeta}

\newcommand{\func }{\mathord{:}}
\newcommand{\VN}[1]{\widehat{#1}}
\newcommand{\GN}[1]{\Dot{#1}}
\newcommand{\CN}[1]{\Tilde{#1}}
\newcommand{\satisfies}{\vDash}
\newcommand{\union}{\cup}
\newcommand{\bigunion}{\bigcup}
\newcommand{\intersect}{\cap}
\newcommand{\bigintersect}{\bigcap}
\newcommand{\forces}{\mathrel\Vdash}
\newcommand{\incompatible}{\perp}
\newcommand{\compatible}{\parallel}
\newcommand{\decides}{\mathrel\Vert}

\newcommand{\subelem}{\prec}

\newcommand{\append}{\mathop{{}^\frown}}
\newcommand{\restricted}{\mathord{\restriction}}
\newcommand{\upto}{\mathord{<}}
\newcommand{\downto}{\mathord{>}}
\newcommand{\power}[1]{\lvert#1\rvert}
\newcommand{\Pset}{{\mathcal{P}}}
\newcommand{\On}{\ensuremath{\text{On}}}
\newcommand{\Card}{\ensuremath{\text{Card}}}
\newcommand{\ordered}[1]{\ensuremath{\langle #1 \rangle}}

\newcommand{\set}[1]{\ensuremath{\{ #1 \}}}
\newcommand{\setof}[2]{\ensuremath{\{ #1 \mid #2 \}}}
\newcommand{\ordof}[2]{\ensuremath{\ordered{ #1 \mid #2 }}}
\newcommand{\formula}[1]{{}^{\ulcorner} #1 {}^{\urcorner}}
\newcommand{\formulal}[1]{{}^{\ulcorner} #1 }
\newcommand{\formular}[1]{ #1 {}^{\urcorner}}

\DeclareMathOperator{\len}{l}

\DeclareMathOperator{\crit}{crit}
\DeclareMathOperator{\dom}{dom}
\DeclareMathOperator{\cf}{cf}
\DeclareMathOperator{\Col}{Col}
\DeclareMathOperator{\C}{C}
\DeclareMathOperator{\supp}{supp}
\DeclareMathOperator{\Ult}{Ult}
\DeclareMathOperator{\limdir}{lim\ dir}
\DeclareMathOperator{\id}{id}
\DeclareMathOperator{\Suc}{Suc}
\DeclareMathOperator{\Lev}{Lev}
\DeclareMathOperator{\ran}{ran}
\DeclareMathOperator*{\dintersect}{\triangle}
\newcommand{\dsintersect}{\sideset{}{^0}\dintersect}

\DeclareMathOperator{\mc}{mc}

\newcommand{\Age}{\mathrel{\mathord{\ge}_{\cA}}}
\newcommand{\Ale}{\mathrel{\mathord{\le}_{\cA}}}
\newcommand{\Egt}{\mathrel{\mathord{>}_{\Es}}}
\newcommand{\Elt}{\mathrel{\mathord{<}_{\Es}}}
\newcommand{\Ege}{\mathrel{\mathord{\ge}_{\Es}}}

\newcommand{\gas}{{\ensuremath{\Bar{\ga}}\/}}
\newcommand{\gbs}{{\ensuremath{\Bar{\gb}}\/}}
\newcommand{\ggs}{{\ensuremath{\Bar{\gga}}\/}}
\newcommand{\ges}{{\ensuremath{\Bar{\gee}}\/}}
\newcommand{\gnv}{{\ensuremath{\Vec{\gn}}\/}}
\newcommand{\gmv}{{\ensuremath{\Vec{\gm}}\/}}
\newcommand{\gns}{{\ensuremath{\Bar{\gn}}\/}}
\newcommand{\gms}{{\ensuremath{\Bar{\gm}}\/}}

\newcommand{\Es}{{\ensuremath{\bar{E}}\/}}
\newcommand{\Fs}{{\ensuremath{\bar{F}}\/}}
\newcommand{\Pe}{{\ensuremath{P_{\ges}\/}}}

\newcommand{\Peii}{{\ensuremath{P_{\ges_2}\/}}}
\newcommand{\PE}{{\ensuremath{P_{\Es}\/}}}
\newcommand{\PES}{{\ensuremath{P^*_{\Es}\/}}}

\newcommand{\ugk}{{\upto \gk}}
\newcommand{\dgk}{{\downto \gk}}
\newcommand{\VS}{V^*}
\newcommand{\MS}{{M^*}}
\newcommand{\NE}{{N^\Es}}
\newcommand{\NS}{{N^*}}
\newcommand{\NSE}{{N^{*\Es}}}
\newcommand{\MSt}{{M^*_\gt}}
\newcommand{\Mt}{{M_\gt}}
\newcommand{\MtS}{{M^*_\gt}}

\newcommand{\MtE}{{M^{\Es}_\gt}}

\newcommand{\MStE}{{M^{* \Es}_\gt}}

\newcommand{\ME}{{M_{\Es}}}
\newcommand{\MES}{{M^*_{\Es}}}
\newcommand{\MSE}{{M^*_{\Es}}}

\begin{document}

\begin{abstract}
We give an application of our extender based
Radin forcing to cardinal arithmetic.
Assuming $\gk$ is a large enough cardinal
we construct a model satisfying $2^{\gk} = \gk^{+n}$
together with $2^\gl=\gl^{+n}$
for each cardinal $\gl < \gk$,
where $0 < n < \gw$.
The cofinality of $\gk$ can be set arbitrarily or
$\gk$ can remain inaccessible.

When $\gk$ remains an inaccessible, $V_\gk$ is a model of ZFC satisfying
$2^\gl = \gl^{+n}$ for all cardinals $\gl$.
\end{abstract}
\ifnum\article=1
\maketitle
\else
\onehalfspacing
\begin{titlepage}
\vfill
\begin{center}
\Huge {\bf A Power Function with a Fixed Finite Gap Everywhere}
\end{center}
\vfill
\begin{center}
\large \bf 
        Thesis submitted for the degree ``Doctor of Philosophy''
        \\
        by
\end{center}
\begin{center}
\Large \bf Carmi Merimovich
\end{center}
\vfill
\begin{center}
\large \bf 
        Submitted to the Senate of Tel-Aviv University
        \\
        November 1999
\end{center}
\end{titlepage}
\begin{titlepage}
\vfill
\begin{center}
\Huge {\bf A Power Function with a Fixed Finite Gap Everywhere}
\end{center}
\vfill
\begin{center}
\large \bf 
        Thesis submitted for the degree ``Doctor of Philosophy''
        \\
        by
\end{center}
\begin{center}
\Large \bf Carmi Merimovich
\end{center}
\vfill
\begin{center}
\large \bf 
        Submitted to the Senate of Tel-Aviv University
        \\
        November 1999
\end{center}
\end{titlepage}
\begin{titlepage}
\vspace*{\fill}
\begin{center}
        \Large \bf
        This work was carried out under the supervision of 
        \\
        Prof.\  Moti Gitik
\end{center}
\vspace*{\fill}
\end{titlepage}
\begin{titlepage}
\hspace*{\fill}
\parbox{0.4\textwidth}{
        \large
        I thank Prof.\  Moti Gitik
        for the infinite amount of ideas, suggestions
        and help he gave me
        while I was working on this research.
}
\end{titlepage}
\tableofcontents
\fi
\ifnum\article=0
\newpage
\fi
\section{Introduction}
Investigation of the power function is as old as set theory itself.
Already Georg Cantor \cite{Cantor2} proposed CH (that is $2^{\ha_0} = \ha_1$).
As is well known, Cantor was not able to prove his hypothesis.
When Kurt G\"{o}del \cite{Godel} introduced his constructible universe $L$,
the first of set theory's \emph{inner models}, he was able to
prove that GCH (that is $2^{\ha_\ga} = \ha_{\ga+1}$ for all ordinals $\ga$)
is consistent with ZFC. That is, it is not
`unsafe' to assume GCH. At least not more so than ZFC.
Still, strictly speaking, the power function behavior was not determined.
We note that
CH is a \emph{local} hypothesis on the power function while GCH is 
a \emph{global}
one. This difference is quite important from current day point of view.

When Paul Cohen \cite{Cohen} ushered  Forcing into set theory, 
he was able to show, among other things,
that CH is independent of ZFC. In fact, his technique showed that
if $\gl$ is a regular cardinal we can make $2^\gl$ almost any
cardinal satisfying very weak restrictions. This result, more or less,
rendered the local behavior of the power function on regular cardinals
`uninteresting': Practically everything goes.

Somewhat before Forcing technology came into scene Dana Scott \cite{Scott}
proved that if $\gl$ is a measurable cardinal violating GCH then 
GCH is violated below $\gl$
on a measure $1$ set. Hence large cardinal axioms impose some
structure on the power function. We note that the above result of
Cohen does not take any additional structures into consideration.
That is, for example, if $\gl$ is measurable and we enlarge its power set
we loose its measurability in the generic extension.

Combining many instances of Cohen's construction in a clever way, 
William Easton \cite{Easton} gave a global result:
We can set the power function on all regular cardinals to practically
everything. In Easton's model the size of the power set of 
the singular cardinals 
is the lowest possible. This behavior became known as the Singular
Cardinal Hypothesis.

Generating a gap on a singular cardinal needed a more advanced methods.
The first result in this direction was a combination of
methods by Jack Silver and Karel Prikry. Silver 
        \cite{ReverseEastonNotes, IteratedForcing} 
showed how to trade a supercompact
for a measurable violating GCH. Prikry \cite{Prikry} showed how to change 
the cofinality of a measurable cardinal to $\gw$ without collapsing cardinals.
Hence we have the first example of a singular strong limit cardinal
violating GCH. We note that this violation is rather `far away'.

Menachem Magidor \cite{SingularCardinalProbI, SingularCardinalProbII} 
showed, using a supercompact cardinal, that it is 
possible that GCH will fail on
the first singular strong limit cardinal. Moreover, starting with
a huge cardinal we can have
the first failure of GCH at a singular strong limit.

At this point the general impression was that there can be arbitrary
behavior on the singular cardinals.
Then Silver \cite{SingularDepends} gave the following 
surprising
result: If $\gl$ is a singular cardinal of uncountable cofinality 
with GCH holding below $\gl$ then GCH holds at $\gl$.
Improving on the above,
Fred Galvin and Andr\'{a}s Hajnal \cite{GalvinHajnal} showed that the behavior of the 
power function
on a singular cardinal of uncountable cofinality, $\gl$, is tightly linked 
to the behavior of the power function below $\gl$.

As can be seen, the methods to get a gap on a singular cardinal
started with some large cardinal. Ronald Jensen \cite{NoLargeCardinal} 
proved that this is a necessary starting point. Specifically he proved 
that if $\lnot 0^{\dagger}$ then SCH.

With these results
the investigation of the power function has transformed its form
to the current day view:
The behavior of the power function on the singular cardinal is linked
to the existence of large cardinals. Our aim is to find restrictions
on the power function where there are ones.
When there are no restrictions we should find equiconsistency
results between existence of large cardinals
and possible behaviors of the power function.

We outline some of the known facts.

Hugh Woodin was first to use hyper measurable cardinals for 
$\lnot \text{SCH}$ results.
Continuing at the same level,
Moti Gitik and Magidor \cite{PrikryExtender} presented a forcing notion,
that used a strong cardinal $\gl$,
to  blow $2^\gl$ to whatever size prescribed together
with making $\cf \gl = \gw$ and keeping GCH below $\gl$.
Hence without some further assumptions we can not restrict the size of
the power set of singulars of countable cofinality.
Indeed a modification of Gitik and Magidor's forcing  to $\ha_\gw$
gives Magidor's original result: GCH below $\ha_\gw$ and
$2^{\ha_\gw} = \ha_{\ga+1}$ for $\ga < \gw_1$. 
(albeit from a considerably lower
large cardinal assumption).
It is still not known if it is possible to get such a model
with $2^{\ha_\gw} \geq \ha_{\gw_1}$.

On the other hand, we do  have restrictions at $\ha_\gw$.
Saharon Shelah showed that 
	$2^{\ha_\gw} < \min (\ha_{(2^{\ha_0})^+}, \ha_{\gw_4})$.

An interesting twist, connected with the above, is a work of Gitik and
Bill Mitchell \cite{GitikM1996} showing that if there is no inner model with a strong cardinal
then $2^{\ha_\gw} < \ha_{\gw_1}$. Hence, if Shelah's bound, $\ha_{\gw_4}$,
is optimal we need a stronger large cardinal to approach this bound than
the ones used to get up to $\ha_{\gw_1}$.

Mitchell \cite{KSeqMeasures} showed how to get high order measurable cardinals
from $\lnot \text{SCH}$ and later on
Gitik, building on results of Mitchell, Shelah and Woodin,
 had pinpoint $\lnot \text{SCH}$ to be equiconsistent with
$\text{o}(\gk) = \gk^{++}$.

Just by looking at these few results --- there are many more known
and we still lack a lot of information --- 
we can see that the situation
on the singular cardinals is much more complicated than the situation
on the regular cardinals.
Of course getting a full result as in Easton's one is beyond our
reach at this point. There are, however, several results
concerning the global behavior on all the cardinals.

Matt Foreman and Woodin \cite{GCHNoWhere}, starting from a supercompact,
constructed a model
in which GCH fails everywhere. In their model the gap between $\gl$ and $2^\gl$
is infinite and is not fixed for all $\gl$. We know that we can not have
a fixed infinite gap for all $\gl$. A later  unpublished modification of the 
construction, due to Woodin, gave a model satisfying $2^\gl = \gl^{++}$
for all cardinals $\gl$. A referee of our thesis brought to our
attention that Woodin generalized the construction and for each
$1 < n < \gw$ got a model
having a power function with a fixed gap $n$.
Unfortunately, this result, also, was never published.

James Cummings \cite{JamesPublishedPhd}, starting from a $\Pset^3(\gk)$-hypermeasurabe constructed 
a model satisfying
$2^\gl = \gl^{++}$ for $\gl$ limit cardinals and GCH everywhere else.

\vspace{\baselineskip}
As the basic idea of our work is a generalization of
Gitik and Magidor's one we 
elaborate more on it.
Until Gitik and Magidor's work the major theme in generating a gap
on a singular cardinal was as follows.
The power set of a large cardinal is blown up sacrificing some of its size
but not all of it. (That is, starting from a supercompact we are left
with a measurable).
Then one of the known forcings for changing cofinality
is applied to this cardinal retaining the gap on the size of the power set.
Gitik and Magidor forcing does \emph{both} tasks in \emph{one} step.
In essence they found a method to add many Prikry sequences at once
using an extender of a prescribed size. Hence in one step they
blow up the size of the power set and change the cofinality to $\gw$.

While Prikry forcing can change the cofinality of a measurable
cardinal to $\gw$ we have
Magidor forcing \cite{ChangeCardinalsCofinality} that adds a new club of a prescribed order type, $\ga$,
using a coherent sequence of measures.
Miri Segal \cite{MiriMsc} combined this forcing with the extender idea of
Gitik and Magidor,
hence she was able to add a prescribed number of $\ga$-sequences.
That is in one step the power set is blown up and the cofinality is changed
to $\cf \ga$.
Of course her starting point was a coherent sequence of extenders.

The marriage of Radin forcing \cite{RadinPublishedPhd, WeakClub} with extenders 
was done by us \cite{MeselfPublishedPhdI}.
Radin forcing adds to Magidor's one a new ingredient.
 It enables us to add a club to a large cardinal $\gk$ while
keeping $\gk$ inaccessible (and even more). Indeed this behavior remains
possible when we add many Radin sequences using an extender.
This gave us the possibility to control the power function
on a club while keeping the cardinal inaccessible. All this in
\emph{one} step.
The different properties of Radin forcing are related to the
length of the measure sequence used to define the forcing.
In our extender based Radin forcing we start from an extender sequence
and the properties are controlled by the length of the extender \emph{sequence}
(controlling the properties of $\gk$)
and the length of the extender (controlling $2^\gk$).

\vspace{0.125 in}
This paper is a step in the investigation of the global behavior
of the power function. 
The forcing we present should be viewed as a template enabling
the construction of models with many different power functions.
We stress that our main point is \emph{not}
a backwards Easton iteration for blowing up power sets cardinal by
cardinal, followed by choosing cardinals after their power set was enlarged.
We go the other way: We choose cardinals and then we blow up their power
set. I.e., we do not blow up the power set of cardinals which are
of no interest for us.

As a specific example we show that for each $1 < n < \gw$
it is possible to have
a power function satisfying $2^\gl = \gl^{+n}$ for all cardinals $\gl$.
It is possible to get this behavior assuming there is $\gk$ that
is $\gk^{+n+1}$-strong. (A somewhat weaker assumption is enough).
In section \ref{ConcludingRemarks} we describe the general behavior
possible with our technique and give more examples.

Our starting point is the extender based Radin forcing.
As mentioned earlier it enables us to control the behavior
of the power function on a club of $\gk$  while keeping $\gk$
inaccessible. Hence $V_\gk$ of the generic extension is a model
of ZFC with a club of cardinals satisfying a prescribed power function.

In this work we control also the cardinals that are not in the club.
We either set their power to the prescribed value or collapse them.

The specific example we construct is $2^\gl = \gl^{+3}$
everywhere. Taking any $n < \gw$ is exactly the same.
Of course other behaviors are possible with this method.

So we start with $\gk$ large enough as witnessed by
$j \func V \to M$.
We construct from $j$ extender based Radin forcing that 
sets $2^\gk = \gk^{+3}$, shoots a club
through $\gk$ and for each $\gl$ in the club we have $2^\gl = \gl^{+3}$.
We are careful to make sure that $\gk$ remains inaccessible in the
generic extension.
In order to control the behavior outside the club we 
add along the Radin sequence generated by the normal measure
other forcing notions. That is if $\gm_1, \gm_2$ are $2$ successive points
in this Radin sequence then we force with 
$\Col(\gm_1^{+6}, \gm_2) \times \C(\gm_1^{+4}, \gm_2^+) \times
        \C(\gm_1^{+5}, \gm_2^{++}) \times
        \C(\gm_1^{+6}, \gm_2^{+3})$ as defined in some inner model.
Note that we actually collapse cardinals in the Radin generic sequence.
In order to allow for a Prikry like condition we need a generic
filter from which the above forcing conditions will come.
That is we will have $I \in V$ that is
$\Col(\gk^{+6}, j(\gk)) \times \C(\gk^{+4}, j(\gk)^+) \times
        \C(\gk^{+5}, j(\gk)^{++}) \times
        \C(\gk^{+6}, j(\gk)^{+3})$ (as defined in some inner model of $M$)
generic over $M$. It was pointed out
by Woodin that a generic filter for forcing of this type 
can be generated through the normal measure (e.g. if $j$ is witnessing
that $\gk$ is $\gk^{+3}$-strong and not $\gk^{+4}$-strong then
the normal measure
generates a $\Col(\gk^{+4}, j(\gk))_M$-generic filter over $M$. However,
it does not generate a $\Col(\gk^{+3}, j(\gk))_M$-generic filter over $M$.).

At this stage $V_\gk$ of the generic extension almost satisfies 
our requirements. The only problem is that we do not have a gap of 
$3$ on $\gm_1^+$,
$\gm_1^{++}$, $\gm_1^{+3}$. 
Of course the naive approach is to add also
	$\C(\gm_1^{+}, \gm_1^{+4}) \times
        \C(\gm_1^{++}, \gm_1^{+5}) \times
        \C(\gm_1^{+3}, \gm_1^{+6})$ along the normal Radin sequence.
However, in order to have the Prikry condition we need a
$\C(\gk^{+}, \gk^{+4}) \times
        \C(\gk^{++}, \gk^{+5}) \times
        \C(\gk^{+3}, \gk^{+6}))_M$-generic over $M$.
As this forcing is `below' the extender length, hence `sees' much of
$V$, we do not have such a generic.
One solution to this problem is to force such a generic into $V$.
We did a similar thing in \cite{MeselfPublishedMsc} and the amount of
technical difficulties we had to overcome  was overwhelming.
(And to a large extent blurred the simple idea).
A second solution, adopted here, is to do
a preparation forcing  making a gap of $3$
on $\gl^+$, $\gl^{++}$, $\gl^{+3}$ for each $\gl < \gk$ inaccessible.
By making sure that the normal Radin sequence pass only through 
inaccessibles we will get the prescribed behavior.

This work clarified to the author many points from the simpler
\cite{MeselfPublishedPhdI}. As a result many of the proofs 
appearing here, in more complex setting, are  simpler than the ones 
appearing in \cite{MeselfPublishedPhdI}.

The reader should be fluent with forcing technology and large cardinals
methods. We assume  that \cite{PrikryExtender} is known. Especially
the `nice extender' built there.
Knowledge of \cite{MeselfPublishedPhdI}
makes reading of this work easier.

The structure of this work is as follows.
We start from a universe $V^*$ satisfying GCH that has
a suitably large cardinal $\gk$ as witnessed by an elementary embedding $j$.
In section \ref{ExtenderSequenceDefinition}, taken almost verbatim
from \cite{MeselfPublishedPhdI}, we define extender sequences.
In section \ref{PreparationForcing} we construct $V$ as a generic extension of
$V^*$ with a changed power function on
the $3$ first successors of inaccessible cardinals below $\gk$. 
We locate in $V$ a generic filter that is used later to 
change the power function on the next $3$ cardinals and collapse
all others.
In section \ref{ExtenderSequenceRedefinition} we incorporate the generic
filters located in the previous section into the definition of
extender sequence. This revised extender sequence is the one we
use in the rest of the paper. In section \ref{PEForcing} we define
our Modified Extender Based Radin Forcing, $\PE$.
In section \ref{BasicProperties} we give some basic properties
of the just defined forcing notion.
Section \ref{HomogenDense} is dedicated to the proof of the
homogeneity of dense open subsets of $\PE$. This property
plays a central role in later analysis. 
In section \ref{Prikry'sCondition} we use the homogeneity to prove
Prikry's condition for $\PE$.
Section \ref{Properness} is used to show how to get
generic filter over elementary submodels.
In section \ref{CardinalStructure} we show what cardinals
are not collapsed by $\PE$ and what their power is.
Section \ref{ConsistencyTheorem} is just the stating of the
consistency theorem we proved.
Section \ref{ConcludingRemarks} is a list of points for later
research and indication of preliminary work we have.
\ifnum\article=0
\newpage
\fi
\section{Extender sequences} \label{ExtenderSequenceDefinition}
\ifnum\article=0
\enlargethispage*{10pt}
\fi
Suppose we have an elementary embedding $j\func V^* \to M^* \supset 
V^{V^*}_\gl$, $\crit(j)=\gk$.
The value of $\gl$ is determined later, according to the different
applications we will have.

Construct from $j$ a nice extender as  in \cite{PrikryExtender}:
\begin{align*}
E(0) = \ordered{\ordof{E_\ga(0)}{\ga \in \cA},
                \ordof{\gp_{\gb,\ga}}{\gb \Age \ga, \  \ga,\gb \in \cA}}.
\end{align*}
We remind the reader what are the properties of this extender:
\begin{enumerate}
\item $\cA \subseteq \power{V^{V^*}_\gl}\setminus\gk$,
\item $\power{\cA} = \power{V^{V^*}_\gl}$,
\item $\ordered{\cA, \Ale}$ is a $\gk^+$-directed partial order,
\item $\forall \ga,\gb\in \cA$ $\gb \Age \ga \implies$
	$\gp_{\gb,\ga}\func V_\gk \to V_\gk$,
\item $\gk \in \cA$,
\item $\forall \ga \in \cA$ $\gk \Ale \ga$. We write $\gp_{\ga,0}$ instead
        of $\gp_{\ga,\gk}$,
\item $\forall \ga,\gb\in \cA$ $\forall \gn < \gk$
        $\gn^0 = \gp_{\ga,0}(\gn) = \gp_{\gb,0}(\gn)$,
\item $\forall \ga,\gb \in \cA$ $\gb \Age \ga \implies$ $\forall \gn <\gk$
        $\gp_{\gb,0}(\gn) = \gp_{\ga,0}(
                                \gp_{\gb,\ga}(\gn))$,
\item $\forall \ga,\gb,\gga \in \cA$ $\gga \Age \gb \Age \ga \implies$
        $\exists X \in E_\gga(0)$
        $\forall \gn \in X$
        $\gp_{\gga,\ga}(\gn) = \gp_{\gb,\ga}(
                                \gp_{\gga,\gb}(\gn))$.
\end{enumerate}
If, for example, we need $\power{E(0)} = \gk^{+3}$ then, under $\text{GCH}$,
we require $\gl = \gk+3$.
A typical big set in this extender concentrates on singletons.

If $j$ is not sufficiently closed, then $E(0) \notin M^*$ and the construction
stops. We  set
\begin{align*}
\forall \ga \in \cA \ \Es_\ga = \ordered{\ga,E(0)}.
\end{align*}
We say that $\Es_\ga$ is an extender sequence of length $1$.
($\len(\Es_\ga)=1$)

If, on the other hand, $E(0) \in M^*$ we can construct for each $\ga \in \dom E(0)$
the following ultrafilter
\begin{align*}
A \in E_{\ordered{\ga,E(0)}}(1) \iff \ordered{\ga, E(0)} \in j(A).
\end{align*}
Such an $A$ concentrates on elements of the form $\ordered{\gx, e(0)}$
where $e(0)$ is an extender on $\gx^0$
and $\gx \in \dom e(0)$. Note that $e(0)$ concentrates on singletons
below $\gx^0$. If, for example, $\power{E(0)} = \gk^{+3}$ then on a large
set we have $\power{e(0)} = (\gx^{0})^{+3}$.

We define $\gp_{\ordered{\gb,E(0)},\ordered{\ga,E(0)}}$ as
\begin{align*}
\gp_{\ordered{\gb,E(0)},\ordered{\ga,E(0)}}(\ordered{\gx, e(0)}) = 
        \ordered{\gp_{\gb,\ga}(\gx), e(0)}.
\end{align*}
From this definition we get
\begin{align*}
j(\gp_{\ordered{\gb,E(0)},\ordered{\ga,E(0)}})(\ordered{\gb, E(0)}) = 
        \ordered{\ga, E(0)}.
\end{align*}
Hence we have here an extender
\begin{align*}
E(1) = \ordered{
                \ordof{E_{\ordered{\ga,E(0)}}(1)}{\ga \in {\cA}},
                \ordof{\gp_{\ordered{\gb,E(0)},\ordered{\ga,E(0)}}}
                        {\gb \Age \ga, \  \ga,\gb \in \cA}
        }.
\end{align*}
Note that the difference between $\gp_{\gb,\ga}$ and
$\gp_{\ordered{\gb,E(0)},\ordered{\ga,E(0)}}$ is quite superficial.
We can define $\gp_{\ordered{\gb,E(0)},\ordered{\ga,E(0)}}$
in a uniform way for both extenders. Just project the first element of
the argument using $\gp_{\gb,\ga}$.

If $\ordered{E(0),E(1)} \notin M^*$ then the construction stops.
In this case we set
\begin{align*}
\forall \ga \in \cA \ \Es_\ga = \ordered{\ga, E(0), E(1)}.
\end{align*}
We say that $\Es_\ga$ is an extender sequence of length $2$.
($\len(\Es_\ga)=2$.)

If $\ordered{E(0),E(1)} \in M^*$ then we construct the extender
$E(2)$ in the same way as we constructed $E(1)$ from $E(0)$.

The above special case being worked out we continue with the general case.
Assume we have constructed
\begin{align*}
\ordof  {E(\gt')}
        {\gt' < \gt}.
\end{align*}
If $\ordof {E(\gt')} {\gt' < \gt} \notin M^*$ then our construction stops
here. We set
\begin{align*}
\forall \ga \in \cA \ \Es_\ga = \ordof  {\ga, E(\gt')}
                        {\gt' < \gt}.
\end{align*}
and we  say that $\Es_\ga$ is
an extender sequence of length $\gt$. ($\len(\Es_\ga)= \gt$.)

If, on the other hand, $\ordof {E(\gt')} {\gt' < \gt} \in M^*$ then 
we construct
\begin{multline*}
A \in E_{\ordof{\ga, E(0),\dotsc, E(\gt'),\dotsc} {\gt' < \gt}}(\gt) \iff
\\
        \ordof{\ga, E(0), \dotsc, E(\gt'), \dotsc} {\gt' < \gt}
                                         \in j(A).
\end{multline*}
Defining $\gp_{\ordof{\gb, E(0),\dotsc, E(\gt'),\dotsc}{\gt' < \gt},
\ordof{\ga, E(0),\dotsc, E(\gt'),\dotsc,}{\gt' < \gt}}$ using the first
coordinate as before gives the needed projection.

We are quite casual in writing the indices of the projections and
ultrafilters.
By this we mean that we sometimes write $\pi_{\gb,\ga}$ when we
should have written
$\gp_{\ordof{\gb, E(0),\dotsc, E(\gt'),\dotsc}{\gt' < \gt},
\ordof{\ga, E(0),\dotsc, E(\gt'),\dotsc,}{\gt' < \gt}}$
and $E_\ga(\gt)$ when we should have written
$E_{\ordof{\ga, E(0),\dotsc, E(\gt'),\dotsc,}{\gt' < \gt}}(\gt)$.

With this abuse of notation, the projection we just defined satisfies
\begin{multline*}
j(\gp_{\gb,\ga})
        (\ordof{\gb, E(0), \dotsc, E(\gt'), \dotsc} {\gt' < \gt})=
                \ordof{\ga, E(0), \dotsc, E(\gt'), \dotsc} {\gt' < \gt}.
\end{multline*}
So we have the extender
\begin{align*}
E(\gt) = \ordered{\ordof{E_\ga(\gt)}{\ga \in {\cA}},
                \ordof{\gp_{\gb,\ga}}{\gb \Age \ga, \  \ga,\gb \in \cA}}.
\end{align*}

We let the construction run until it stops due to the extender
sequence not being in $M^*$.

\begin{definition}
$\gms$ is an extender sequence if there is an elementary embedding
$j\func V^* \to M^*$ such that $\gns$ is an extender sequence generated 
as above and $\gms = \gns \restricted \gt$ for $\gt \leq \len(\gns)$.
$\gk(\gms)$ is the ordinal at the beginning of
the sequence. (i.e. $\gk(\bar{E}_\ga)=\ga$) $\gk^0(\gms)$ is
$(\gk(\gms))^0$. (i.e. $\gk^0(\bar{E}_\ga)=\gk$)
\end{definition}
That is, we do not have to construct the extender sequence until it is
not in $M^*$. We can stop anywhere on the way.
\begin{definition}
A sequence of extender sequences $\ordered{\gms_1, \dotsc, \gms_n}$ is
be called $^0$-increasing if $\gk^0(\gms_1) < \dotsb <\gk^0(\gms_n)$.
\end{definition}
\begin{definition}
Let $\ordered{\gms_1, \dotsc, \gms_n}$ be $^0$-increasing.
An extender sequences $\gms$ is called permitted to 
$\ordered{\gms_1, \dotsc, \gms_n}$ if
         $\gk^0(\gms_n) < \gk^0(\gms)$.
\end{definition}
%
%
\begin{definition}
We  say $A \in \bar{E}_\ga$ if $\forall \gx < \len(\bar{E}_\ga) \ A\in E_\ga(\gx)$.
\end{definition}
\begin{definition}
$\Es =\ordof{\Es_\ga}{\ga \in \cA}$ is an extender sequence system
if there is an elementary embedding $j\func V^* \to M^*$ such that all
$\Es_\ga$ are extender sequences generated from $j$ as prescribed above
and $\forall \ga,\gb\in{\cA} \ \len(\Es_\ga)=\len(\Es_\gb)$. This
common length is called the length of the system, $\len(\Es)$. We
write $\Es(\gms)$ for the 
extender sequence system to which $\gms$ belongs (i.e. $\Es(\Es_\ga)=\Es$).
\end{definition}
We point out that there is a $\gk^+$-directed partial order on the
system $\Es$ inherited from $\cA$. That is $\Es_\gb \Ege \Es_\ga \iff
\gb \Age \ga$. Of course, this implies that there is $\min \Es$,
namely $\Es_\gk$. From now on we use only the order $\Ege$ (even for $\cA$)
and we write $\dom \Es$ for $\cA$.

$\Es_\ga$ is the generalization of the measure on the $\ga$ coordinate in
Gitik-Magidor forcing \cite{PrikryExtender}.
\ifnum\article=0
\newpage
\fi
\section{Preparation forcing} \label{PreparationForcing}
\ifnum\article=0
\enlargethispage*{10pt}
\fi
Suppose that we have GCH and an elementary embedding 
	$j\func  \VS \to \MS \supset \VS_{\gk+3}$, 
		$\crit(j)=\gk$.
(If we need to get, say, $2^\gk = \gk^{+5}$ then we should
start with $j\func  \VS \to \MS \supset \VS_{\gk+5}$.)
Construct from $j$ an extender sequence system, $\Es$,
and define the following embeddings 
$\forall \gt' < \gt < \len(\Es)$
\begin{align} \label{E-system}
& j_\gt\func  \VS \to \MSt \simeq \Ult(\VS, E(\gt)),
\notag \\
\qquad & k_{\gt}(j_\gt(f)(\Es_\ga\restricted \gt)) = 
		j(f)(\Es_\ga \restricted \gt),
\notag\\
& \Es \restricted \gt \in \MSt,
\\
\notag & i_{\gt', \gt}(j_{\gt'}(f)(\Es_\ga \restricted \gt')) = 
	j_\gt(f)(\Es_\ga \restricted \gt'),
\\
\notag & \ordered{\MSE,i_{\gt, \Es}} = \limdir \ordered {
		\ordof{\MSt} {\gt < \len(\Es)},
                \ordof{i_{\gt',\gt}} {\gt' \leq \gt < \len(\Es)}
		},
\end{align}
giving rise to the following commutative diagram
\begin{align*}
\begin{diagram}
\node{\VS} 
		\arrow[3]{e,t}{j} 
		\arrow{sse,l}{j_{\gt'}}
		\arrow[2]{se,l}{j_\gt} 
		\arrow{seee,t,1}{j_{\Es}}
	\node{}
	\node{}
	\node{\MS}
\\
\node{}
	\node{}
	\node{}
        \node{\MSE} 
		\arrow{n,r}{k_{\Es}}
\\
\node{}
	\node{M^*_{\gt'}}
		\arrow[2]{ne,t,3}{k_{\gt'}} 
		\arrow{nee,t,2}{i_{\gt', \Es}} 
		\arrow{e,b}{i_{\gt', \gt}} 
	\node{\MSt = \Ult(\VS, E(\gt))} 
		\arrow[1]{ne,b}{i_{\gt, \Es}} 
		\arrow{nne,b,1}{k_{\gt}} 
\end{diagram}
\end{align*}
Note that in general $\len(\Es)$ can be very large. However, the
requirement $\Es \restricted \gt \in \MSt$ imposes limit on
$\len(\Es)$. In our case to $\len(\Es) \leq \gk^{+4}$.
In this work we are more strict than that and require
$\len(\Es) < \gk^{+4}$. 
In fact, if $\len(\Es) = \gk^{+4}$ then the forcing notion we define
using it, $\PE$,
is isomorphic to the forcing defined with $\Es \restricted \gt$
for $\gt < \gk^{+4}$.
For the specific result we are aiming to
it is enough to have $\len(\Es) = \gk^+$. 

We use the following,
quite standard, notation
\begin{align*}
& \C(\gl, \gm) = \setof {f} {f \func A \to 2,\ A \in [\gm]^{\upto \gl}},
\\
& \Col(\gl, \gm) = \setof {f} {f \func A \to \gm,\ A \in [\gl]^{\upto \gl}}.
\end{align*}
Note that there is a change from previous works of this type. We  use
$\Col(\gl, \gm)$ and not the Levy collapse $\Col(\gl, \upto \gm)$.
As will be seen (much) later this helped us a lot.

Now that the setting is clear and before we start with the preparation,
a note is in order.
The main point of this work is the forcing $\PE$ described in 
section \ref{PEForcing} and the technicalities of the current section are
somewhat `off track'. A reader willing to accept a somewhat
weaker result than the one we stated can do without the preparation.
For example, to prove
\begin{theorem*}
If there is $\Es$ such that $\power{\Es} = \gk^{+3}$,
$\cf \len(\Es) > \gk$ then there is
a model containing a (class) club $C$ such that the cardinals
in the model are 
$\setof {\gm^+, \gm^{++}, \gm^{+3}, \gm^{+4}} {\gm \in C} \union
	\lim C$
and the power function is
\begin{align*}
2^\gl = \begin{cases}
	\gm^{+3} & \gm \in \lim C,\ \gl \in \set{\gm, \gm^{+}, \gm^{++}} 
	\\
	\gl^{+} & \text{otherwise}
	\end{cases}.
\end{align*}
\end{theorem*}
it is enough to locate in $V^*$ a $\Col(\gk^{+4}, j_\Es(\gk))_{\MES}$-generic
filter over $\MES$ and to jump to the next section (skimming through the
definitions of subsection \ref{IterateUltrapower}).
Similarly we can prove
\ifnum\article=0
\enlargethispage*{10pt}
\fi
\begin{theorem*}
If there is $\Es$ such that $\power{\Es} = \gk^{+3}$,
$\cf \len(\Es) > \gk$ then there is
a model containing a (class) club $C$ such that the  cardinals
in the model are $\setof{\gm^+, \dotsc, \gm^{+6}} {\gm \in C} \union \lim C$
and the power function is
\begin{align*}
2^\gl = \begin{cases}
	\gm^{+3} & \gm \in \lim C,\ \gl \in \set{\gm, \gm^{+}, \gm^{++}} 
	\\
	\gl^{+3} & \gm \in C,\ \gl \in \set{\gm^{+4}, \gm^{+5}, \gm^{+6}} 
	\\
	\gl^{+} & \text{otherwise}
	\end{cases}.
\end{align*}
\end{theorem*}
it is enough to find in $V^*$ a 
$(\Col(\gk^{+6}, j_\Es(\gk)) \times \C(\gk^{+4}, j_\Es(\gk)^+) \times
	\C(\gk^{+5}, j_\Es(\gk)^{++}) \times
	\C(\gk^{+6}, j_\Es(\gk)^{+3}))_{\MES}$-generic
filter over $\MES$ and to jump to the next section
(again, skimming through the
definitions of subsection \ref{IterateUltrapower}).
It is easy to construct both generic filters by going
through $\Ult(V^*, E_\gk(0))$.

The above said, we continue to the preparation.
It will give us a very rough approximation of the power set
function we seek. For each $\gn \leq \gk$ an inaccessible,
we will have
\begin{align*}
& 2^{(\gn^{+})} = \gn^{+4},
\\
& 2^{(\gn^{++})} = \gn^{+5},
\\
& 2^{(\gn^{+3})} = \gn^{+6}.
\end{align*}
In the generic extension there is
a generic filter for the forcing notion
\begin{align*}
(\Col(\gk^{+6}, j_\gt(\gk)) \times &
	\C(\gk^{+6}, j_\gt(\gk)^{+3}) \times
	\C(\gk^{+5}, j_\gt(\gk)^{++}) \times
	\C(\gk^{+4}, j_\gt(\gk)^{+}))_{\MSt[G_{\upto \gk}]}
\end{align*}
over $\MSt[G^\gt][H^\gt]$, where $G_{\upto \gk}$ $G^\gt$, $H^\gt$
will be defined later.
\subsection{Reverse Easton forcing for pulling in the needed generic 
		filters}
We make a reverse Easton forcing and lift the diagram
\begin{align*}
\begin{diagram}
\node{\VS} 
		\arrow[3]{e,t}{j_{\Es}}
		\arrow{se,b}{j_{\gt'}}
		\arrow{see,t}{j_\gt}
        \node{}
        \node{}
        \node{\MSE}
\\
\node{}
	\node{M^*_{\gt'}}
		 \arrow{nee,t,3}{i_{\gt', \Es}} 
		 \arrow{e,b}{i_{\gt', \gt}} 
	\node{\MSt} 
		\arrow[1]{ne,b}{i_{\gt, \Es}} 
\end{diagram}
\end{align*}
We define the following reverse Easton iteration 
$\ordof{P_\gn, \GN{Q}_\gn}{\gn \leq \gk}$:
When $\gn$ is accessible $\GN{Q}_\gn = \GN{1}$
and when $\gn$ is inaccessible
\begin{align*}
\GN{Q}_\gn = \C(\VN{\gn}^+, \VN{\gn}^{+4}) \times
                \C(\VN{\gn}^{++}, \VN{\gn}^{+5}) \times
                \C(\VN{\gn}^{+3}, \VN{\gn}^{+6}).
\end{align*}
We set $\GN{P}_{\downto \gn}$ to be the forcing notion name satisfying
\begin{align*}
P_\gk = P_\gn * \GN{Q}_\gn * \GN{P}_{\downto \gn}.
\end{align*}
We factor through the normal ultrafilter to get
\begin{align*}
\begin{aligned}
\begin{diagram}
\node{\VS} 
		\arrow{e,t}{j_\Es}
		\arrow{se,t}{j_\gt}
		\arrow{s,l}{i_U}
        \node{\MSE}
\\
\node{\NS \simeq \Ult(\VS, U)}
		 \arrow{e,b}{i_{U, \gt}}
		 \arrow{ne,b}{i_{U, \Es}}
        \node{\MSt}
		 \arrow{n,b}{i_{\gt, \Es}}
\end{diagram}
\end{aligned}
\begin{aligned}
\qquad
\begin{split}
& U = E_\gk(0),
\\
& i_U \func  \VS \to \NS \simeq \Ult(\VS, U),
\\
& i_{U, \gt}(i_U(f)(\gk)) = j_\gt(f)(\gk),
\\
& i_{U, \Es}(i_U(f)(\gk)) = j_\Es(f)(\gk).
\end{split}
\end{aligned}
\end{align*}
For later convenience we set
\begin{align*}
& \ordof{P^U_\gn, \GN{Q}^U_\gn}{\gn \leq i_U(\gk)} = 
	i_U(\ordof{P_\gn, \GN{Q}_\gn}{\gn \leq \gk}),
\\
& P^U_{i_U(\gk)} = P^U_\gk * \GN{Q}^U_\gk * \GN{P}^U_{\downto \gk},
\\
& \ordof{P^\gt_\gn, \GN{Q}^\gt_\gn}{\gn \leq j_\gt(\gk)} = 
	j_\gt(\ordof{P_\gn, \GN{Q}_\gn}{\gn \leq \gk}),
\\
& P^\gt_{j_\gt(\gk)} = P^\gt_\gk * \GN{Q}^\gt_\gk * \GN{P}^\gt_{\downto \gk},
\\
& \ordof{P^\Es_\gn, \GN{Q}^\Es_\gn}{\gn \leq j_\Es(\gk)} = 
	j_\Es(\ordof{P_\gn, \GN{Q}_\gn}{\gn \leq \gk}),
\\
& P^\Es_{j_\Es(\gk)} = P^\Es_\gk * \GN{Q}^\Es_\gk * \GN{P}^\Es_{\downto \gk}.
\end{align*}
We note that $P_\gk = P^U_\gk = P^\gt_\gk = P^\Es_\gk$.
Let $G_{\upto \gk}$ be $P_\gk$-generic over $\VS$.
As $\crit(i_{U, \Es}) = \crit(i_{U, \gt}) = (\gk^{++})_{\NS}$ and 
$i''_{U, \Es}G_{\upto \gk} = i''_{U, \gt} G_{\upto \gk}= G_{\upto \gk}$ 
we have the lifting
\begin{align*}
\begin{diagram}
\node{} 
        \node{\MSE[G_{\upto \gk}]}
\\
\node{\NS[G_{\upto \gk}]}
		 \arrow{e,b}{i_{U, \gt}}
		 \arrow{ne,b}{i_{U, \Es}}
        \node{\MSt[G_{\upto \gk}]}
		 \arrow{n,b}{i_{\gt, \Es}}
\end{diagram}
\end{align*}
In order to find a generic filter for stage $\gk$ we need 
\ref{EmbeddedForcing}.
The following series of propositions provides a proof of \ref{EmbeddedForcing}.
\begin{proposition}
If $\gl < j_\Es(\gk)$ then $\cf \gl \leq \gk^{+3}$.
\end{proposition}
\begin{proof}
As $\len(\Es) < \gk^{+4}$, $\power{\Es} = \gk^{+3}$ we get that
$j_\Es(\gk) < \gk^{+4}$. Hence $\gl < \gk^{+4}$.
\end{proof}
\begin{proposition}
$\C(\gk^{+}, \gk^{+3})_{V^*[G_\ugk]} \simeq 
	\C(\gk^{+}, \gk^{+4})_{\MSE[G_\ugk]}$.
\end{proposition}
\begin{proof}
We carry the proof inside $V^*[G_\ugk]$.
The first point to observe is that 
$\C(\gk^{+}, \gk^{+3}) = \C(\gk^{+}, \gk^{+3})_{\MSE[G_\ugk]}$.
So if $A \in \MSE[G_\ugk]$ and $\power{A}_\MSE=\gk^{+3}$ then
$\C(\gk^{+}, \gk^{+3}) \simeq \C(\gk^{+}, A)_{\MSE[G_\ugk]}$.
Let $\setof {\ga_\gx} {\gx < \gk^{+3}}$ be a club of $(\gk^{+4})_\MSE$.
Hence
	$\power{[\ga_\gx, \ga_{\gx+1})}_{\MSE[G_\ugk]} \leq \gk^{+3}$.
For simplicity of notation let us assume 
 $\C(\gk^{+}, \gk^{+3}) \simeq 
	\C(\gk^{+}, [\ga_\gx, \ga_{\gx+1}))_{\MSE[G_\ugk]}$.
The proposition follows from
\begin{align*}
\C(\gk^{+}, \gk^{+3}) \simeq 
	\linebreak[0] 
	\prod_{\gx < {\gk^{+3}}}\C(\gk^{+}, \gk^{+3}) \simeq 
	\linebreak[0]
	\prod_{\gx < \gk^{+3}} \C(\gk^{+}, [\ga_\gx, \ga_{\gx+1}))_{\MSE[G_\ugk]}.
\end{align*}
\end{proof}
\begin{proposition}
$\C(\gk^{++}, \gk^{+3})_{V^*[G_\ugk]} \simeq \C(\gk^{++}, \gk^{+4})_{\MSE[G_\ugk]}$.
\end{proposition}
\begin{proof}
We carry the proof inside $V^*[G_\ugk]$.
The first point to observe is that 
$\C(\gk^{++}, \gk^{+3}) = \C(\gk^{++}, \gk^{+3})_{\MSE[G_\ugk]}$.
So if $A \in \MSE[G_\ugk]$, $\power{A}_\MSE=\gk^{+3}$ then
$\C(\gk^{++}, \gk^{+3}) \simeq \C(\gk^{++}, A)_{\MSE[G_\ugk]}$.
Let $\setof {\ga_\gx} {\gx < \gk^{+3}}$ be a club of $(\gk^{+4})_\MSE$.
Then $\power{[\ga_\gx, \ga_{\gx+1})}_{\MSE[G_\ugk]} \leq \gk^{+3}$.
For simplicity of notation let us assume
 $\C(\gk^{++}, \gk^{+3}) \simeq \C(\gk^{++}, [\ga_\gx, \ga_{\gx+1}))_{\MSE[G_\ugk]}$.
The proposition follows from
\begin{align*}
\C(\gk^{++}, \gk^{+3}) \simeq 
	\linebreak[0] 
	\prod_{\gx < {\gk^{+3}}}\C(\gk^{++}, \gk^{+3}) \simeq 
	\linebreak[0]
	\prod_{\gx < \gk^{+3}} \C(\gk^{++}, [\ga_\gx, \ga_{\gx+1}))_{\MSE[G_\ugk]}.
\end{align*}
\end{proof}
\begin{proposition} \label{Same:2:5}
$\C(\gk^{++}, \gk^{+3})_{V^*[G_\ugk]} \simeq \C(\gk^{++}, \gk^{+5})_{\MSE[G_\ugk]}$.
\end{proposition}
\begin{proof}
We carry the proof inside $V^*[G_\ugk]$.
The point to observe is that 
if $A \in \MSE[G_\ugk]$, $\power{A}_\MSE=\gk^{+4}$ then
$\C(\gk^{++}, \gk^{+3}) \simeq \C(\gk^{++}, A)_{\MSE[G_\ugk]}$.
Let $\setof {\ga_\gx} {\gx < \gk^{+3}}$ be a club of $(\gk^{+5})_\MSE$.
Then $\power{[\ga_\gx, \ga_{\gx+1})}_{\MSE[G_\ugk]} \leq \gk^{+4}$.
For simplicity of notation let us assume
 $\C(\gk^{++}, \gk^{+3}) \simeq \C(\gk^{++}, [\ga_\gx, \ga_{\gx+1}))_{\MSE[G_\ugk]}$.
The proposition follows from
\begin{align*}
\C(\gk^{++}, \gk^{+3}) \simeq 
	\linebreak[0] 
	\prod_{\gx < {\gk^{+3}}}\C(\gk^{++}, \gk^{+3}) \simeq 
	\linebreak[0]
	\prod_{\gx < \gk^{+3}} \C(\gk^{++}, [\ga_\gx, \ga_{\gx+1}))_{\MSE[G_\ugk]}.
\end{align*}
\end{proof}
\begin{proposition}
$\C(\gk^{+3}, \gk^{+3})_{V^*[G_\ugk]} \simeq \C(\gk^{+3}, \gk^{+6})_{\MSE[G_\ugk]}$.
\end{proposition}
\begin{proof}
We carry the proof inside $V^*[G_\ugk]$.
The first point to observe is that 
$\C(\gk^{+3}, \gk^{+3}) = \C(\gk^{+3}, \gk^{+3})_{\MSE[G_\ugk]}$.
The same technique as for the proof of \ref{Same:2:5} is to be used.
In $3$ steps we can show that
$\C(\gk^{+3}, \gk^{+3}) \simeq \C(\gk^{+3}, \gk^{+4})_{\MSE[G_\ugk]}$,
$\C(\gk^{+3}, \gk^{+3}) \simeq \C(\gk^{+3}, \gk^{+5})_{\MSE[G_\ugk]}$, 
$\C(\gk^{+3}, \gk^{+3}) \simeq \C(\gk^{+3}, \gk^{+6})_{\MSE[G_\ugk]}$.
\end{proof}
\begin{corollary}
 $(\C(\gk^{+},  \gk^{+4}) \times
 	\C(\gk^{++}, \gk^{+5}) \times
	\C(\gk^{+3}, \gk^{+6}))_{\MSE[G_\ugk]}$
is isomorphic to
$\C(\gk^{+},  \gk^{+3}) \times
	\C(\gk^{++}, \gk^{+3}) \times \linebreak[0]
	\C(\gk^{+3}, \gk^{+3})_{V^*[G_\ugk]}$.
\end{corollary}
\begin{corollary} \label{EmbeddedForcing}
There is $\gj \in V^*[G_\ugk]$ such that
$\gj \func (\C(\gk^{+},  \gk^{+4}) \times
	\C(\gk^{++}, \gk^{+5}) \times
	\C(\gk^{+3}, \gk^{+6}))_{\MSE[G_\ugk]} \to
\C(\gk^{+},  \gk^{+4}) \times
	\C(\gk^{++}, \gk^{+5}) \times
	\C(\gk^{+3}, \gk^{+6})_{V^*[G_\ugk]}$ is a complete embedding.
\end{corollary}
Choosing $H$,
a $\GN{Q}_\gk[G_{\upto \gk}]$-generic filter over $\VS[G_{\upto \gk}]$,
 is done
with some caution. For this we 
set
	$S_1 = \C(\gk^+, \gk^{+4})_{\VS[G_{\upto \gk}]}$,
	$S_2 = (\C(\gk^{++}, \gk^{+5}) \times \C(\gk^{+3}, \gk^{+6}))
                                _{\VS[G_{\upto \gk}]}$.
That is $\GN{Q}_\gk[G_{\upto \gk}] = S_1 \times S_2$.
In the same fashion we
factor $\GN{Q}^\Es_\gk[G_{\upto \gk}]$:
$S^\Es_1 = \C(\gk^+, \gk^{+4})_{\MES[G_{\upto \gk}]}$,
$S^\Es_2 = (\C(\gk^{++}, \gk^{+5}) \times \C(\gk^{+3}, \gk^{+6}))
                                _{\MES[G_{\upto \gk}]}$.
We point out that $\gj$ embeds $S^\Es_1$, $S^\Es_2$ into  $S_1$, $S_2$.
We set
	$S^U_1 = \C(\gk^+, \gk^{+4})_{\NS[G_{\upto \gk}]}$,
	$S^U_2 = (\C(\gk^{++}, \gk^{+5}) \times \C(\gk^{+3}, \gk^{+6}))
                                _{\NS[G_{\upto \gk}]}$.
Of course,
	$i_{U, E}(S^U_1 \times S^U_2) = S^\Es_1 \times S^\Es_2$.
$\NS[G_{\upto \gk}]$ contains only $\gk^+$ anti-chains of $S^U_2$ and
$S^U_2$ is a $\gk^+$-closed forcing notion in $\VS[G_{\upto \gk}]$.
Hence, there is a decreasing sequence $\ordof{p_\gx}{\gx < \gk^+} 
        \subset \NS[G_{\upto \gk}]$ that
give rise to an $S^U_2$-generic filter over $\NS[G_{\upto \gk}]$.
The crucial point is that $\ordof{i_{U, \Es}(p_\gx)}{\gx < \gk^+} \subset
	S^\Es_2$ 
and $\ordof{\gj(i_{U, \Es}(p_\gx))}{\gx < \gk^+} \subset
	S_2$.
As $\VS[G_\ugk] \satisfies \formula {S_2 \text{ is } \gk^{++}\text{-closed} }$
there is
$p \in \GN{Q}_\gk[G_\ugk]$ such that $\forall \gx < \gk^+$ $p \leq 
\gj(i_{U, \Es} (p_\gx))$.

We take $H_2$ to be $S_2$-generic over $\VS[G_{\upto \gk}]$ with $p \in H_2$. 
Let $H_1$ be $S_1$-generic over $\VS[G_{\upto \gk}][H_2]$.
We set $H = H_1 \times H_2$. Then $H$ is $\GN{Q}_\gk[G_{\upto \gk}]$-generic
over $\VS[G_{\upto \gk}]$.

We find, in $\VS[G_{\upto \gk}][H]$, filters $G^\gt$, $H^\gt$,
$G^\Es$, $H^\Es$ such that
\begin{enumerate}
\item
	$G^\gt \times H^\gt$ is $P^\gt_{j_\gt(\gk)} * 
		\GN{Q}^\gt_{j_\gt(\gk)}$-generic 
	over $\MSt$,
\item
	$j''_\gt(G_{\upto \gk} * H) \subseteq
		 G^\gt * H^\gt$,
\item
	$G^\Es \times H^\Es$ is $P^\Es_{j_\Es(\gk)} * \GN{Q}^\Es_{j_\Es(\gk)}$-generic 
	over $\MSE$,
\item
	$j''_\Es(G_{\upto \gk} * H) \subseteq
		 G^\Es * H^\Es$,
\item
        $\forall \gt' < \gt$ $i''_{\gt', \gt} (G^{\gt'} \times H^{\gt'})
                \subseteq G^\gt \times H^\gt$,
\item 
        $i''_{\gt, \Es} (G^{\gt} \times H^{\gt})
                \subseteq G^\Es \times H^\Es$.
\end{enumerate}
Under these conditions the lifting of $i_{\gt', \gt}$, $i_{\gt, \Es}$, 
$j_\gt$, $j_\Es$ is possible.
Moreover, the lifting $j_\gt$ is generated by an extender that 
continues $E(\gt)$. 

There is $G^\Es_\gk \in V^*[G_\ugk][H]$
which is $\GN{Q}^\Es_\gk[G_{\upto \gk}]$-generic over
	$\MES[G_\ugk]$, namely the filter generated by
	$\gj^{-1 \prime \prime} H$.
For later purpose we do the factoring $G^\Es_\gk = H^\Es_1 \times H^\Es_2$
such that $H^\Es_1 \subseteq H_1$, $H^\Es_2 \subseteq H_2$.
Note that $\ordof{i_{U, \Es}(p_\gx)}{\gx < \gk^+} \subset H^\Es_2$.

Let $G^\gt_\gk$ be the filter generated by 
$i^{-1 \prime \prime}_{\gt, \Es} G^\Es_\gk$.
We claim that it is 	$\GN{Q}^\gt_\gk[G_{\upto \gk}]$-generic over
	$\MtS[G_\ugk]$.
So, let $A \in \MtS[G_\ugk]$ be maximal anti-chain in
		$\GN{Q}^\gt_\gk[G_{\upto \gk}]$.

Then $i_{\gt, \Es} (A) \in \MES[G_\ugk]$ is a maximal anti-chain
in $\GN{Q}^\Es_\gk[G_{\upto \gk}]$.
As $\crit (i_{\gt, \Es}) = (\gk^{+4})_{\MtS}$ we have
	$i_{\gt, \Es} (A) = i_{\gt, \Es}'' A$.
Hence there is $a \in A$ such that
	$i_{\gt, \Es}(a) \in i_{\gt, \Es}(A) \intersect G^\Es_\gk$.
By its definition $a \in A \intersect G^\gt_\gk$.

Let $G^U_\gk$ be the filter generated by 
$i^{-1 \prime\prime}_{U, \Es} G^\Es_\gk$.
We claim that $G^U_\gk$ is $\GN{Q}^U_\gk[G_{\upto \gk}]$-generic over 
$\NS[G_{\upto \gk}]$.
For the purpose of the proof we let $G^U_\gk = H^U_1 \times H^U_2$
so that $i''_{U, \Es} H^U_1 \subseteq H^\Es_1$, 
                $i''_{U, \Es} H^U_2 \subseteq H^\Es_2$.
We  start by showing that $H^U_2$ is $Q^U_2$-generic over
$\NS[G_{\upto \gk}]$. 

As $\setof{i_{U, \Es}(p_\gx)} {\gx < \gk^+} \subseteq 
H^\Es_2$ we get $\setof{p_\gx} {\gx < \gk^+} \subseteq H^U_2$. The sequence
$\setof{p_\gx} {\gx < \gk^+}$ generates a $Q^U_2$-generic filter over
$\NS[G_\ugk]$. Hence $H^U_2$ is $Q^U_2$-generic over $\NS[G_\ugk]$.

We  show that $H^U_1$ is $Q^U_1$-generic over
$\NS[G_{\upto \gk}][H^U_2]$. As $\NS[G_\ugk] \satisfies \formula{
Q^U_2 \text{ is } \linebreak[0] \gk^{++} \text{-closed and }Q^U_1 
	\text{ is }\gk^{++} \text{c.c.}}$ it is 
enough to show
genericity over $\NS[G_{\upto \gk}]$.
So, let $A \in \NS[G_{\upto \gk}]$ be a maximal anti-chain in 
$Q^U_1$.

Then $i_{U, \Es} (A) \in \MS[{G_\ugk}]$ is a maximal anti-chain in
$Q_1$. 
This time we get $i_{U, \Es} (A) = i''_{U, \Es} A$ for free.
So there is $i_{U, \Es} (a) \in H^\Es_1 \intersect i_{U, \Es} (A)$ and
by its definition $a \in H^U_1$. Hence $H^U_1 \intersect A \not=
\emptyset$. With this we showed that
$G^U_\gk = H^U_1 \times H^U_2$ is $\GN{Q}^U_\gk[G_{\upto \gk}]$-generic over 
        $\NS[G_{\upto \gk}]$.

Of course $i''_{U, \gt} G^U_\gk \subseteq G^\gt_\gk$, hence we have the 
lifting
\begin{align*}
\begin{diagram}
\node{} 
        \node{\MSE[G_{\upto \gk}][G^\Es_\gk]}
\\
\node{\NS[G_{\upto \gk}][G^U_\gk]}
		 \arrow{e,b}{i_{U, \gt}}
		 \arrow{ne,b}{i_{U, \Es}}
        \node{\MSt[G_{\upto \gk}][G^\gt_\gk]}
		 \arrow{n,b}{i_{\gt, \Es}}
\end{diagram}
\end{align*}
\ifnum\article=0
\enlargethispage*{10pt}
\fi
We set
\begin{multline*}
R_U = (\Col(\gk^{+6}, i_U(\gk)) \times 
	\C(\gk^{+4}, i_U(\gk)^+) \times
		\C(\gk^{+5}, i_U(\gk)^{++}) \times
	\\
	\C(\gk^{+6}, i_U(\gk)^{+3}))
				_{\NS[G_{\upto \gk}]}.
\end{multline*}
Let $i^2_U$ be the iterate of $i_U$. We choose a function, $R(-,-)$, such that
\begin{align*}
& R_U = i^2_U(R)(\gk, i_U(\gk)).
\end{align*}
We are going to find  $I_U$,
an $R_U $-generic 
filter, over
$\NS[G_{\upto \gk}]$. We do not force with $I_U$.
Anticipating its later usage in the definition of $\PE$
 we need it to be in $\MSE[G_\ugk][G^\Es_\gk]$. For this we work as follows.

As $U \in \MSE$ and $\forall \gt < \len(\Es)$ $E(\gt) \in \MSE$
we have the following diagram
\begin{align*}
\begin{aligned}
\begin{diagram}
\node{\MSE} 
		\arrow{s,l}{i_U^\Es}
		\arrow{se,t}{j_{\gt}^\Es}
\\
\node{\NSE}
		\arrow{e,b}{i_{U, \gt}^\Es}
	\node{\MStE}
\end{diagram}
\end{aligned}
\begin{aligned}
\begin{split}
& U = E_\gk(0),
\\
& i^\Es_U \func  \MSE \to \NSE \simeq \Ult(\MSE, U),
\\
& j^\Es_\gt \func  \MSE \to \MStE \simeq \Ult(\MSE, E(\gt)),
\\
& i^\Es_{U, \gt}(i^\Es_U(f)(\gk)) = j^\Es_{\gt}(f)(\gk),
\end{split}
\end{aligned}
\end{align*}
As $V^{V^*}_{\gk+3} = V^\MSE_{\gk+3}$ we get that the following
are
generic extensions:
$\MSE[G_\ugk][G^\Es_\gk]$,
$\NSE[G_\ugk][G^N_\gk]$.
We set
\begin{multline*}
R^\Es_U = (\Col(\gk^{+6}, i^\Es_U(\gk)) \times 
	\C(\gk^{+4}, i^\Es_U(\gk)^+) \times
		\C(\gk^{+5}, i^\Es_U(\gk)^{++}) \times
	\\
	\C(\gk^{+6}, i^\Es_U(\gk)^{+3}))
				_{\NSE[G_{\upto \gk}]}.
\end{multline*}
$R_U$, $R^\Es_U$  and their
anti-chains are coded in
$V^{\NS[G_\ugk]}_{i_U(\gk)+3}$,
$V^{\NSE[G_\ugk]}_{i^\Es_U(\gk)+3}$ respectively.
$V^{\NS[G_\ugk]}_{i_U(\gk)+3}$, 
$V^{\NSE[G_\ugk]}_{i^\Es_U(\gk)+3}$
are determined by $V^{V^*}_{\gk+3}$, $V^{\MSE}_{\gk+3}$ (and $U$, of course).
As $U \in \MSE$ and $V^{V^*}_{\gk+3} = V^\MSE_{\gk+3}$ we get that
	$R_U =R^\Es_U$
and each anti-chain of $R_U$
appearing in $\NS[G_\ugk]$
is also an anti-chain of 
$R^\Es_U$
 appearing in $\NSE[G_\ugk]$.

Hence, if $I_U$ is 
$R^\Es_U$-generic filter over 
$\NSE[G_\ugk]$ then it is also 
$R_U$-generic 
filter over $\NS[G_\ugk]$.
Construction of such $I_U$ is done as follows. 
We set
\begin{align*}
& R^U_1 = \Col(\gk^{+6}, i_U^\Es(\gk))_{\NSE[G_{\upto \gk}]},
\\
& R^U_2 = \C(\gk^{+6}, i_U^\Es(\gk)^{+3})_{\NSE[G_{\upto \gk}]},
\\
& R^U_3 = \C(\gk^{+5}, i_U^\Es(\gk)^{++})_{\NSE[G_{\upto \gk}]},
\\
& R^U_4 = \C(\gk^{+4}, i_U^\Es(\gk)^+)_{\NSE[G_{\upto \gk}]},
\end{align*}
so that $R^U = R^U_1 \times R^U_4 \times R^U_3 \times
		R^U_2$.

We claim that there is $I^U_1 \in \MSE[G_{\upto \gk}]$ 
which is $R^U_1$-generic 
over $\NSE[G_{\upto \gk}]$. 
This is immediate due to $R^U_1$
being $\gk^+$-closed in $\MSE[G_{\upto \gk}]$
and $\NSE[G_{\upto \gk}]$ containing only $\gk^+$ maximal anti-chains of 
$R^U_1$.

The next step is to show that there is $I^U_2 \in \MSE[G_{\upto \gk}][G^\Es_\gk]$ 
which is 
$R^U_2$-generic over $\NSE[G_{\upto \gk}][I^U_1]$.
In $\MSE[G_{\upto \gk}]$, $R^U_2 \simeq 
			\C(\gk^+, \gk^{+3})_{\MSE[G_{\upto \gk}]}$.
As $G^\Es_\gk$ is generic over $\MSE[G_{\upto \gk}]$ and $I^U_1 \in 
			\MSE[G_{\upto \gk}]$
we get that there is $I^U_2 \in \MSE[G_{\upto \gk}][G^\Es_\gk]$ which is
$R^U_2$-generic over $\NSE[G_{\upto \gk}][I^U_1]$.

Locating $I^U_3 \in \MSE[G_{\upto \gk}][G^\Es_\gk]$ which is
$R^U_3$-generic over $\NSE[G_{\upto \gk}][I^U_1 \times I^U_2]$ is done
as follows.
As $\NSE[G_{\upto \gk}] \satisfies \formula{ R^U_1 \times R^U_2 \text{ is } 
\gk^{+6}\text{-closed} }$, finding $I^U_3$ which is generic over 
	$\NSE[G_{\upto \gk}]$
is enough to ensure genericity over $\NSE[G_{\upto \gk}][I^U_1 \times I^U_2]$.
In $\MSE[G_{\upto \gk}]$, $R^U_3 \simeq 
			\C(\gk^+, \gk^{++})_{\MSE[G_{\upto \gk}]}$. 
Once more, as $G^\Es_\gk$ is
generic over $\MSE[G_{\upto \gk}]$, there is
$I^U_3 \in \MSE[G_{\upto \gk}][G^\Es_\gk]$ which $R^U_3$-generic over 
						$\NSE[G_{\upto \gk}]$.

$I^U_4$ is constructed in the same way.

Let us set $I_U = I^U_1 \times I^U_2 \times I^U_3 \times I^U_4$.
The above yield that $I_U$ is $R_U$-generic over $\NSE[G_{\upto \gk}]$,
hence over $\NS[G_\ugk]$.
As $\NS[G_{\upto \gk}] \satisfies \formula{ R^U \text{ is } 
						\gk^{+4} \text{-closed} }$,
$G^U_\gk$ is generic over $\NS[G_{\upto \gk}][I_U]$.
Hence $I_U$ is $R_U$-generic over $\NS[G_{\upto \gk}][G^U_\gk]$.
We need to construct $G^U_\dgk$ with some care
in order to allow $I_U$ to be generic over $\NS[G_{\upto \gk}][G^U_\gk][G^U_\dgk]$.

For this we look at $\GN{P}^U_{\dgk}[G_\ugk * G^U_\gk]$
in $\NS[G_\ugk][G^U_\gk][I^U_1][I^U_2]$.
This forcing, as seen by $\VS[G_\ugk][H]$, is $\gk^+$-closed of
size $\gk^+$. Hence it is isomorphic to $\C(\gk^+, \gk^+)_{\VS[G_\ugk]}$.
The main point is that we do not need all of $H$ in order to see 
$\NS[G_\ugk][G^U_\gk][I^U_1][I^U_2]$. Namely, we set
$\gl = (\gk^{+4})_{\MES}$ and factor
$\GN{Q}_\gk[G_\ugk]$ as
$(C(\gk^+, \gl) \times C(\gk^{++}, \gk^{+5}) \times C(\gk^{+3}, \gk^{+6}))
 \times   \C(\gk^+, \gk^{+4} \setminus \gl)$.
Then we factor $H$ appropriately as $H' \times H''$.
Then $\NS[G_\ugk][G^U_\gk][I^U_1][I^U_2]$ is definable inside
$\VS[G_\ugk][H']$.
So, $\GN{P}^U_{\dgk}[G_\ugk * G^U_\gk]$ is embeddable in 
$\C(\gk^+, \gk^{+4} \setminus \gl)$ and all of its dense sets appearing
in $\NS[G_\ugk][G^U_\gk][I^U_1][I^U_2]$ are coded in $\VS[G_\ugk][H']$.
Hence, there is $G^U_\dgk \in \VS[G_\ugk][H'][H'']$
which is $\GN{P}^U_{\dgk}[G_\ugk * G^U_\gk]$-generic over 
	$\NS[G_\ugk][G^U_\gk][I^U_1][I^U_2]$.

We consider $R^U_3 \times R^U_4$
in $\NS[G_\ugk][G^U_\gk][I^U_1][I^U_2][G^U_\dgk]$.
Evidently each anti chain of $R^U_3 \times R^U_4$
in $\NS[G_\ugk][G^U_\gk][I^U_1][I^U_2][G^U_\dgk]$ already appears
in $\NS[G_\ugk][G^U_\gk]$.
Hence $I^U_3 \times I^U_4$ is $R^U_3 \times R^U_4$-generic over
	$\NS[G_\ugk][G^U_\gk][I^U_1][I^U_2][G^U_\dgk]$.

Of course, all of this means that $I_U$ is $R_U$-generic over
	$\NS[G_\ugk][G^U_\gk][G^U_\dgk]$.

We set $G^U = G_{\upto \gk} * G^U_\gk *  G^U_{\downto \gk}$.
Then $G^U$
is $P^U_{i_U(\gk)}$-generic over $\NS$. As $i''_U G_{\upto \gk} = 
G_{\upto \gk}$ we get that
$i''_U G_{\upto \gk} \subseteq G^U$, hence we have the lifting
$i_U \func  \VS[G_{\upto \gk}] \to \NS[G^U]$. 
Note that
$i_U$ is defined in $\VS[G_\ugk][H]$ and it is the natural embedding
defined by a $\VS[G_{\ugk}]$-ultrafilter extending $U$.

Let $H^U$ be the filter generated by $i''_U H$.
This definition is possible as $H \subset \VS[G_{\upto \gk}]$ and we just 
lifted $i_U$ to $\VS[G_{\upto \gk}]$.
We claim that $H^U$ is $\GN{Q}^U_{i_U(\gk)}[G^U]$-generic
over	$\NS[G^U]$. Let, then, 
$D \in \NS[G^U]$ be dense open in 
$\GN{Q}^U_{i_U(\gk)}[G^U]$.

Then there is $f \in \VS[G_{\upto \gk}]$ such that $i_U(f)(\gk) = D$.
On a big set, in $\VS[G_{\upto \gk}][H]$ sense, $f(\gn)$ is a dense open subset
of $\GN{Q}_\gk[G_{\upto \gk}]$. $\GN{Q}_\gk[G_{\upto \gk}]$ is $\gk^+$-closed 
in $\VS[G_{\upto \gk}]$ hence
$D^* = \bigintersect_{\gn < \gk} f(\gn) \in \VS[G_{\upto \gk}]$ is dense open in it 
and by its
definition $i_U(D^*) \subseteq i_U(f)(\gk)$.
Choose $p \in D^* \intersect H$. 
$i_U(p) \in i_U(D^*) \intersect i_U''H$. Hence, 
$D \intersect H^U \not= \emptyset$.

So we can lift $i_U$ to $i_U\func \VS[G_{\upto \gk}][H] \to \NS[G^U][H^U]$.
This embedding is definable inside $\VS[G_{\upto \gk}][H]$.

We note that $I_U$ is $R_U$-generic over $\NS[G^U][H^U]$
as $\NS[G^U] \satisfies \formula {\GN{Q}^U_{i_U(\gk)}[G^U] \text{ is } 
	\linebreak[0]
	i_U(\gk)^+ 
	\linebreak[0]
	\text{-closed}}$,
hence adds to new anti-chains to $R_U$.

Let $G^\gt_{\downto \gk}$ be the filter generated by
$i''_{U, \gt} G^U_{\downto \gk}$. We claim that $G^\gt_{\downto \gk}$ is
$\GN{P}^\gt_{\downto \gk}[G_{\upto \gk} * G^\gt_\gk]$-generic over 
	$\MSt[G_{\upto \gk}][G^\gt_\gk]$.
So, let $D \in \MSt[G_{\upto \gk}][G^\gt_\gk]$ be dense open in 
$\GN{P}^\gt_{\downto \gk}[G_{\upto \gk} * G^\gt_\gk]$.

Let $\GN{D} \in \MSt$ be a $P_\gk * \GN{Q}^\gt_\gk$-name for $D$.
Then there is $f \in \VS$ such that $j_\gt(f)(\Es_\ga \restricted \gt) = 
\GN{D}$.
Hence
$\setof{\gns} {
	\forces_{P_{\gk^0(\gns)} * \GN{Q}_{\gk^0(\gns)}} 
		\formula{ f(\gns) 
			\text{ is dense open in } \GN{P}_{\downto \gk^0(\gns)} 
		\text{ which }\linebreak[0] \text{is } \linebreak[0] 
		\gk^0(\gns)^{+4} \text{-closed}}
	} \in E_\ga(\gt)$.
Hence, there is a name, $f^*(\gm)$,  such that
$\forces_{P_{\gm} * \GN{Q}_{\gm}} \formula{ f^*(\gm) \subseteq f(\gns) 
	\text{ is dense open}}$
whenever $\gk^0(\gns) = \gm$. Hence
\ifnum\article=1
\enlargethispage*{10pt}
\fi
\begin{align*}
& \MtS \satisfies \formula {
	\forces_{P_{\gk} * \GN{Q}^\gt_{\gk}} \formula{ j_\gt(f^*)(\gk) 
		\subseteq
		\GN{D} \text{ is dense open in } 
			\GN{P}^\gt_{\downto \gk} }
	},
\\
& \NS \satisfies \formula {
	\forces_{P_{\gk} * \GN{Q}^U_{\gk}} \formula{ i_U(f^*)(\gk) 
			\text{ is dense open in } \GN{P}^U_{\downto \gk} }
	}.
\end{align*}
So there is $i_U(g)(\gk) \in i_U(f^*)(\gk)[G_{\upto \gk}][G^U_\gk] \intersect 
	G^U_{\downto \gk}$.
We get $D \intersect G^\gt_{\downto \gk} \not= \emptyset$ by noticing
that
$j_\gt(g)(\gk) \in \linebreak[0] j_\gt(f^*)(\gk)[G_{\upto \gk}]
		\linebreak[0]
		[G^\gt_\gk] 
		\intersect 
		\linebreak[0]
		i''_{U, \gt}G^U_{\downto \gk}$.

Let $G^\gt = G_{\upto \gk} * G^\gt_\gk * G^\gt_{\downto \gk}$.
$j''_\gt G_{\upto \gk} = G_{\upto \gk}$, so $j''_\gt G_{\upto \gk} \subseteq G^\gt$. Hence we can lift
$j_\gt$ to
	$j_\gt\func  \VS[G_{\upto \gk}] \to \MSt[G^\gt]$.
We note that this lift is defined in $\VS[G_{\upto \gk}][H]$ and it is the natural
embedding of a $\VS[G_{\upto \gk}]$-extender continuing $E(\gt)$.

Let $H^\gt$ be the filter generated by $j_\gt''H$.
We claim that $H^\gt$ is $\GN{Q}^\gt_{j_\gt(\gk)}[G^\gt]$-generic
over $\MSt[G^\gt]$. Let $D \in \MSt[G^\gt]$ be dense open in
$\GN{Q}^\gt_{j_\gt(\gk)}[G^\gt]$.

Then there is $f \in \VS[G_{\upto \gk}]$ such that 
$j_\gt(f)(\Es_\ga \restricted \gt) = D$.
As
$A = \setof {\gns}
	{
	 f(\gns) \text{ is} \linebreak[0] \text{ dense} \linebreak[0]
			\text{ open} \linebreak[0]
			\text{ in }  \linebreak[0]
			\GN{Q}_\gk[G_{\upto \gk}]
	} \in E_\ga(\gt)$
and $\GN{Q}_\gk[G_{\upto \gk}]$ is $\gk^+$-closed we get that
$D^* = \bigintersect_{\gns \in A} f(\gns) \in \VS[G_{\upto \gk}]$ is dense open in
$\GN{Q}_\gk[G_{\upto \gk}]$. So there is $p \in D^* \intersect H$.
Hence $j_\gt(p) \in D^* \intersect j_\gt''H$.
Yielding, $D \intersect H^\gt \not= \emptyset$.

So we can do the lift $j_\gt \func  \VS[G_{\upto \gk}][H] \to \MSt[G^\gt][H^\gt]$.
In order to finish we need to build generic filters over $\MES$.
We split the handling into $2$ cases:
\begin{enumerate}
\item 
	$\len(\Es) = \gt + 1$:
	In this case we have $\MES = \MtS$, so by setting
	$G^\Es = G^\gt$, $H^\Es = H^\gt$ we have the needed filters.
\item 
	$\len(\Es)$ is limit:
	We let $G^\Es$, $H^\Es$ be the filters generated by
	$\bigunion_{\gt < \len(\Es)} i''_{\gt, \Es} G^\gt$,
	$\bigunion_{\gt < \len(\Es)} i''_{\gt, \Es} H^\gt$
	respectively.
\end{enumerate}
After the forcing the generic extension the power function we have is
\begin{align*}
2^\gm = 
	\begin{cases}
	\gm^{+3}	&	\text{if }\gm \in \set{\gn^+, \gn^{++}, \gn^{+3}}
				\text{ where }\gn \leq \gk \text{ is inaccessible}
	\\
	\gm^+		&	\text{Otherwise}
	\end{cases},
\end{align*}
and we still have \eqref{E-system}.

The new diagram we have after all the liftings is
\begin{align*}
\begin{diagram}
\node{V = \VS[G_{\upto \gk}][H]} 
		\arrow[2]{e,t}{j_\Es}
		\arrow{s,l}{i_{U}}
		\arrow{se,b}{j_{\gt'}}
		\arrow{see,b}{j_\gt}
	\node{}
        \node{\ME = \MSE[G^\Es][H^\Es]}
\\
	\node{N = \NS[G^U][H^U]}
		 \arrow{e,b}{i_{U, \gt'}} 
	\node{M_{\gt'} = M^*_{\gt'}[G^{\gt'}][H^{\gt'}]}
		 \arrow{ne,t,3}{i_{\gt', \Es}} 
		 \arrow{e,b}{i_{\gt', \gt}} 
	\node{\Mt = \MSt[G^\gt][H^\gt]} 
		\arrow[1]{n,b}{i_{\gt, \Es}} 
\end{diagram}
\end{align*}
\subsection{Cardinal structure in $N[I_U]$}
We claim that in $N[I_U]$ there are
no cardinals in $[\gk^{+7}, i_U(\gk)]$ and all other $N$-cardinals
are preserved. The power function differs from the power
function of $N$ at the following points:
	$2^{\gk^{+4}} = i_U(\gk)^+$,
	$2^{\gk^{+5}} = i_U(\gk)^{++}$,
	$2^{\gk^{+6}} = i_U(\gk)^{+3}$.

Before continuing we recall that a forcing notion $P$ is called
$\gl$-dense if the intersection of less than $\gl$ dense open subsets
of $P$ is dense open. We use the following obvious proposition.
\begin{proposition} \label{DenseFromClose}
Let $P, Q$ be forcing notions being $\gl$-closed, $\gl$-c.c. respectively.
Then
	$\forces_{Q} \formula{ P \text{ is } \gl \text{-dense}}$.
\end{proposition}
Our construction of $N$ from $N^*$ means that
$N[I_U] = N^*[G_\ugk][G^U_\gk][G^U_\dgk][I_U]$. 
If we write $N^*[G_\ugk][I_U][G^U_\gk]$ 
as $N^*[G_\ugk][I^1_U \times I^2_U][I^3_U][I^4_U][G^U_\gk]$ 
then the usual arguments for product forcing
show that the claim is satisfied in this model.
However $N^*[G_\ugk][I^1_U \times I^2_U][I^3_U][I^4_U][G^U_\gk]
	\satisfies \formula{
		\GN{P}^U_\dgk[G_\ugk][G^U_\gk] \text{ is }
				\gk^{+4} \text{-closed}
	}$.
Hence we might loose $\gk^{+5}$, $\gk^{+6}$ in
$N^*[G_\ugk][I^1_U \times I^2_U][I^3_U][I^4_U][G^U_\gk][G^U_\dgk]$.
In order to show that $\gk^{+5}$ is preserved we reorder as follows
$N[I_U] = N^*[G_\ugk][I^1_U \times I^2_U][I^3_U][G^U_\gk][I^4_U][G^U_\dgk]$.

We start by showing
$N^*[G_\ugk][I^1_U \times I^2_U][I^3_U][G^U_\gk] \satisfies
	\formula{
		R^4_U \text{ is } \gk^{+5} \text{-c.c.}
	}$.
Towards a contradiction, let us assume that there is an  anti-chain
of length $\gk^{+5}$.
Until further notice we work inside 
	$N^*[G_\ugk][I^1_U \times I^2_U][I^3_U]$.
Let $\ordof{\GN{a}_\gx} {\gx < \gk^{+5}} $ be a $\GN{Q}^U_\gk[G_\ugk]$-name of 
this anti-chain.
Let $\GN{d}_\gx$ be the name of $\dom \GN{a}_\gx$.
As  $\GN{Q}^U_\gk[G_\ugk]$ is $\gk^{+4}$-c.c. we get that there are
$d_\gx$
such that
$\power{d_\gx} < \gk^{+4}$
and
$\forces_{\GN{Q}^U_\gk[G_\ugk]} \formula {\GN{d}_\gx \subseteq \VN{d}_\gx }$. 
By invoking the $\gD$-lemma 
we can find a subset of $\ordof{d_\gx} {\gx < \gk^{+5}}$ of size
$\gk^{+5}$ which is a $\gD$-system. 
Hence, without loss of generality,
we can assume that 
	$\forces_{\GN{Q}^U_\gk[G_\ugk]} \formula 
			{\dom \GN{a}_\gx = \VN{d}_\gx }$
and $\forall \gx_1 \not= \gx_2\ d_{\gx_1} \intersect d_{\gx_2} = d$.
We can have this situation only if $2^{\gk^{+3}} \geq \gk^{+5}$, 
contradicting $2^{\gk^{+3}} = \gk^{+4}$.
At this point we stop working inside $N^*[G_\ugk][I^1_U 
					\times I^2_U][I^3_U]$.
We know now that
\begin{align*}
& N^*[G_\ugk][I^1_U \times I^2_U][I^3_U][G^U_\gk] \satisfies \formula{
	R^4_U \text{ is } \gk^{+5}\text{-c.c.}
},
\\
& N^*[G_\ugk][I^1_U \times I^2_U][I^3_U][G^U_\gk] \satisfies \formula{
	\GN{P}^U_\dgk[G_\ugk][G^U_\gk] \text{ is } \gk^{+5} \text{-closed}
}.
\end{align*}
By \ref{DenseFromClose} we get
$N^*[G_\ugk][I^1_U \times I^2_U][I^3_U][G^U_\gk][I^4_U] \satisfies \formula{
	\GN{P}^U_\dgk[G_\ugk][G^U_\gk] \text{ is } \gk^{+5} \text{-dense}
}$,
hence $\gk^{+5}$ remains a cardinal in 
	$N^*[G_\ugk][I^1_U \times I^2_U][I^3_U][G^U_\gk][I^4_U][G^U_\dgk]$.

For the preservation of $\gk^{+6}$ we start with
the reorder
$N[I_U] = N^*[G_\ugk][I^1_U \times I^2_U][G^U_\gk][I^3_U][I^4_U][G^U_\dgk]$.
By the same method as before we get
\begin{align*}
& N^*[G_\ugk][I^1_U \times I^2_U][G^U_\gk] \satisfies \formula{
	R^3_U \times R^4_U \text{ is } \gk^{+6}\text{-c.c.}
},
\\
& N^*[G_\ugk][I^1_U \times I^2_U][G^U_\gk] \satisfies \formula{
	\GN{P}^U_\dgk[G_\ugk][G^U_\gk] \text{ is } \gk^{+6} \text{-closed}
}.
\end{align*}
By \ref{DenseFromClose} we get
$N^*[G_\ugk][I^1_U \times I^2_U][G^U_\gk][I^3_U][I^4_U] \satisfies \formula{
	\GN{P}^U_\dgk[G_\ugk][G^U_\gk] \text{ is } \gk^{+6} \text{-dense}
}$,
hence $\gk^{+6}$ remains a cardinal in 
	$N^*[G_\ugk][I^1_U \times I^2_U][G^U_\gk][I^3_U][I^4_U][G^U_\dgk]$.

By observing that $N^*[G_\ugk][I_U][G^U_\gk] \satisfies
	\formula{
		\power{\GN{P}^U_\dgk[G_\ugk][G^U_\gk]} = \gk^{+6}
	}$ we see that the power function of
$N^*[G_\ugk][I_U][G^U_\gk][G^U_\dgk]$ is the same as the power function
of $N^*[G_\ugk][I_U][G^U_\gk]$.
\subsection{Locating the needed generic filters}
As $U \in \ME$ and $\forall \gt < \len(\Es)$ $E(\gt) \in \ME$
we have the following diagram
\begin{align*}
\begin{aligned}
\begin{diagram}
\node{\ME} 
		\arrow{s,l}{i_U^\Es}
		\arrow{se,t}{j_\gt^\Es}
\\
\node{\NE}
		 \arrow{e,b}{i_{U, \gt}^\Es}
        \node{\MtE}
\end{diagram}
\end{aligned}
\begin{aligned}
\begin{split}
& U = E_\gk(0),
\\
& i^\Es_U \func  \ME \to \NE \simeq \Ult(\ME, U),
\\
& j^\Es_\gt \func  \ME \to \MtE \simeq \Ult(\ME, E(\gt)),
\\
& i^\Es_{U, \gt}(i^\Es_U(f)(\gk)) = j^\Es_\gt(f)(\gk).
\end{split}
\end{aligned}
\end{align*}
Recall that we have $I_U \in \ME$ which is $R_U$-generic over $\NE$.
\subsubsection{Generic over $\Mt$ when $\gt + 1 < \len(\Es)$} \label{MtGeneric}
Consider the following forcing notion:
\begin{multline*}
 R_\gt = (\Col(\gk^{+6}, j_\gt(\gk)) \times 
	 \C(\gk^{+4}, j_\gt(\gk)^+) \times
		\C(\gk^{+5}, j_\gt(\gk)^{++}) \times
	\\
	 \C(\gk^{+6}, j_\gt(\gk)^{+3}))
				_{\MSt[G_{\upto \gk}]}.
\end{multline*}
We show that there is $I_\gt \in \ME$, an $R_\gt$-generic
filter over $\Mt$. 
Moreover, whenever $\gt' < \gt < \len(\Es)$ we have
	$i''_{\gt', \gt}I_{\gt'} \subseteq I_\gt$.
For this we set
\begin{multline*}
 R^\Es_\gt = (\Col(\gk^{+6}, j^\Es_\gt(\gk)) \times 
	 \C(\gk^{+4}, j^\Es_\gt(\gk)^+) \times
		\C(\gk^{+5}, j^\Es_\gt(\gk)^{++}) \times
	\\
	 \C(\gk^{+6}, j^\Es_\gt(\gk)^{+3}))
				_{\MStE[G_{\upto \gk}]}.
\end{multline*}
$R_\gt$, $R^\Es_\gt$  are coded in $V^{\MSt[G_\ugk]}_{j_\gt(\gk)+3}$, 
	$V^{\MStE[G_\ugk]}_{j^\Es_\gt(\gk)+3}$ respectively.
$V^{\MSt[G_\ugk]}_{j_\gt(\gk)+3}$, $V^{\MStE[G_\ugk]}_{j^\Es_\gt(\gk)+3}$
are determined by $V^{V^*}_{\gk+3}$, $V^{\MSE}_{\gk+3}$ (and $E(\gt)$, 
of course).
As $E(\gt) \in \MSE$ and $V^{V^*}_{\gk+3} = V^\MSE_{\gk+3}$ we get that
$R_\gt = R^\Es_\gt$.

By the same reasoning, each anti-chain of $R_\gt$ appearing in 
$\Mt$
is also an anti-chain of $R^\Es_\gt$ appearing in $\MtE$.
Hence, if $I_\gt \in \ME$ is 
$R^\Es_\gt$-generic filter over $\MtE$
then it is also $R_\gt$-generic filter over $\Mt$.

Let $I_\gt \in \ME$ be the filter generated by 
	$i^{\Es \prime \prime}_{U, \gt} I_U$.
We have the natural factoring $I_\gt = I^\gt_1 \times \dotsb \times I^\gt_4$,
We claim that $I_\gt$ is $R_\gt$-generic
over $\MtE$. We start by showing genericity over $\MStE[G_\ugk]$.
So, let $D \in \MStE[G_\ugk]$ be dense open in $R_\gt$.

Let $\GN{D} \in \MStE$ be a $P_\gk$-name for $D$.
Choose $f \in \MSE$ such that $j^\Es_\gt(f)(\Es_\ga \restricted \gt) = \GN{D}$.
So in $\MSE$ we have
\begin{multline*}
 A = \setof {\gns}
     {
        \MSE \satisfies \formula {
	  \forces_{P_{\gk^0(\gns)}} \formula {
                 f(\gns) \text{ is dense open in }
		 \Col(\gk^0(\gns)^{+6}, \gk)
				\times
		\\
				\C(\gk^0(\gns)^{+4}, \gk^{+})
				\times
				\C(\gk^0(\gns)^{+5}, \gk^{++})
				\times
				\C(\gk^0(\gns)^{+6}, \gk^{+3})
           }
	}
     }
	\in E_\ga(\gt).
\end{multline*}
The standard observation 
        $B = \setof {\gm}
	{
	  \power{\setof{\gns \in A} {\gk^0(\gns) = \gm} } \leq \gm^{+3}
	} \in E_\gk(0)$
yields that for each $\gm \in B$ there is a $P_\gm$-name, $f^*(\gm)$, such that
for all $\gns \in A$ with $\gk^0(\gns) = \gm$ we have
$\MSE \satisfies \formulal{} \forces_{P_\gm} \formula{ f^*(\gm) \subseteq f(\gns) 
				\text { is dense open}} \formular{}$.
Hence
\begin{align*}
& \NSE \satisfies \formula{
	\forces_{P_\gk} \formula{ i^\Es_U(f^*)(\gk) \text { is dense open}}
	 },
\\
& \MStE \satisfies \formula{
	\forces_{P_\gk} \formula{ j^\Es_\gt(f^*)(\gk) 
				\subseteq j^\Es_\gt(f)(\Es_\ga \restricted \gt)
                }
        }.
\end{align*}
So there is $g \in \MSE$ such that
$i^\Es_U(g)(\gk)[G_{\upto \gk}] \in i^\Es_U(f^*)(\gk)[G_{\upto \gk}] 
		\intersect I_U$.
Noting that we have $i^\Es_{U, \gt}\func \NSE[G_{\upto \gk}] \to 
		\MStE[G_{\upto \gk}]$, we get
$j^\Es_\gt(g)(\gk)[G_{\upto \gk}] \in j^\Es_\gt(f^*)(\gk)[G_{\upto \gk}] 
		\intersect i^{\Es\prime \prime}_{U, \gt}I_U$.
That is $j^\Es_\gt(f)(\Es_\ga \restricted \gt)[G_{\upto \gk}] \intersect  
		I_\gt \not= 
	\emptyset$.
By this we showed that $I_\gt$ is $R_\gt$-generic over
$\MStE[G_\ugk]$. Hence it is $R_\gt$-generic over
$\MSt[G_\ugk]$.

So, we can consider $\MSt[G_{\upto \gk}][I_\gt]$.
As $\MSt \satisfies \formula{ R_\gt \text{ is } \gk^{+4}\text{-closed}}$ we
get that $G^\gt_\gk$ is $\GN{Q}^\gt_\gk[G_{\upto \gk}]$-generic
over $\MSt[G_{\upto \gk}][I_\gt]$. By commutativity of product forcing
we get that $I_\gt$ is $R_\gt$-generic over $\MSt[G_{\upto \gk}][G^\gt_\gk]$.
Hence we can consider $\MSt[G_{\upto \gk}][G^\gt_\gk][I_\gt]$.
We want to show that $G^\gt_\dgk$ is generic over
	$\MSt[G_{\upto \gk}][G^\gt_\gk][I_\gt]$.
For this we lift
$i_{U, \gt}$ to 
	$i^*_{U, \gt}\func \NS[G_{\upto \gk}][G^U_\gk][I_U] \to 
		\MSt[G_{\upto \gk}][G^\gt_\gk][I_\gt]$
which is possible by recalling how we generated $I_\gt$ from $I_U$.
Let $D \in \MSt[G_{\upto \gk}][G^\gt_\gk][I_\gt]$ be dense open in
$\GN{P}^\gt_\dgk[G_\ugk][G^\gt_\gk]$.

Then there is $\GN{D} \in \MSt[G_{\ugk}][G^\gt_\gk]$ 
such that $\GN{D}[I_\gt] = D$.
Of course $\GN{D} \in \MSt[G^\gt]$ as
$\MSt[G^\gt] = \MSt[G_{\ugk}][G^\gt_\gk][G^\gt_\dgk] \supset 
		\MSt[G_{\ugk}][G^\gt_\gk]$.
Hence there is
$f \in \VS[G_\ugk]$ such that $j_\Es(f)(\Es_\ga \restricted \gt) = \GN{D}$.
Then in $\VS[G_\ugk][H]$ we have
\begin{multline*}
A = \setof {\gns}
	{\VS[G_{\upto \gk^0(\gns)}][G_{\gk^0(\gns)}]
		\satisfies
		\formula{
		\forces_{R(\gk^0(\gns), \gk)} 
\\
		\formula{
		f(\gns) \text{ is dense open in }
		\GN{P}_{\downto \gk^0(\gns)}
		}
	}}
	\in \Es_\ga(\gt).
\end{multline*}
We note that  for $\gm$ inaccessible we have
	$f \restricted \gm^{+3} \in 
	\VS[G_{\upto \gm}][G_{\gm}]$
and
$\VS[G_{\upto \gm}][G_{\gm}] \satisfies \formula{
	\forces_{R(\gm, \gk)} \formula{
	\GN{P}_{\downto \gm}[G_{\upto \gm}][G_{\gm}] 
	 \text{ is }\gm^{+4} \text{-closed}
}}$. 
Hence there is $f^* \in \VS[G_\ugk]$ such that
$\forall \gns \in A$
$\VS[G_{\upto \gk^0(\gns)}][G_{\gk^0(\gns)}] \satisfies \formula{
	\forces_{R(\gk^0(\gns), \gk)} \formula{
	f^*(\gk^0(\gns)) \subseteq f(\gns)
		\text{ is dense open}
}}$. This implies
\begin{align*}
& \NS[G_{\ugk}][G^U_{\gk}][I_U] \satisfies \formula{
	i_U(f^*)(\gk)[I_U] \text{ is dense open}
	},
\\
& \MtS[G_{\ugk}][G^\gt_{\gk}][I_\gt] \satisfies \formula{
	j_\gt(f^*)(\gk)[I_\gt] \subseteq 
		j_\gt(f)(\Es_\ga \restricted \gt)[I_\gt]
	}. 
\end{align*}
By genericity of $G^U_\dgk$ over $\NS[G_\ugk][G^U_\gk][I_U]$
we get that  there is $g \in \VS[G_\ugk][H]$ such that $i_U(g)(\gk) \in 
	i_U(f^*)(\gk)[I_U] \intersect G^U_\dgk$.
Hence  $j_\gt(g)(\gk)\in 
	j_\gt(f)(\Es_\ga \restricted \gt)[I_\gt] \intersect 
				i''_{U, \gt }G^U_\dgk$.
That is $D \intersect G^\gt_\dgk \not= \emptyset$.

Of course this implies that $I_\gt$ is generic over $\MtS[G^\gt]$.
As $H^\gt$ adds no new anti-chains to $R_\gt$ we get that
$I_\gt$ is generic over $\MtS[G^\gt][H^\gt]$.

As for the moreover: 
	$V^{\NSE}_{i^\Es(\gk)+3} = V^\NS_{i(\gk)+3}$
and
                $i^\Es_{U, \gt} \restricted 
				V^\NSE_{i^\Es(\gk)+3} =
                	i_{U, \gt} \restricted 
				V^\NS_{i(\gk)+3}$,
                $i^\Es_{\gt', \gt} \restricted 
				V^\NSE_{i^\Es(\gk)+3} =
                	i_{\gt', \gt} \restricted 
				V^\NS_{i(\gk)+3}$.
By their definition
       $i^{\Es \prime \prime}_{\gt', \gt} \circ 
		i^{\Es \prime \prime}_{U, \gt'} I^U =
                i^{\Es \prime \prime}_{U, \gt} I^U$.
Hence
$i''_{\gt', \gt} I_{\gt'} \subseteq I_\gt$.
\subsubsection{Generic over $\ME$}
\begin{multline*}
 R_\Es = (\Col(\gk^{+6}, j_\Es(\gk)) \times 
	 \C(\gk^{+4}, j_\Es(\gk)^+) \times
		\C(\gk^{+5}, j_\Es(\gk)^{++}) \times
	\\
	 \C(\gk^{+6}, j_\Es(\gk)^{+3}))
				_{\MSE[G_{\upto \gk}]}.
\end{multline*}
We are going to show that there is $I_\Es \in V$,
which is $R_\Es$-generic over $\ME$.
Moreover, whenever $\gt < \len(\Es)$ we have
	$i''_{\gt, \Es}I_{\gt} \subseteq I_\Es$.
There are $2$ different cases
\begin{enumerate}
\item 
        $\len(\Es)$ is limit:
	We set $I_\Es$ to be the filter generated by 
	$\bigunion_{\gt < \len(\Es)} i''_{\gt, \Es} I_\gt$.
        Its'  genericity is rather straightforward. Let $D \in \ME$ be
        dense open in $R_\Es$.

        Then there is $D_\gt \in \Mt$ such that $i_{\gt, \Es} (D_\gt) = D$.
        By genericity of $I_\gt$, there is $p \in D_\gt \intersect I_\gt$.
        Hence $i_{\gt, \Es}(p) \in D \intersect i''_{\gt, \Es} I_\gt$.
        That is $D \intersect I_\Es \not= \emptyset$.
\item
	$\len(\Es) = \gt + 1$: Just follow the proof of \ref{MtGeneric}
	with $V^*$, $\MSt$ instead of $\MSE$, $\MStE$.
\end{enumerate}
\subsection{Cardinal structure in $M_\gt[I_\gt]$, $\ME[I_\Es]$}
\label{prototype}
The following lifting says everything which we can possibly say.
\begin{align*}
\begin{diagram}
\node{} 
	\node{}
        \node{\ME[I_\Es]}
\\
	\node{N[I_U]}
		 \arrow{e,b}{i^*_{U, \gt'}} 
	\node{M_{\gt'}[I_{\gt'}]}
		 \arrow{ne,t,3}{i^*_{\gt', \Es}} 
		 \arrow{e,b}{i^*_{\gt', \gt}} 
	\node{\Mt[I_\gt]} 
		\arrow[1]{n,b}{i^*_{\gt, \Es}} 
\end{diagram}
\end{align*}
We use these embeddings only in this subsection and not carry them on.

The forcing notion we define later, $\PE$, adds a club to $\gk$.
For each $\gn_1, \gn_2$ successive points in the club
the cardinal structure and power function in the range $[\gn_1^+, \gn_2^{+3}]$
of the generic extension
is the same as the cardinal structure and power function in the range
$[\gk^+, j_\Es(\gk)^{+3}]$ of $\ME[I_\Es]$.
\subsection{Generic filters over iterated ultrapowers} \label{IterateUltrapower}
We iterate $j_\Es$ and consider the following diagram
\begin{align*}
\begin{diagram}
\node{V}
		\arrow[2]{e,t}{j_\Es = j^{0,1}_\Es}
		\arrow{se,t,1}{j_{\gt_1}}
		\arrow{s,l}{i_{U}}
	\node{}
        \node{\ME}
		\arrow[2]{e,t}{j^{1,2}_\Es}
		\arrow{se,t,1}{j^2_{\gt_2}}
		\arrow{s,l}{i^2_{U}}
	\node{}
        \node{M_\Es^2}
		\arrow[2]{e,t}{j^{2,3}_\Es}
		\arrow{se,t,1}{j^3_{\gt_3}}
		\arrow{s,l}{i^3_{U}}
	\node{}
        \node{M_\Es^3}
		\arrow[1]{e,..}
\\
	\node{N}
		\arrow[1]{e,b}{i_{U, \gt_1}} 
		\arrow[1]{nee,t,3}{i_{U, \Es}} 
	\node{M_{\gt_1}}
		\arrow[1]{ne,b,1}{i_{\gt_1, \Es}} 
	\node{N^2}
		\arrow[1]{e,b}{i^2_{U, \gt_2}} 
		\arrow[1]{nee,t,3}{i^2_{U, \Es}} 
	\node{M^{2}_{\gt_2}}
		\arrow[1]{ne,b,1}{i^2_{\gt_2, \Es}} 
	\node{N^3}
		\arrow[1]{e,b}{i^3_{U, \gt_3}} 
		\arrow[1]{nee,t,3}{i^3_{U, \Es}} 
	\node{M^{3}_{\gt_3}}
		\arrow[1]{ne,b,1}{i^3_{\gt_3, \Es}} 
\end{diagram}
\end{align*}
\begin{align*}
& j^0_\Es = \id,\ j^{n}_\Es = j^{0, n}_\Es,
\\
& M^0_\Es = V,\ M^1_\Es = M_\Es,
\\
& G^{\Es, 0} = G_\ugk, \ H^{\Es, 0} = H,
\\
& G^{\Es, n} = j^n_\Es(G^{\Es, 0}), \ H^{\Es, n} = j^n_\Es(H^{\Es, 0}),
\\
& G^{\Es, n} = G^{\Es, n}_{\upto \gk_{n-1}} \times G^{\Es, n}_{\gk_{n-1}} \times 
				G^{\Es, n}_{\downto \gk_{n-1}},
\\
& \gk_0 = \gk, \ \gk_n = j^n_\Es(\gk),
\\
& j^{n, n+1}_{\gt_{n+1}} \func M^{n}_\Es \to M^{n+1}_{\gt_{n+1}} \simeq 
				\Ult(M^{n}_\Es, j_\Es^n(E)(\gt_{n+1})),
\\
& M^{n}_\Es = M^{*n}_\Es[G^{\Es, n}][H^{\Es, n}],
\\
& j^{m, n}_\Es = j^{n-1, n}_\Es \circ \dotsb \circ j^{m+1, m+2} \circ 
			j^{m ,m+1}_\Es,
\\
& R^1_\Es = R_\Es, \ R^n_\Es = j^{n-1}_\Es(R_\Es),\ 
 	 R^1_U = R_U, \ R^n_U = j^{n-1}_U(R_U),
\\
& I^1_\Es = I_\Es, \ I^{n}_\Es = j^{n - 1}_\Es(I_\Es),\ 
	 I^1_U = I_U, \ I^{n}_U = j^{n - 1}_U(I_U).
\end{align*}
We note that
\begin{multline*}
 R^{n}_\Es = (\Col(j^{n - 1}_\Es(\gk)^{+6}, j^{n}_\Es(\gk)) \times 
	 \C(j^{n - 1}_\Es(\gk)^{+4}, j^n_\Es(\gk)^+) \times
		\C(j^{n - 1}_\Es(\gk)^{+5}, j^n_\Es(\gk)^{++}) \times
	\\
	 \C(j^{n - 1}_\Es(\gk)^{+6}, j^n_\Es(\gk)^{+3}))
			_{M^{* n + 1}_\Es[G^{\Es, n}_{\upto \gk_{n-1}}]}.
\end{multline*}
We claim that $I_\Es \times I^2_\Es \times \dotsb \times I^n_\Es$ is
$R_\Es \times R^2_\Es \times \dotsb \times R^n_\Es$-generic over
	$M^{n+1}_\Es$. Of course genericity over
		$M^{n}_\Es$ is more than enough for this.
Moreover, if $D \in M^{n+1}_\Es$ is dense open in $R^n_\Es$
then there is $j_\Es(f)(\Es_\ga \restricted \gt) \in I_\Es$ such that
$j^n_\Es(f)(j^{n-1}_\Es(\Es_\ga \restricted \gt)) \in D \intersect I^n_\Es$.

By construction, $I_\Es$ is $R_\Es$-generic over $M_\Es$. 
Formally we have
\begin{align*}
& V \satisfies \formula{
	I_\Es \text{ is } R_\Es\text{-generic over } \ME
}.
\end{align*}
Applying $j^{n-1}_\Es$ we get
\begin{align*}
& M^{n-1}_\Es \satisfies \formula{
	I^{n}_\Es \text{ is } R^{n}_\Es \text{-generic over } M^{n}_\Es
}.
\end{align*}
As for the product forcing we note that
\begin{align*}
& G^{\Es, n}_{\upto \gk_{n-1}} = 
	G_\ugk \times G^{\Es, 1}_\gk \times G^{\Es, 1}_\dgk \times
		G^{\Es, 2}_{\gk_1} \times G^{\Es, 2}_{\downto \gk_1} \times
		\dotsb \times
		G^{\Es, n-1}_{\gk_{n-2}} \times G^{\Es, n-1}_{\downto \gk_{n-2}},
\\
& G^{\Es, n} = G^{\Es, n}_{\upto \gk_{n-1}} 
		\times
		G^{\Es, n}_{\gk_{n-1}} \times G^{\Es, n}_{\downto \gk_{n-1}}.
\end{align*}
Let use assume, by induction, that $I_\Es \times I^2_\Es\times \dotsb
\times I^{n-1}_\Es$ is $R_\Es \times R^2_\Es \times \dotsb 
\times R^{n-1}_\Es$-generic
over $M^{n-1}_\Es$. Then it is also generic over
$M^{n}_\Es$. 
As
$
M^{*n}_\Es[G^{\Es, n}_{\upto \gk_{n-1}}] \satisfies \formula {
	(\GN{Q}^{\Es, n}_{\gk_{n-1}} *
	\GN{P}^{\Es, n}_{\downto \gk_{n-1}} *
	\GN{Q}^{\Es, n}_{\gk_{n}})
				[G^{\Es, n}_{\upto \gk_{n-1}}]
	 \times
		R^n_\Es \text{ is } \gk_{n-1}^{+} \text{-closed}
}
$, all anti-chains of $R_{\Es} \times \dotsb \times R^{n-1}_{\Es}$
appearing in 	$M^{*n}_\Es [G^{\Es, n}_{\upto \gk_{n-1}}]
			    [G^{\Es, n}_{\gk_{n-1}}]
			    [G^{\Es, n}_{\downto \gk_{n-1}}]
			    [H^{\Es, n}]
			     [I^n_\Es]$
are already in $M^{*n}_\Es [G^{\Es, n}_{\upto \gk_{n-1}}]$.
That is $I^1_\Es \times \dotsb \times I^{n-1}_\Es$ is
	$R^1_\Es \times \dotsb \times R^{n-1}_\Es$-generic
over
	$M^{*n}_\Es [G^{\Es, n}_{\upto \gk_{n-1}}]
		    [G^{\Es, n}_{\gk_{n-1}}]
		    [G^{\Es, n}_{\downto \gk_{n-1}}]
		    [H^{\Es, n}]
		    [I^n_\Es]$.
So we get what we need: $I^1_\Es \times \dotsb \times I_\Es$ is
	$R^1_\Es \times \dotsb \times R^{n}_\Es$-generic
over
	$M^{n}_\Es$.
We are left to prove the `moreover' part.

We start by showing that $j^{n-1 \prime \prime}_\Es I_\Es$ is dense in
$j^{n-1}_\Es(I_\Es) = I^n_\Es$. 
For this we point out that $i''_{U, \Es} I_U$ is dense $I_\Es$.
By elementarity we get that $i^{n \prime \prime}_{U, \Es} I^n_U$ is dense in
	$I^n_\Es$.
So, it is enough to show that $j^{n-1 \prime \prime}_\Es I_U$ is dense in
$j^{n-1}_\Es (I_U) = I^n_U$.

The proof is by induction and we start with 
$n = 2$.
Let $p \in I^2_U$. Choose $f$ such that 
				$j_\Es(f)(\Es_\ga \restricted \gt) = p$.
Then $A = \setof {\gns} {f(\gns) \in I_U} \in E_\ga(\gt)$.
Let $B = \setof {f(\gns)} {\gns \in A}$.
As $N \supset N^\gk$ we get that $B \in N$.
As $N \satisfies \formula {R_U \text{ is } \gk^+ \text{-closed}}$
we get that there is $q \in I_U$ such that $\forall \gns \in A$ 
$q \leq f(\gns)$.
Hence $j_\Es(q) \leq j_\Es(f)(\Es_\ga \restricted \gt)$.
By this we proved $j''_\Es I_U$ is dense in $I^2_U$.

Of course, by elementarity $j^{n-1,n \prime \prime}_\Es I^{n}_U$ is dense
in $j^{n-1, n}_\Es (I^n_U) = I^{n + 1}_U$. Our induction hypothesis is that
	$j^{n-1 \prime \prime}_\Es I_U$ is dense in $I^{n}_U$.
Hence $j^{n-1, n \prime \prime} \circ j^{n-1 \prime \prime}_\Es I_U$ is 
dense in $I^{n + 1}_U$ as needed.

By the above, if $D \in M^{n}_\Es$ is dense open in $R^n_\Es$, then
$D \intersect j^{n-1 \prime \prime}_\Es I_\Es \not= \emptyset$.
By showing that $j^{n - 1 \prime \prime}_\Es I_\Es$ is of the required 
structure we finish the claim. So, 
let $p \in I_\Es$. 

Then $p = j_\Es(f)(\Es_\ga \restricted \gt)$.
Formally, $V \satisfies \formula{ p = j_\Es(f)(\Es_\ga \restricted \gt)}$.
Applying $j^{n-1}_\Es$ we get 
	$M^{n-1}_\Es \satisfies \formula{ j^{n-1}_\Es(p) = 
		j^{n-1,n}_\Es(j^{n-1}_\Es(f))(j^{n-1}_\Es(\Es_\ga 
				\restricted \gt))}$.
That is, $j^{n-1}_\Es(p) = j^{n}_\Es(f)(j^{n-1}_\Es(\Es_\ga \restricted 
			\gt))$.
\ifnum\article=0
\newpage
\fi
\section{Redefining extender sequences} \label{ExtenderSequenceRedefinition}
Our starting assumption in this section is the models and embeddings
constructed in the previous section.
The extender sequences we define here are based on the old ones and
have the same length. They differ in that we add generic filters
into the sequence.

We construct a new extender sequence system, $\Fs$,
from $\Es$. 
If $\len(\Es) = 0$ then
\begin{align*}
& \forall \ga \in \dom \Es\ \Fs_\ga = \ordered{\ga}.
\end{align*}

If $\len(\Es) = 1$  we set $F(0) = E(0)$. 
According to the previous section construction, there is $I(0) \in V$
which is $R_0$-generic over $M_0$. We set
\begin{align*}
& \forall \ga \in \dom \Es\ \Fs_\ga = \ordered{\ga, F(0), I(0)}.
\end{align*}
By $I(\Fs)$ we mean $I(0)$.

We continue by induction. Assume we have defined 
$\ordof{F(\gt'), I(\gt')} {\gt' < \gt}$.

If $\gt = \len(\Es)$ then we set
\begin{align*}
& \forall \ga \in \dom \Es\ 
	\Fs_\ga = \ordof{\ga, F(0), I(0), \dotsc, F(\gt'), I(\gt'), \dotsc}
			{\gt' < \gt}.
\end{align*}
We define $I(\Fs)$ as follows
\begin{enumerate}
\item
	$\gt$ is limit:
	By $I(\Fs)$ we mean the filter generated by
	$\bigunion_{\gt' < \gt} i''_{\gt', \Es} I(\gt')$.
	This filter is $R_\Es$-generic over $\ME$.

\item
	$\gt = \gt' + 1$: $I(\Fs)$ is $I(\gt')$. Note that in this case
	$\ME = \Mt$, so $I(\Fs)$ is $R_\Es$-generic over $\ME$.
\end{enumerate}
If $\gt < \len(\Es)$ then 
we define
\begin{align*}
& A \in F_\ga(\gt) \iff
	\ordof{\ga, F(0), I(0), \dotsc, F(\gt'), I(\gt'), \dotsc}
		{\gt' < \gt} \in j_\Es(A).
\end{align*}
If $\gt + 1 < \len(\Es)$ then there is $I(\gt) \in \ME$ which is $R_\gt$-generic
over $\Mt$. If $\gt + 1 = \len(\Es)$ then there is $I(\gt) \in V$ which
is $R_\gt$-generic over $\Mt$.

By this we finished the definition of the new extender sequence.
We point out that $\Ult(V, E(\gt)) = \Ult(V, F(\gt))$.
This is due to $I(\gt') \in \Mt$ when $\gt' < \gt$.
Hence, the $F(\gt)$'s, do not `pull' into $\Mt$ sets which were not
already pulled in by the $E(\gt)$.

From now on we continue with this new definition of extender sequence
and we use $\Es$ for the new extender sequence system constructed.
%
%
\begin{definition}
We  say $T \in \bar{E}_\ga$ if $\forall \gx < \len(\bar{E}_\ga) \ 
T\in E_\ga(\gx)$.
\end{definition}
\begin{note}
In \cite{MeselfPublishedPhdI} we defined here what is an $\Es_\ga$-tree. In this work we  use
just sets which are in $\Es_\ga$. We use the letters $S$, $T$, $R$, etc.
(used for trees in \cite{MeselfPublishedPhdI}) for these sets.
\end{note}
The operations defined next are the substitute for $T_{\ordered{\gns}}$,
and $T(\gns)$ from \cite{MeselfPublishedPhdI}.
\begin{definition}
\begin{align*}
& T \setminus \gns = T \setminus V_{\gk^0(\gns)},
\\
& T \restricted \gns = T \intersect V_{\gk^0(\gns)}.
\end{align*}
\end{definition}
We define next a form of diagonal intersection which works well
also for the non-normal measures.
\begin{definition}
Assume $\forall \gx<\gk$ $T^\gx \subseteq V_{\gk}$ such that
the elements of $T^\gx$ are extender sequences.
Then $\dsintersect_{\gx < \gk} T^\gx= \setof {\gns \in V_\gk} 
		{\forall \gx < \gk^0(\gns) \ \gns \in T^\gx}$.
\end{definition}
Obviously if $\forall \gx<\gk$ $T^\gx \in \Es_\ga$ then
$\dsintersect_{\gx < \gk} T^\gx \in \Es_\ga$.
\begin{definition}
$S \subseteq [V_\gk]^m$ is called an $\Es_\ga$-fat tree if
\begin{enumerate}
\item 
	$\exists \gx < \len(\Es_\ga)$ $\Lev_0(S) \in \Es_\ga(\gx)$,
\item
	$\forall \ordered{\gns_1, \dotsc, \gns_n} \in S$ 
	$\exists \gx < \len(\Es_\ga)$ $\Suc_S(\ordered{\gns_1, \dotsc, \gns_n} ) 
			\in \Es_\ga(\gx)$.
\end{enumerate}
\end{definition}
%
%
%
%
%
%
%
%
\ifnum\article=0
\newpage
\fi
\section{$P_{\Es}$-Forcing} \label{PEForcing}
The following definition is for the degenerate case. It is
the equivalent of adding a singleton in Radin forcing.
Note that extension in this forcing means only extending
the functions.
%
%
\begin{definition}
Assume $\Es$ is an extender sequence system such that 
$\len(\Es)=0$. A condition $p$ in $P^*_{\Es}$ is of the form
\begin{equation*}
\set{\ordered{\Es_\gk, p^{\Es_\gk}}} \union
        \set{\ordered{\Es_\ga, f}}
\end{equation*}
where
\begin{enumerate}
\item
	$p^{0} \in V_{\gk^0(\Es)}$ is an extender sequence
	(We allow $p^{\Es_\gk}=\emptyset$.). We write $p^0$ for
        $p^{\min s}$ (i.e. $p^{\Es_\gk}$),
\item 
        $f \in R(\gk(p^0), \gk^0(\Es))$.
	If $p^0 = \emptyset$ 	then $f = \emptyset$.
\end{enumerate}
We write $f^p$, $\mc(p)$, $\supp p$, for $f$, $\Es_\ga$, $\set{\gk}$,
respectively.
\end{definition}
%
%
\begin{definition}
Assume $\len(\Es)=0$. 
Let $p, q \in P^*_{\Es}$. We say that $p$ is a Prikry extension of $q$
($p \leq^* q$) if
\begin{enumerate}
\item $\mc(p) = \mc(q)$,
\item $p^0 = q^0$,
\item $f^p \leq f^q$.
\end{enumerate}
\end{definition}
\par\noindent
When $\len(\Es)=0$ the order $\leq^{**}$ disappears:
\begin{definition}
Assume $\len(\Es)=0$. 
Let $p, q \in P^*_{\Es}$. We say that $p \leq^{**} q$
 if $p = q$.
\end{definition}
\par\noindent
Clearly $\ordered{P_{\Es}^*/p,\leq^*} \simeq R(\gk(p^0), \gk^0(\Es))/f^p$.
\begin{definition}
Assume $\len(\Es) > 0$. A condition $p$ in $P_{\Es}^*$ is of the form
\begin{equation*}
\setof{\ordered{\ggs, p^\ggs}}{\ggs \in s} \union
        \set{\ordered{\Es_\ga, T, f, F}}
\end{equation*}
where
\begin{enumerate}
\item
	$s \subseteq \Es$, $\power{s} \leq \kappa$,
		$\min \Es\in s$,
\item
	$p^{\Es_\gk} \in V_{\gk^0(\Es)}$ is an extender sequence
	(We allow $p^{\Es_\gk}=\emptyset$.). We write $p^0$ for
	$p^{\min s}$ (i.e. $p^{\Es_\gk}$),
\item
	$\forall \ggs \in s \setminus \set{\min s}$
	$p^\ggs \in [V_{\gk^0(\Es)}]^{\upto \gw}$ is a sequence of 
        extender sequences where $\gk^0(p^\ggs)$ is increasing.
	(We allow $p^\ggs=\emptyset$.),
\item
	$\forall \ggs \in s$ $\gk(p^0) \leq \max \gk^0(p^{\ggs})$,
\item
	$T \in \Es_{\ga}$,
\item
	$\forall \gns \in T \ 
        \power{ \setof{\ggs \in s}{\max \gk^0(p^\ggs) < \gk^0(\gns)}}
                                                \leq \gk^0(\gns)$,
\item
	$\forall \ggs \in s$ $\Es_\ga \Ege \ggs$,
\item
	$\forall \gbs,\ggs \in s$
	$\forall \gns \in T$ 
        $\max \gk^0(p^\gbs), \max \gk^0(p^\ggs) < \gk^0(\gns) \implies$
        $\gp_{\Es_\ga, \gbs}(\gns) \not =
         \gp_{\Es_\ga, \ggs}(\gns)$,
\item
	$f$ is a function such that
	\begin{enumerate}
	\item
		$\dom f = \setof{\gns \in T} {\len(\gns) = 0}$,
	\item
		$f(\gn_1) \in R(\gk(p^0), \gn_1^0)$.
		If $p^0 = \emptyset$ then $f(\gn_1) = \emptyset$.
	\end{enumerate}
\item
	$F$ is a function such that
	\begin{enumerate}
	\item
		$\dom F = [\setof{\gns \in T} {\len(\gns) = 0}]^2$,
	\item
		$F(\gn_1,\gn_2) \in R(\gn_1^0, \gn_2^0)$,
	\item
		$j^2_{\Es}(F)(\ga, j_{\Es}(\ga)) \in I(\Es)$.
	\end{enumerate}
\end{enumerate}
As usual we write $f^p$, $F^p$, $T^p$, $\mc(p)$, $\supp p$ for 
$f$, $F$, $T$, $\Es_\ga$, $s$ respectively.
Note that we do not
require $\Es_\ga \in s$. That is, we do not consider $\mc(p)$ a part of the
support.
\begin{figure}[htb]
\begin{align*}
\begin{diagram}
\node[4]{} \node{}       \node[5]{}           \node[5]{}           \node[5]{}           \node[5]{\gms_{4,3}} \node[5]{} 
\\
\node[4]{} \node{}       \node[5]{\gms_{1,2}} \node[5]{}           \node[5]{\gms_{3,2}}  \node[5]{\gms_{4,2}} \node[5]{} 
\\
\node[4]{} \node{}       \node[5]{\gms_{1,1}} \node[5]{\gms_{2,1}} \node[5]{\gms_{3,1}} \node[5]{\gms_{4,1}} \node[5]{} 
\\
\node[4]{} \node{\gms_0} \node[5]{\gms_{1,0}} \node[5]{\gms_{2,0}} \node[5]{\gms_{3,0}} \node[5]{\gms_{4,0}} \node[5]{T, f, F} 
\\
\node[4]{} \arrow[27]{e,-} \node[5]{}
\\
\node{\text{Support}} \node[3]{} \node{\Es_\gk} \node[5]{\Es_{\ga_1}} \node[5]{\Es_{\ga_2}} \node[5]{\Es_{\ga_3}} \node[5]{\Es_{\ga_4}} \node[5]{\Es_{\ga_5}= \mc} 
\end{diagram}
\end{align*}
\caption{An example of $p_0 \in P^*_{\Es}$.}
\label{ExampleCondition*}
\end{figure}
We note that the properties of $R_\Es = j^2_\Es(R)(\gk, j_\Es(\gk))$ we 
use are the $\gk^{+}$-closedness and $j_{\Es}(\gk)^+$-c.c.
Any forcing notion satisfying these requirements can be used
instead.
\end{definition}
%
%
\begin{definition} \label{PrikryOrder}
Assume $\len(\Es)>0$.
Let $p, q \in P^*_{\bar{E}}$. We say that $p$ is a Prikry extension of $q$
($p \leq^* q$) if
\begin{enumerate}
\item
	$\supp p \supseteq \supp q$,
\item
	$\mc(p) \Ege \mc(q)$,
\item
	If $\mc(p) \Egt \mc(q)$ then $\mc(q) \in \supp p$,
\item
	$\forall \ggs \in \supp q \ p^\ggs = q^\ggs$,
\item   
        \label{NewPrikryRestriction}
	$\forall \ggs \in \supp p \setminus \supp q \ 
		\max \gk^0(p^\ggs) > \union \union j_\Es(f^{\prime q})(\gk(\mc(q)))$
	where $f^{\prime q}$ is the collapsing part of $f^{q}$,
\item
	$T^p \subseteq \gp^{-1}_{\mc(p), \mc(q)}\,T^q$,
\item
	$\forall \ggs \in \supp q$  $\forall \gns \in T^p$
	\begin{align*}
	\max \gk^0(p^\ggs) < \gk^0(\gns) \implies
	\gp_{\mc(p), \ggs}(\gns) = 
		\gp_{\mc(q), \ggs}(\gp_{\mc(p), \mc(q)}(\gns)),
	\end{align*}
\item
	$\forall \gn_1 \in \dom f^p$
		$f^p(\gn_1) \leq f^q \circ \gp_{\mc(p), \mc(q)}(\gn_1)$,
\item
	$\forall \ordered{\gn_1, \gn_2} \in \dom F^p$ 
	 	$F^p(\gn_1, \gn_2) \leq 
			F^q \circ \gp_{\mc(p), \mc(q)}(\gn_1, \gn_2)$.
\end{enumerate}
The requirement \ref{NewPrikryRestriction} is essential for the proof
of the homogeneity of dense open subsets and hence to the
proof of Prikry's condition.
\end{definition}
\begin{figure}[htb]
\begin{multline*}
\begin{diagram}
\node[3]{} \node{}       \node[2]{\gns_{0,3}}    \node[3]{}           \node[3]{}              \node[3]{}           \node[5]{}           \node[5]{\gms_{4,3}} \node[5]{}              \node[5]{}
\\
\node[3]{} \node{}       \node[2]{\gns_{0,2}}    \node[3]{\gms_{1,2}} \node[3]{}              \node[3]{}           \node[5]{\gms_{3,2}} \node[5]{\gms_{4,2}} \node[5]{}              \node[5]{\gns_{2,2}}
\\
\node[3]{} \node{}       \node[2]{\gns_{0,1}}    \node[3]{\gms_{1,1}} \node[3]{\gns_{1,1}}    \node[3]{\gms_{2,1}} \node[5]{\gms_{3,1}} \node[5]{\gms_{4,1}} \node[5]{\gms_{5,1}}    \node[5]{\gns_{2,1}}
\\
\node[3]{} \node{\gms_0} \node[2]{\gns_{0,0}}    \node[3]{\gms_{1,0}} \node[3]{\gns_{1,0}}    \node[3]{\gms_{2,0}} \node[5]{\gms_{3,0}} \node[5]{\gms_{4,0}} \node[5]{\gms_{5,0}}    \node[5]{\gns_{2,0}}
\\
\node[3]{} \arrow[32]{e,-} \node[5]{}
\\
\node[3]{} \node{\Es_{\gk}} \node[2]{\Es_{\gb_0}} \node[3]{\Es_{\ga_1}} \node[3]{\Es_{\gb_1}} \node[3]{\Es_{\ga_2}} \node[5]{\Es_{\ga_3}} \node[5]{\Es_{\ga_4}} \node[5]{\Es_{\ga_5}} \node[5]{\Es_{\gb_2}}
\end{diagram}
\\
\begin{diagram}
\node[10]{\gp_{\gb_3,\ga_5}^{-1}T, f \circ \gp_{\gb_3,\ga_5}, F \circ \gp_{\gb_3,\ga_5}}
\\
\arrow[16]{e,-}
\\
\node[5]{\Es_{\gb_3}= \mc} 
\end{diagram}
\end{multline*}
\caption{A direct extension of $p_0$ from
                figure \ref{ExampleCondition*}.}
\end{figure}
The order $\leq^{**}$ we  define now allows only shrinkage
of the measure $1$ set. Everything else is the same.
\begin{definition}
Assume $\len(\Es)>0$.
Let $p, q \in \PES$. We say $p \leq^{**} q$
if
\begin{enumerate}
\item
	$\supp p = \supp q$,
\item
	$\mc(p) = \mc(q)$,
\item
	$\forall \gga \in \supp q \ p^\gga=q^\gga$,
\item
	$T^p \subseteq T^q$,
\item
	$\forall \gn_1 \in \dom f^p$ 
		$f^p(\gn_1) = f^q(\gn_1)$,
\item
	$\forall \ordered{\gn_1, \gn_2} \in \dom F^p$ 
 	$F^p(\gn_1, \gn_2) = F^q (\gn_1, \gn_2)$.
\end{enumerate}
\end{definition}
%
%
\begin{definition}
A condition in $P_{\Es}$ is of the form
\begin{equation*}
p_n \append \dotsb \append p_0
\end{equation*}
where
\begin{itemize}
\item
	$p_0 \in P^*_{\Es}$,
	$\gk^0(p_0^0) \geq \gk^0(\gms_1)$,
\item 
	$p_1 \in P^*_{\gms_1}$, 
	$\gk^0(p_1^0) \geq \gk^0(\gms_2)$,
\item 
	$\vdots$
\item 
	$p_n \in P^*_{\gms_n}$,
\end{itemize}
where $\Es, \gms_1, \dotsc, \gms_n$ are extender sequence
systems satisfying
\begin{align*}
	\gk^0(\gms_n) < \dotsb < \gk^0(\gms_1) < \gk^0(\Es).
\end{align*}
When $p = p_n \append \dotsb \append p_0$ we use the short cut
$p_{k..l}$ for $p_k \append \dotsb \append  p_l$.
\end{definition}
\begin{figure}[htb]
\begin{multline*}
\begin{diagram}
\node{\gt_0}  \node[4]{g}
\\
\arrow[8]{e,-}
\\
\node{\gms_{0,0}} \node[4]{\gms_5=\mc}
\end{diagram}
\\
\begin{diagram}
\node{}           \node[3]{}           \node[3]{\gms_{2,3}} \node[3]{}           \node[3]{} \node[4]{}
\\
\node{}           \node[3]{\gms_{1,2}} \node[3]{\gms_{2,2}} \node[3]{}           \node[3]{\gms_{4,2}} \node[4]{}
\\
\node{}           \node[3]{\gms_{1,1}} \node[3]{\gms_{2,1}} \node[3]{\gms_{3,1}} \node[3]{\gms_{4,1}} \node[4]{}
\\
\node{\gms_{0,0}} \node[3]{\gms_{1,0}} \node[3]{\gms_{2,0}} \node[3]{\gms_{3,0}} \node[3]{\gms_{4,0}} \node[4]{S, h, H}
\\
\arrow[20]{e,-} \node[13]{}
\\
\node{\gns_{0,0}} \node[3]{\gns_{1,0}} \node[3]{\gns_2} \node[3]{\gns_3} \node[3]{\gns_{4,0}} \node[4]{\gns_5=\mc}
\end{diagram}
\\
\begin{diagram}
\node{}           \node[3]{}           \node[3]{}           \node[3]{\gns_{6,3}} \node[3]{}           \node[4]{} 
\\
\node{}           \node[3]{}           \node[3]{\gns_{5,2}} \node[3]{\gns_{6,2}} \node[3]{}           \node[4]{} 
\\
\node{}           \node[3]{}           \node[3]{\gns_{5,1}} \node[3]{\gns_{6,1}} \node[3]{}           \node[4]{} 
\\
\node{\gns_{0,0}} \node[3]{\gns_{1,0}} \node[3]{\gns_{5,0}} \node[3]{\gns_{6,0}} \node[3]{\gns_{4,0}} \node[4]{T,f,F} 
\\
\arrow[20]{e,-} \node[2]{}
\\
\node{\Es_{\gk}} \node[3]{\Es_{\ga_1}} \node[3]{\Es_{\ga_2}} \node[3]{\Es_{\ga_3}} \node[3]{\Es_{\ga_4}} \node[4]{\Es_{\ga_5} = \mc} 
\end{diagram}
\end{multline*}
\caption{An Example of a Condition in $P_{\Es}$.}
\end{figure}
\begin{definition}
Let $p,q \in P_{\Es}$. We say that $p$ is a Prikry extension of $q$
($p \leq^* q$)
if $p,q$ are of the form
\begin{equation*}
\begin{split}
p &= p_n \append \dotsb \append p_0,
\\
q &= q_n \append \dotsb \append q_0,
\end{split}
\end{equation*}
and
\begin{itemize}
\item $p_0,q_0 \in P^*_{\Es},\ p_0 \leq^* q_0$,
\item $p_1,q_1 \in P^*_{\gms_1},\ p_1 \leq^* q_1$,
\item $\vdots$
\item $p_n,q_n \in P^*_{\gms_n},\ p_n \leq^* q_n$.
\end{itemize}
\end{definition}
\begin{definition}
Let $p,q \in P_{\Es}$. We say $p \leq^{**} q$ 
if $p,q$ are of the form
\begin{equation*}
\begin{split}
p &= p_n \append \dotsb \append p_0,
\\
q &= q_n \append \dotsb \append q_0,
\end{split}
\end{equation*}
and
\begin{itemize}
\item $p_0,q_0 \in P^*_{\Es},\ p_0 \leq^{**} q_0$,
\item $p_1,q_1 \in P^*_{\gms_1},\ p_1 \leq^{**} q_1$,
\item $\vdots$
\item $p_n,q_n \in P^*_{\gms_n},\ p_n \leq^{**} q_n$.
\end{itemize}
\end{definition}
%
%
$p_{0\ordered{\gns}}$, defined now, is the basic non-direct extension 
in $\PE$ of the condition $p_0$, which adds the extender sequence 
$\gns \in T^{p_0}$ (and hence a condition
$p'_1 \in P_{\gns}$) to the finite sequence.
\begin{definition}
Let $p_0 \in P^*_{\Es}$, $\gns \in T^{p_0}$, 
$\union \union j_\Es(f^{\prime p_0})(\gk(\mc(p_0))) < \gk^0(\gns)$,
where $f^{\prime p_0}$ is the collapse part of $f^{p_0}$. 
We define
$p_{0\ordered{\gns}}$ to be
$p'_1 \append p'_0$ where
\begin{enumerate}
\item
	$\supp p'_0 = \supp p_0$,
\item
	$\forall \ggs \in \supp p'_0 \ p_0^{\prime\ggs}=
        \begin{cases}
        \gp_{\mc(p_0),\ggs}(\gns) & 
                                \max \gk^0(p^\ggs_0) < \gk^0(\gns),\ 
                                \len(\gns) > 0.
        \\
        \gp_{\mc(p_0),\ggs}(\gns) & 
                                \max \gk^0(p^\ggs_0) < \gk^0(\gns),\ 
                                \len(\gns) = 0,\ 
                \\
                                 &   \ggs = \Es_\gk.
        \\
        p^\ggs \append \gp_{\mc(p_0),\ggs}(\gns) & 
                                \max \gk^0(p^\ggs_0) < \gk^0(\gns),\ 
                                \len(\gns) = 0,\ 
                \\
                                &    \ggs \not= \Es_\gk.
        \\
        p^\ggs_0  &       \text{Otherwise.}
        \end{cases}$,
\item
	$\mc(p'_0) = \mc(p_0)$,
\item
	$T^{p'_0}=T^{p_0} \setminus \gns$,
\item
	$\forall \gn_2 \in T^{p'_0}$ $f^{p'_0}(\gn_2) = F^{p_0}(\gk(\gns), \gn_2)$,
\item
	$F^{p'_0} = F^{p_0}$,
\item
	$\supp p'_1=
        \begin{cases}
            \setof{\gp_{\mc(p_0), \ggs}(\gns)}
              {\ggs \in \supp p_0, \, \max \gk^0(p^\ggs_0) < \gk^0(\gns)}
                        &
                                \len(\gns) > 0
                        \\
                \set{\gp_{\mc(p_0),0}(\gns)} & \len(\gns) = 0
        \end{cases}$,
\item
	$\forall \ggs \in \supp p_0$ 
	$\forall \gp_{\mc(p_0), \ggs}(\gns) \in \supp p'_1$ 
                $p_1^{\prime\gp_{\mc(p_0), \ggs}(\gns)}=p^\ggs_0$,
\item
	$\mc(p'_1) = \gns$,
\item
	If $\len(\gns) = 0$ then
	\begin{enumerate}
	\item
		$T^{p'_1} = \emptyset$,
	\item
		$f^{p'_1} = f^{p_0}(\gk(\gns))$,
	\item
		$F^{p'_1} = \emptyset$.
	\end{enumerate}
\item
	If $\len(\gns) > 0$ then
	\begin{enumerate}
	\item
		$T^{p'_1} = T^{p_0} \restricted \gns$,
	\item
		$f^{p'_1} = f^{p_0} \restricted \gns$,
	\item
		$F^{p'_1} = F^{p_0} \restricted \gns$.
	\end{enumerate}
\end{enumerate}
\begin{figure}[htb]
\begin{multline*}
\begin{diagram}
\node{}       \node[5]{}                        \node[6]{}                       \node[5]{\gms_{4,3}}              \node[5]{} 
\\
\node{}       \node[5]{\gms_{1,2}}              \node[6]{\gms_{3,2}}             \node[5]{\gms_{4,2}}              \node[5]{} 
\\
\node{}       \node[5]{\gms_{1,1}}              \node[6]{\gms_{3,1}}             \node[5]{\gms_{4,1}}              \node[5]{}
\\
\node{\gms_0} \node[5]{\gms_{1,0}}              \node[6]{\gms_{3,0}}             \node[5]{\gms_{4,0}}              \node[5]{T\restricted \gns, f\restricted \gns, F\restricted \gns} 
\\
\arrow[23]{e,-} \node[18]{}
\\
\node{\gns^0} \node[5]{\gp_{\ga_5,\ga_1}(\gns)} \node[6]{\gp_{\ga_5,\ga_3}(\gns)} \node[5]{\gp_{\ga_5,\ga_4}(\gns)} \node[5]{\gns}
\end{diagram}
\\
\begin{diagram}
\node{}          \node[4]{}                                  \node[4]{\gms_{2,1}} \node[4]{}                        \node[5]{}                        \node[6]{}
\\
\node{\gns^0}    \node[4]{\gp_{\ga_5,\ga_1}(\gns)}           \node[4]{\gms_{2,0}} \node[4]{\gp_{\ga_5,\ga_3}(\gns)} \node[5]{\gp_{\ga_5,\ga_4}(\gns)} \node[6]{T \setminus \gns, f, F}
\\
\arrow[26]{e,-}  \node[18]{}
\\
\node{\Es_{\gk}} \node[4]{\Es_{\ga_1}}                       \node[4]{\Es_{\ga_2}} \node[4]{\Es_{\ga_3}}             \node[5]{\Es_{\ga_4}}            \node[6]{\Es_{\ga_5} = \mc}
\end{diagram}
\end{multline*}
\caption{$p_{0 \ordered{\protect\gns}}$ for $p_0$ of figure \protect\ref{ExampleCondition*} for $\len(\protect\gns) > 0$.}
\end{figure}
\begin{figure}[htb]
\begin{multline*}
\begin{diagram}
\node{\gms_0} \node[5]{f(\gns)}
\\
\arrow[8]{e,-}
\\
\node{\gns^0} \node[5]{\gns}
\end{diagram}
\\
\begin{diagram}
\node{}       \node[5]{}                         \node[5]{}           \node[5]{}                          \node[5]{\gp_{\ga_5, \ga_4}(\gns)} \node[5]{} 
\\
\node{}       \node[5]{\gp_{\ga_5, \ga_1}(\gns)} \node[5]{}           \node[5]{\gp_{\ga_5, \ga_3}(\gns)}  \node[5]{\gms_{4,3}}               \node[5]{} 
\\
\node{}       \node[5]{\gms_{1,2}}               \node[5]{}           \node[5]{\gms_{3,2}}                \node[5]{\gms_{4,2}}               \node[5]{} 
\\
\node{}       \node[5]{\gms_{1,1}}               \node[5]{\gms_{2,1}} \node[5]{\gms_{3,1}}                \node[5]{\gms_{4,1}}               \node[5]{} 
\\
\node{\gns^0} \node[5]{\gms_{1,0}}               \node[5]{\gms_{2,0}} \node[5]{\gms_{3,0}}                \node[5]{\gms_{4,0}}               \node[5]{T \setminus \gns, f, F} 
\\
\arrow[27]{e,-} \node[5]{}
\\
\node{\Es_\gk} \node[5]{\Es_{\ga_1}}             \node[5]{\Es_{\ga_2}} \node[5]{\Es_{\ga_3}}              \node[5]{\Es_{\ga_4}}              \node[5]{\Es_{\ga_5}= \mc} 
\end{diagram}
\end{multline*}
\caption{$p_{0 \ordered{\protect\gns}}$ for $p_0$ of figure \protect\ref{ExampleCondition*} for $\len(\protect\gns) = 0$.}
\end{figure}
We use the following notation:
$(p_{0\ordered{\gns}})_0 = p'_0$,
$(p_{0\ordered{\gns}})_1 = p'_1$.
If $\ordered{\gns_1, \gns_2} \in [T^{p_0}]^2$ then
$p_{0\ordered{\gns_1, \gns_2}} = (p_{0\ordered{\gns_1}})_1 \append
(p_{0\ordered{\gns_1}})_{0\ordered{\gns_2}}$ and so on.
\end{definition}
\begin{definition}
Let $p,q \in P_{\Es}$. We say that $p$ is a $1$-point extension of $q$
($p \leq^1 q$) if
$p,q$ are of the form
\begin{equation*}
\begin{split}
p &= p_{n+1} \append p_n \append \dotsb \append p_0,
\\
q &= q_n \append \dotsb \append q_0,
\end{split}
\end{equation*}
and there is $0 \leq k \leq n$ such that
\begin{itemize}
\item $p_i,q_i \in P^*_{\gms_i},\ p_i \leq^* q_i$ for $i=0,\dotsc,k-1$,
\item $p_{i+1},q_i \in P^*_{\gms_i},\ p_{i+1} \leq^* q_i$ for $i=k+1,\dotsc,n$,
\item There is $\ordered{\gns} \in T^{q_k}$ such that
        $p_{k+1} \append p_k \leq^* q_{k \ordered{\gns}}$.
\end{itemize}
\end{definition}
%
%
%
\begin{definition}
Let $p,q \in P_{\Es}$. We say that $p$ is an $n$-point extension of $q$
($p \leq^n q$) if there are $p^n, \dotsc, p^0$ such that
\begin{align*}
        p=p^n \leq^1 \dotsb \leq^1 p^0=q.
\end{align*}
\end{definition}
We consider Prikry extension to be a $0$-point extension. That is
\begin{definition}
Let $p,q \in P_{\Es}$. We say that $p$ is a $0$-point extension of $q$
($p \leq^0 q$) if $p \leq^* q$.
\end{definition}
\begin{definition}
Let $p,q \in P_{\Es}$. We say that $p$ is an extension of $q$
($p \leq q$) if there is an $n$ such that $p \leq^n q$.
\end{definition}
\par\noindent Later on by $\PE$ we  mean $\ordered{\PE,\leq}$.
\begin{note}
When $\len(\Es)=1$ the forcing $P_{\Es}$ is
as the forcing defined in \cite{MeselfPublishedMsc}.
\end{note}
\begin{definition}
Let $\ges$ be an extender sequence such that $\gk^0(\ges)< \gk^0(\Es)$.
\begin{align*}
\PE/\Pe = \setof {p}{\exists q\in \Pe \  q\append p \in \PE}.
\end{align*}
\end{definition}
When $\gt_1 < \gt_2$ we have $A \in \Es_\ga\restricted \gt_2
\implies A \in \Es_\ga\restricted \gt_1$. So it is our convention
that $P_{\Es \restricted \gt_2} \subseteq P_{\Es \restricted \gt_1}$.
%
%
%
\ifnum\article=0
\newpage
\fi
\section{Basic properties of $P_{\Es}$} \label{BasicProperties}
\begin{claim}
$P_\Es$ satisfies $\gk^{++}$-c.c.
\end{claim}
\begin{proof}
Let $\setof {p^\gx} {\gx < \gk^{++}} \subset \PE$.
As $\power{ \setof{p^\gx_{n_\gx} \append \dotsb \append p^\gx_1} {\gx < \gk^{++}} } = \gk$
we can assume without loss of generality that there is 
	$p_n \append \dotsb \append p_1$ such that $\forall \gx < \gk^{++}$
$p_n \append \dotsb \append p_1 = p^\gx_{n_\gx} \append \dotsb \append p^\gx_1$.
Hence we can ignore these lower parts and assume that the set of conditions
we start with is $\setof {p^\gx_0} {\gx < \gk^{++}}$.

Let $d_\gx = \supp p_0^\gx \union \set{\mc(p_0^\gx)}$.
As $\power {d_\gx} < \gk^+ < \gk^{++}$, $(\gk^{+})^{\gk} = \gk^+ < \gk^{++}$
we can invoke the $\gD$-lemma. Hence, without loss of generality,
there is $d$ such that $\forall \gx_1 \not= \gx_2$ 
	$d_{\gx_1} \intersect d_{\gx_2} = d$.
As $\gk^\gk = \gk^+ < \gk^{++}$ we can assume, without loss of generality, 
$\forall \gas \in d$
$\forall \gx_1, \gx_2$ $(p^{\gx_1}_0)^\gas = (p^{\gx_2}_0)^\gas$.

The $T^{p^\gx_0}$'s, $F^{p^\gx_0}$'s are always compatible.
Hence we are left with handling of
	$\setof {f^{p^\gx_0}} {\gx < \gk^{++}}$.
$\setof {j_\Es(f^{p^\gx_0})(\gk(\mc(p^\gx_0)))} {\gx < \gk^{++}} \subset
	j_\Es(R)(\gk((p^\gx_0)^0), \gk)$ and
$j_\Es(R)(\gk((p^\gx_0)^0), \gk)$ satisfies $\gk^+$-c.c.
Hence there are $\gx_1, \gx_2$ such that 
\begin{align*}
	j_\Es(f^{p^{\gx_1}_0})(\gk(\mc(p^{\gx_1}_0))) \compatible
		j_\Es(f^{p^{\gx_2}_0})(\gk(\mc(p^{\gx_2}_0))).
\end{align*}
Hence $p^{\gx_1}_0 \compatible p^{\gx_2}_0$.
\end{proof}
%
%
%
%
%
%
%
%
\begin{lemma} \label{ToPreDense-0}
Let $p=p_1 \append p_0 \in \PE$. 
Assume that we have
$S^0$, $t'_1(\gns_1)$
such that
\begin{enumerate}
\item 
	$S^0 \subseteq T^{p_0}$,
\item 
	$S^0  \in \mc(p_0)(\gx_0)$,
\item 
	$\forall \ordered{\gns_{1}}
                \in S^0\ 
	t'_{1}(\gns_1)  \leq^{**}
	         (p_{0 \ordered{\gns_{1}}})_{1}$.
\end{enumerate}
Then there are $p^*_0 \leq^{**} p_0$, $S^{0*} \subseteq S^0$,
such that
\begin{align*}
\setof {
	p_1 \append t'_{1}(\gns_1) \append
         (p^*_{0 \ordered{\gns_{1}}})_0
	}
	{
	\ordered{\gns_1} \in S^{0*}
	}
\end{align*}
is pre-dense below $p_1 \append p^*_0$.
\end{lemma}
\begin{proof}
Let $\Es_\ga = \mc(p_0)$. Set
\begin{align*}
& T_{\upto \gx_0} = j_{\Es}(T^{t'_1})^{(\Es_\ga \restricted \gx_0)},
\\
& T_{ \gx_0} =
	\setof{
		\gns_1 \in S^0
	}
	{
		T_{\upto \gx_0} \restricted \gns_1 = T^{t'_1(\gns_1)}
	},
\\
& T_{\downto \gx_0} = 
	\setof{
		\gns \in T^{p_0}
	}
	{
		\exists \gt < \len(\gns)\ T_{\gx_0} \restricted \gns 
				\in \gns(\gt)
	}.
\end{align*}
We set $S^{0*} = T_{\gx_0}$.
It is clear that $T_{\upto \gx_0} \in \Es_\ga \restricted \gx_0$,
$T_{\gx_0} \in E_\ga(\gx_0)$, $\forall \gx_0 < \gx < \len(\Es)\ 
T_{\downto \gx_0} \in E_\ga(\gx)$. 
Let $T^* = T_{\upto \gx_0} \union T_{\gx_0} \union T_{\downto \gx_0}$.
We define the condition $p^*_0$ to be $p_0$ with $T^*$ substituted
for $T^{p_0}$. 
Let $r \append q \leq p_1 \append p^*_0$ such that
$r \leq p_1$, $q \leq p^*_0$. Of course  the part below
$p_1$ poses no problem. So we are left to show that there is
$\gns \in S^{0*}$ such that $q \compatible t'_1(\gns_1) \append 
(p^*_{0 \ordered{\gns_1}})_0$.
By the definition of $\leq$ there is $\ordered{\gns_1, \dotsc, \gns_n}
\in [T^{p^*_0}]^n$ such that $q=q_n \append \dotsb \append q_0 \leq^* 
p^*_{0 \ordered{\gns_1, \dotsc, \gns_n}}$. 
Let $k$ be the last such that $\ordered{\gns_1, \dotsc, \gns_k} \in
[T_{\upto \gx_0}]^k$. ($k$ can be $0,1,\dotsc,n$.)
If $k = n$ we choose $\gns_{k+1} \in T_{\gx_0} \setminus \gns_{n}$.
We split the handling according to where $\gns_{k+1}$ is.
\begin{enumerate}
\item 
	$\gns_{k+1} \in T_{\gx_0}$:
	$q_n \append \dotsb \append q_{n-k+1} \leq 
		(p^*_{0 \ordered{\gns_{k+1}}})_1$,
	$q_{n-k} \append \dotsb \append q_{0} \leq \linebreak[0]
		(p^*_{0 \ordered{\gns_{k+1}}})_0$.
	As $T^{t'_1(\gns_{k+1})} = T_{\upto \gx_0} \restricted \gns_{k+1}$
	we have
	$\ordered{\gns_1, \dotsc, \gns_k} \in
	[T^{t'_1(\gns_{k+1})}]^k$. Hence
	$q_n \append \dotsb \linebreak[0] \append \linebreak[0] q_{n-k+1} 
		\compatible t'_1(\gns_{k+1})$.
\item
	$\gns_{k+1} \in T_{\downto \gx_0}$:
	So $\exists \gt < \len(\gns_{k+1})\ S^{0*} \restricted \gns_{k+1} \in 
	\gns_{k+1}(\gt)$. In particular
	there is $\gns \in (S^{0*} \intersect 
	T^*) \restricted \gns_{k+1} \setminus \gns_k$.
	Let $s = (p^*_{0 \ordered{\gns_1, \dotsc, \gns_k, \gns}})_1$.
	Then $q_n \append \dotsb \append q_{n-k+1} \append \linebreak[0] s
	\append \linebreak[0]
	q_{n-k} \append \linebreak[0] \dotsb \append q_{0} \leq \linebreak[0]
		p^*_{0 \ordered{\gns}}$.
	 Once more, as $\gns \in T_{\gx_0}$
	we have  $T^{t'_1(\gns)} = T_{\upto \gx_0} \restricted \gns$, hence
	$\ordered{\gns_1, \dotsc, \gns_k} \in T^{t'_1(\gns)}$. That is
	$q_n \append \dotsb \append q_{n-k+1} \append \linebreak[0] s
	\append \linebreak[0]
	q_{n-k} \append \linebreak[0] \dotsb \append q_{0} \compatible \linebreak[0]
		t'_1(\gns) \append (p^*_{0 \ordered{\gns}})_0$.
\end{enumerate}
\end{proof}
\begin{lemma} \label{ToPreDense-2}
\label{ToPreDense}
Let $p=p_1 \append p_0 \in \PE$. 
Assume that we have
$S^0$, $t'_2(\gns_1)$, $t'_1(\gns_1, \gns_2)$
such that
\begin{enumerate}
\item 
	$S^0 \subseteq [T^{p_0}]^2$,
\item 
	$S^0$ is $\mc(p_0)$-fat tree,
\item 
	$\forall \ordered{\gns_1, \gns_2}
                \in S^0\ 
	t'_2(\gns_1) \append t'_1(\gns_1, \gns_2) \leq^{**}
	         (p_{0 \ordered{\gns_1, \gns_2}})_{2..1}$.
\end{enumerate}
Then there are $p^*_0 \leq^{**} p_0$, $S^{0*} \subseteq S^0$,
such that
\begin{align*}
\setof {
	p_1 \append t'_2(\gns_1) \append t'_1(\gns_1, \gns_2)
		\append
         (p^*_{0 \ordered{\gns_1, \gns_2}})_0
	}
	{
	\ordered{\gns_1, \gns_2} \in S^{0*}
	}
\end{align*}
is pre-dense below $p_1 \append p^*_0$.
\end{lemma}
\begin{proof}
Fix $\gns_1 \in \Lev_0(S^0)$. Then we can invoke \ref{ToPreDense-0}
on $t'_2(\gns_1) \append (p_{0 \ordered{\gns_1}})_0$, $t'_1(\gns_1, \gns_2)$,
$\Suc_{S^0}(\gns_1)$ to get 
$r_0(\gns_1) \leq^{**} (p_{0 \ordered{\gns_1}})_0$, $R^{0}(\gns_1) 
\subseteq \Suc_{S^0}(\gns_1)$ such that
\begin{align*}
\setof{
	p_1 \append t'_2(\gns_1) \append t'_1(\gns_1, \gns_2) \append
		(r_0(\gns_1)_{\ordered{\gns_2}})_0
}
{
	\gns_2 \in R^{0}(\gns_1)
}
\end{align*}
is pre-dense below $p_1 \append t'_2(\gns_1) \append 
r_0(\gns_1)$.

We do the above for all $\gns_1 \in \Lev_0(S^0)$.
We let $T^{r_0} = \dsintersect_{\gns_1 \in \Lev_0(S^0)} T^{r_0(\gns_1)}$.
Of course, $r_0$ is $p_0$ with $T^{p_0}$ substituted by $T^{r_0}$.
We point out that $\forall \gns_1 \in \Lev_0(S^0)$ 
$(r_{0 \ordered{\gns_1}})_0 \leq^{**} r_0(\gns_1)$.

We invoke $\ref{ToPreDense-0}$ with $p_1 \append r_0$, $t'_2(\gns_1)$, 
$\Lev_0(S^0)$ to get $p^*_0 \leq^{**} r_0$, $R^{0*} \subseteq \Lev_0(S^0)$ 
such that
\begin{align*}
\setof{
	p_1 \append t'_2(\gns_1) \append
		(p^*_{0 \ordered{\gns_1}})_0
}
{
	\gns_1 \in R^{0*}
}
\end{align*}
is pre-dense below $p_1 \append p^*_0$.
We set $\Lev_0(S^{0*}) = R^{0*}$, $\Suc_{S^{0*}}(\gns_1) = R^0(\gns_1)$.
We claim that
\begin{align*}
\setof{
	p_1 \append t'_2(\gns_1) \append t'_1(\gns_1, \gns_2) \append
		(p^*_{0 \ordered{\gns_1, \gns_2}})_0
}
{
	\ordered{\gns_1, \gns_2} \in S^{0*}
}
\end{align*}
is pre-dense below $p_1 \append p^*_0$.
Let $q \leq p_1 \append p^*_0$.

Then there is $\gns_1 \in R^{0*}$ such that
$q \compatible p_1 \append t'_2(\gns_1) \append (p^*_{0 \ordered{\gns_1}})_0$.
As $(p^*_{0 \ordered{\gns_1}})_0 \leq^{**} r_0(\gns_1)$ we have
$q \compatible p_1 \append t'_2(\gns_1) \append r_0(\gns_1)$.
Choose $s \leq q, p_1 \append t'_2(\gns_1) \append r_0(\gns_1)$.
There is $\gns_2 \in R(\gns_1)$ such that
$s \compatible p_1 \append t'_2(\gns_1) \append t'_1(\gns_1, \gns_2)
	\append (r_0(\gns_1)_{\ordered{\gns_2}})_0$.
As $s \leq p_1 \append p^*_0$, $(p^*_{0 \ordered{\gns_1, \gns_2}})_0 \compatible
(r_0(\gns_1)_{\ordered{\gns_2}})_0$ we get
$s \compatible p_1 \append t'_2(\gns_1) \append t'_1(\gns_1, \gns_2)
	\append (p^*_{0 \ordered{\gns_1, \gns_2}})_0$.
We complete the proof by noting that $S^{0*}$ was constructed such 
that $\ordered{\gns_1, \gns_2} \in S^{0*}$.
\end{proof}
Repeat invocation of the above proof yields
\begin{claim} \label{ToPreDense-n}
Let $p=p_1 \append p_0 \in \PE$. 
Assume we have
$S^0$, $t'_n(\gns_1)$, $\dotsc$, $t'_1(\gns_1, \dotsc, \linebreak[0] \gns_n)$
such that
\begin{enumerate}
\item 
	$S^0 \subseteq [T^{p_0}]^n$,
\item
	$S^0$ is $\mc(p_0)$-fat tree,
\item 
	$\forall \ordered{\gns_1, \dotsc, \gns_n}
                \in S^0\ 
	t'_{n}(\gns_1) \append \dotsb \append t'_1(\gns_1, \dotsc, \gns_n)
		 \leq^{**}
		         (p_{0 \ordered{\gns_1, \dotsc, \gns_n}})_{n..1}$.
\end{enumerate}
Then there are $p^*_0 \leq^{**} p_0$, $S^{0*} \subseteq S^0$,
such that
\begin{align*}
\setof {
	p_1 \append t'_{n}(\gns_1) \append \dotsc \append 
		t'_1(\gns_1, \dotsc, \gns_n)
		\append
         (p^*_{0 \ordered{\gns_1, \dotsc, \gns_n}})_0
	}
	{
	\ordered{\gns_1, \dotsc, \gns_n} \in S^{0*}
	}
\end{align*}
is pre-dense below $p_1 \append p^*_0$.
\end{claim}
%
%
%
%
%
%
\ifnum\article=0
\newpage
\fi
\section{Homogeneity in dense open subsets} \label{HomogenDense}
Our aim in this section is to prove the following theorem. 
Unlike \cite{MeselfPublishedPhdI} we do not carry exact information on the measure $1$ set
of new blocks.
This made us use $\leq^{**}$ in this theorem.
\begin{theorem} \label{FullHomogen}
Let $D$ be dense open in $\PE$, $p = p_{l..0} \in \PE$.
Then there is $p^*  \leq^* p$
such that
\begin{multline*}
\exists S^l\ \forall \ordered{\gns_{l, 1}, \dotsc,\gns_{l, n_l}} \in S^l\ 
\exists t'_{l, n_l} \append \dotsb \append t'_{l, 1} \leq^{**}
                (p^*_{l \ordered{\gns_{l, 1}, \dotsc, \gns_{l, n_l}}})_{n_l..1}
\\
\vdots
\\
\exists S^0\ \forall \ordered{\gns_{0, 1}, \dotsc,\gns_{0, n_0}} \in S^0\ 
\exists t'_{0, n_0} \append \dotsb \append t'_{0, 1} \leq^{**}
                (p^*_{0 \ordered{\gns_{0, 1}, \dotsc, \gns_{0, n_0}}})_{n_0..1}
\\
t'_{l, n_l} \append \dotsb \append t'_{l, 1} \append
                (p^*_{l \ordered{\gns_{l, 1}, \dotsc, \gns_{l, n_l}}})_{0}
		\append
\\
	\vdots
\\
t'_{0, n_0} \append \dotsb \append t'_{0, 1} \append
                (p^*_{0 \ordered{\gns_{0, 1}, \dotsc, \gns_{0, n_0}}})_{0}
		\in D
\end{multline*}
where for each $k=0,\dotsc,l$ $S^k \subseteq [V_{\gk^0(\mc(p^*_k))}]^{n_k}$ 
is an $\mc(p^*_k)$-fat tree.
\end{theorem}
A word of caution is in order. If one of the elements of $p^*$ does not
contain a tree then we mean in the above formula just a direct extension
of it. For example, let $p = p_{2..0} \in \PE$ and $p_2 \in P_{\ges_2}$,
$p_1 \in \Pe$ where $\len(\ges_2) > 0$, $\len(\ges) = 0$. Then the above
formula should be read as
\begin{multline*}
\exists S^2\ \forall \ordered{\gns_{2, 1}, \dotsc,\gns_{2, n_2}} \in S^2\ 
\exists t'_{2, n_2} \append \dotsb \append t'_{2, 1} \leq^{**}
                (p^*_{2 \ordered{\gns_{2, 1}, \dotsc, \gns_{2, n_2}}})_{n_2..1}
\\
\exists S^0\ \forall \ordered{\gns_{0, 1}, \dotsc,\gns_{0, n_0}} \in S^0\ 
\exists t'_{0, n_0} \append \dotsb \append t'_{0, 1} \leq^{**}
                (p^*_{0 \ordered{\gns_{0, 1}, \dotsc, \gns_{0, n_0}}})_{n_0..1}
\\
t'_{2, n_2} \append \dotsb \append t'_{2, 1} \append
                (p^*_{2 \ordered{\gns_{2, 1}, \dotsc, \gns_{2, n_2}}})_{0}
		\append
\\
		p^*_1
		\append
\\
t'_{0, n_0} \append \dotsb \append t'_{0, 1} \append
                (p^*_{0 \ordered{\gns_{0, 1}, \dotsc, \gns_{0, n_0}}})_{0}
		\in D.
\end{multline*}
We prove the theorem by induction on $l$, the number of blocks
in $p$. We  give the proof in a series of lemmas. \ref{DenseHomogen-0},
\ref{DenseHomogenTriv-0} are the case $l = 0$ of the theorem.
%
%
\begin{lemma} \label{FixSupport} \label{AlmostDenseHomogen-0}
Let $D$ be dense open in $\PE/\Pe$, 
$p = p_0 \in \PE/\Pe$, $n < \gw$.
Then there is $p^*_0 \leq^* p_0$
such that one and only one of the following is satisfied
\begin{enumerate}
\item There is  $S \subseteq [T^{p^*_0}]^n$, an 
$\mc(p^*_0)$-fat tree, such that
\begin{multline*}
\forall \ordered{\gns_1,\dotsc,\gns_{n}} \in S\ 
\exists t'_n \append \dotsb \append t'_1 \leq^{**}
                (p^*_{0 \ordered{\gns_1,\dotsc,\gns_{n}}})_{n..1}
\\
                t'_n \append \dotsb \append t'_1 \append
		(p^*_{0 \ordered{\gns_1,\dotsc,\gns_{n}}})_0 \in D.
\end{multline*}
\item
	$\forall \ordered{\gns_1,\dotsc,\gns_{n}} \in 
                        [T^{p^*_0}]^n$
        $\forall q \leq^* p^*_{0 \ordered{\gns_1,\dotsc,\gns_{n}}}$
               $q \notin D$.
\end{enumerate}
\end{lemma}
\begin{proof}
We give the proof for $n=1$. It is essentially the same for all $n$.

Let $\gc$ be large enough so that $H_\gc$ catch `what interests us' and
let
\begin{align*}
& N \subelem H_\gc,
\\
& \power{N} = \gk,
\\
& N \supset \gk,
\\
& N \supseteq N^{\upto \gk},
\\
& \PE, \Pe, D, p_0 \in N.
\end{align*}
Choose $\ga \in \dom \Es$ such that 
\begin{align*}
\forall \gga \in \dom \Es \intersect N \ \gga \Elt \ga,
\end{align*}
and set
\begin{align*}
A = \gp^{-1}_{\ga, \gb} \, T^{p_0}
\end{align*}
where $\Es_\gb = \mc(p_0)$.

Let $\preceq$ be a well ordering of $A$ such that
\begin{align*}
& \forall \gns_1,\gns_2 \in A \ 
	\gns_1 \preceq \gns_2 \implies \gk^0(\gns_1) \leq \gk^0(\gns_2) .
\end{align*}
We shrink $A$ a bit so that the following is satisfied
\begin{align*}
& \forall \gns \in A\ 
	\power{\setof{\gms \in A} {\gk^0(\gms) < \gk^0(\gns)}} 
		\leq \gk^0(\gns).
\end{align*}
We start an induction on $\gns$ in which we  build
\begin{align*}
\ordof{\gas^\gns_0, u^\gns_0, T^\gns_0, F^\gns_0}
      {\gns \in A}.
\end{align*}
Assume that we have constructed
\begin{align*}
\ordof{\gas^\gns_0, u^\gns_0, T^\gns_0, F^\gns_0}
      {\gns \prec \gns_0}.
\end{align*}
Set the following:
\begin{itemize}
\item 
    $\gns_0$ is $\prec$-minimal:
    \begin{align*}
	& q' = p_0 \setminus \set{\ordered{\Es_\gb, T^{p_0}, f^{p_0},
					F^{p_0} }},
	\\
	& \ga' = \gb.
    \end{align*}
\item  
    $\gns_0$ is the immediate $\prec$-successor of $\gns$:
    \begin{align*}
	    & q' = u^\gns_0,
	    \\
	    & \ga' = \ga^\gns.
    \end{align*}
\item 
        $\gns_0$ is $\prec$-limit: Choose $\ga' \in N$ such that
    $\forall \gns \prec \gns_0$ $\ga' \Egt \ga^\gns$ 
    and set
    \begin{align*}
	    q' = \bigunion_{\gns \prec \gns_0} u^\gns_0.
    \end{align*}
\end{itemize}
We start an induction on $i$. 
We construct in it
\begin{align*}
	\ordof{
		\ga^{\gns_0, i}_0,
		u^{\gns_0, i}_0,
		T_0^{\gns_0, i}, 
		f_0^{\gns_0, i},
		F_0^{\gns_0, i},
		\gas^{\gns_0, i}_1,
		u^{\gns_0, i}_1
		T_1^{\gns_0, i}, 
		f_1^{\gns_0, i},
		F_1^{\gns_0, i}
		} 
               {i < \gk}
\end{align*}
such that
$\ordof{j_{\Es}(f_0^{\gns_0, i})({\ga^{\gns_0, i}})} 
{i < \gk}$
is a maximal anti-chain in $j_\Es(R)(\gk^0(\gns_0), \gk)$ below 
$j_\Es(F^{p_0})(\gp_{\ga, \gb}(\gk(\gns_0)), \gb)$.
Assume we have constructed
\begin{align*}
	\ordof{
		\ga^{\gns_0, i}_0,
		u^{\gns_0, i}_0,
		T_0^{\gns_0, i}, 
		f_0^{\gns_0, i},
		F_0^{\gns_0, i},
		\gas^{\gns_0, i}_1,
		u^{\gns_0, i}_1,
		T_1^{\gns_0, i}, 
		f_1^{\gns_0, i},
		F_1^{\gns_0, i}
		} 
               {i < i_0},
\end{align*}
and we do step $i_0$.
\begin{itemize}
\item 
	$i_0 = 0$:
    \begin{align*}
	& q'' = q',
	\\
	& \ga'' = \ga',
	\\
	& f'' = F^{p_0} (\gp_{\ga, \gb}(\gk(\gns_0)), 
					\gp_{\ga'', \gb}(-)).
    \end{align*}
\item 
	$i_0 = i + 1$:
	If $\ordof{j_{\Es}(f_0^{\gns_0, i})(\ga_0^{\gns_0, i})} 
				{i < i_0}$
	is a maximal anti-chain below 
		$j_\Es(F^{p_0})(\gp_{\ga, \gb}(\gk(\gns_0)), \gb)$
	then we finish the induction on $i$.
	Otherwise we choose $f''$, $\gb''$ such that
	$j_\Es(f'')(\gb'') \incompatible 
		j_{\Es}(f_0^{\gns_0, i})(\ga^{\gns_0, i}_0)$
	for all $i < i_0$.
	We choose $f''$ such that 
		$j_\Es(f'')(\gb'') \not\in j_\Es(R)(\gk^0(\gns_0), \gm)$
	for each inaccessible $\gm < i_0$.
    \begin{align*}
	    & q'' = u^{\gns_0, i}_0,
	    \\
	    & \ga'' \Egt \ga^{\gns_0, i}, \gb''.
    \end{align*}
\item 
	$i_0$ is limit: 
	If $\ordof{j_{\Es}(f_0^{\gns_0, i})(\ga_0^{\gns_0, i})} 
				{i < i_0}$
	is a maximal anti-chain below 
		$j_\Es(F^{p_0})(\gp_{\ga, \gb}(\gk(\gns_0)), \gb)$
	then we finish the induction on $i$.
	Otherwise we choose $f''$, $\gb''$ such that
	$j_\Es(f'')(\gb'') \incompatible 
		j_{\Es}(f_0^{\gns_0, i})(\ga_0^{\gns_0, i})$
	for all $i < i_0$.
	We choose $f''$ such that 
		$j_\Es(f'')(\gb'') \not\in j_\Es(R)(\gk^0(\gns_0), \gm)$
	for each inaccessible $\gm < i_0$.
	Choose $\ga'' \in N$ such that
    $\forall i < i_0$ $\gas'' \Egt \gas^{\gns_0, i}, \gb''$ 
    and set
    \begin{align*}
	    q'' = \bigunion_{i < i_0} u^{\gns_0, i}_0.
    \end{align*}
\end{itemize}
Set
    \begin{align*}
	&T''_0 =   \gp^{-1}_{\ga'',\gb} \,
			T^{p_0} \setminus \gp_{\ga, \ga''}(\gns_0),
	\\
	&F''_{0} =
	   F^{p_0}  \circ \gp_{\ga'', \gb},
	\\
	& u''_0 = 
		(q''_{\ordered{\gp_{\ga, \ga''}(\gns_0)}})_0
		 \union \set{\ordered{
			\ga'', 
			T''_0, 
			f'', 
			F''_0
		}}),
\end{align*}
\begin{multline*}
       	u''_1 = 
		(q''_{\ordered{\gp_{\ga, \ga''}(\gns_0)}})_1
		 \union \set{\ordered{
			\gp_{\ga, \ga''}(\gns_0),
			T''_0 \restricted \gp_{\ga, \ga''}(\gns_0), 
	\\
			f^{p_0} \circ \gp_{\ga'', \gb} 
				\restricted \gp_{\ga, \ga''}(\gns_0),
			F''_0 \restricted \gp_{\ga, \ga''}(\gns_0)
		}}.
\end{multline*}
If there is $q_1 \append q_0 \in D \intersect N$ such that
\begin{align*}
q_1 \append q_0 \leq^* u''_1 \append u''_0
\end{align*}
then set
\begin{align*}
& \ga_0^{\gns_0, i_0} = \gk(\mc(q_0)),
\\
& u^{\gns_0, i_0}_0 = q'' \union
	\setof{\ordered{\ggs, q_0^{\ggs}}}
		{\ggs \in \supp q_0 \setminus \supp q''},
\\
& T^{\gns_0, i_0}_0 = T^{q_0},
\\
& f^{\gns_0, i_0}_0 = f^{q_0},
\\
& F^{\gns_0, i_0}_0 = F^{q_0},
\\
& \gas_1^{\gns_0, i_0} = \mc (q_1),
\\
& u^{\gns_0, i_0}_1 = 
	\setof{\ordered{\ggs, q_1^{\ggs}}}
		{\ggs \in \supp q_1},
\\
& T^{\gns_0, i_0}_1 = T^{q_1},
\\
& f^{\gns_0, i_0}_1 = f^{q_1},
\\
& F^{\gns_0, i_0}_1 = F^{q_1},
\end{align*}
otherwise we set
\begin{align*}
& \ga_0^{\gns_0, i_0} = \ga'',
\\
& u^{\gns_0, i_0}_0 = q'',
\\
& T^{\gns_0, i_0}_0 = T''_0,
\\
& f^{\gns_0, i_0}_0 = f''_0,
\\
& F^{\gns_0, i_0}_0 = F''_0,
\\
& \gas_1^{\gns_0, i_0} = \mc(u''_1),
\\
& u^{\gns_0, i_0}_1 = 
	\setof{\ordered{\ggs, u_1^{\prime \prime \ggs}}}
		{\ggs \in \supp u''_1},
\\
& T^{\gns_0, i_0}_1 = T^{u''_1},
\\
& f^{\gns_0, i_0}_1 = f^{u''_1},
\\
& F^{\gns_0, i_0}_1 = F^{u''_1}.
\end{align*}
When the induction on $i$ terminates we have
\begin{align*}
\ordof {\ga_0^{\gns_0, i}, u^{\gns_0, i}_0, T^{\gns_0, i}_0,
                 f_0^{\gns_0, i}, F_0^{\gns_0, i},
	\gas^{\gns_0, i}_1, u^{\gns_0, i}_1, T^{\gns_0, i}_1,
                 f_1^{\gns_0, i}, F_1^{\gns_0, i} }
       {i < \gk}.
\end{align*}
We complete step $\gns_0$ by setting
\begin{align*}
&  \ga^{\gns_0}_0 \in N,
\\
& \forall i < i^{\gns_0} \  \ga^{\gns_0}_0 \Egt \ga^{\gns_0, i}_0,
\\
& u^{\gns_0}_0 = \bigunion_{i < \gk}
		u^{\gns_0, i}_0,
\\
& T^{\gns_0}_0 = \dsintersect_{i < \gk}
		\gp^{-1}_{\ga^{\gns_0}_0, \ga^{\gns_0, i}_0} \, 
		T^{\gns_0, i}_0,
\\
& \forall i < \gk\ 
	F_0^{\gns_0}\leq F_0^{\gns_0, i}
	\circ
		\gp_{\ga^{\gns_0}, \ga^{\gns_0, i}}.
\end{align*}
When the induction on $\gns$ terminates we have
\begin{align*}
\ordof{\ga^\gns_0, u^\gns_0, T^\gns_0, F^\gns_0}
      {\gns \in A}.
\end{align*}
We define the following function with domain $A$:
\begin{align*}
g(\gns)= \ordof{j_\Es(f^{\gns, i}_0)(\ga_0^{\gns, i})} 
			{i < \gk}.
\end{align*}
By the construction $g(\gns) \in \ME$ is a maximal anti-chain 
in $j_\Es(R)(\gk^0(\gns), \gk)$
below 
$j_\Es(F^{p_0})(\gp_{\ga, \gb}(\gk(\gns)), \gb)$.
Hence, $\forall \gt < \len(\Es)$ 
	$j_\Es(g)(\Es_\ga \restricted \gt) \in M^2_\Es$ 
is a maximal anti-chain below
$j^2_\Es(F^{p_0})(\gb, j_\Es(\gb))$.
By genericity of $I(\Es)$ over $M^2_\Es$, there are $g^\gt$'s
such that
$j^2_\Es(g^\gt)(\ga, j_\Es(\ga)) \in I(\Es)$ and
$j^2_\Es(g^\gt)(\ga, j_\Es(\ga))$ is stronger than a condition in
$j_\Es(g)(\Es_\ga \restricted \gt)$.
Let $h^\gt$ be such that $j^2_\Es(g^\gt)(\ga, j_\Es(\ga)) \leq
j_\Es(g)(\Es_\ga \restricted \gt)(j_\Es(h^\gt)(\Es_\ga \restricted \gt))$.
Note that we can use $\ga$ here because the generic was built through
the normal ultrafilter. If we would not have had this property we
would have enlarged $\ga$ to accommodate the intersection and we might
have needed different $\ga$ for each $\gt$.
We combine the information gathered into $\gt+1$ conditions.
\begin{align*}
& T^{p'_0} = \dsintersect_{\gns \in A} \gp^{-1}_{\ga, \ga^\gns} \, 
			T^{\gns}_0 ,
\\
& \forall \gn \in A \ \len(\gn)=0 \implies f^{p'_0}(\gn) = 
			f^{p_0} \circ \gp_{\ga, \gb} (\gn),
\\
& \forall \gns \in A \ F^{p'_0} \leq F_0^\gns \circ 
			\gp_{\ga, \ga^\gns},
\\
& p'_0 = \bigunion_{\gns \in A} u^{\gns}_0.
\end{align*}
For each $\gt < \len(\Es)$ we define
\begin{align*}
& \gas^{\gt}_1 = j_\Es(\gas_1)^{\Es_\ga \restricted \gt, 
		j_\Es(h^\gt)(\Es_\ga \restricted \gt)},
\\
& u^{\prime\gt}_1 = j_\Es(u_1)^{\Es_\ga \restricted \gt, 
		j_\Es(h^\gt)(\Es_\ga \restricted \gt)},
\\
& u^\gt_1 = \setof {\ordered{\Es_\gga, 
			(u^{\prime\gt}_1)^{\Es_\gga\restricted \gt} }} 
		{\Es_\gga \restricted \gt \in \supp u^{\prime\gt}_1 },
\\
& T^{\gt}_1 = j_\Es(T_1)^{\Es_\ga \restricted \gt, 
		j_\Es(h^\gt)(\Es_\ga \restricted \gt)},
\\
& \text{If } \gt > 0 \text{ then }
  f^{\gt}_1 = j_\Es(f_1)^{\Es_\ga \restricted \gt, 
		j_\Es(h^\gt)(\Es_\ga \restricted \gt)},
\\
& \text{If } \gt = 0 \text{ then }
 \forall \gns \in A\ \len(\gns)=0 \implies 
			f^0_1(\gns) = f^{\gns, h^0(\gns)},
\\
& F^{\gt}_1 = j_\Es(F_1)^{\Es_\ga \restricted \gt, 
		j_\Es(h^\gt)(\Es_\ga \restricted \gt)}.
\end{align*}
We note that the construction ensures us that 
$f^\gt_1 \leq f^{p_0}$, $F^\gt_1 \leq F^{p_0}$ when $\gt > 0$ and
$u^\gt_1$ and $p'_0$ do not
contain contradictory information. Hence we can define for each
$\gt < \len(\Es)$ the following
\begin{align*}
& \gb^\gt \Egt \ga_1^\gt, \ga,
\\
& T^{\gt}_0  = \gp^{-1}_{\gb^\gt, \ga} T^{p'_0},
\\
& f^{\gt}_0 = f^\gt_1 \circ \gp_{\gb^\gt, \ga^\gt_1},
\\
& F^{\gt}_0  \leq F^\gt_1 \circ \gp_{\gb^\gt, \ga^\gt_1}, 
		  g^\gt \circ \gp_{\gb^\gt, \ga},
\\
& p^\gt_0 = p'_0 \union u^\gt_1 \union
	 \set{\ordered{
		\Es_{\gb^\gt}, T^\gt_0, f^\gt_0, F^\gt_0
		}}.
\end{align*}
We consider the following sets
\begin{align*}
& A^\gt = \setof {\gns}
		{ \exists t'_1 \leq^{**}
			(p^\gt_{0 {\ordered{\gns}}})_1\ 
			t'_1 \append (p^\gt_{0 {\ordered{\gns}}})_0 \in D
		}.
\end{align*}
\ifnum\article=0
\enlargethispage*{10pt}
\fi
There are $2$ possibilities at this point
\begin{enumerate}
\item
	There is $\gt < \len(\Es_\ga)$ such that 
	$A^\gt \in E_{\gb^\gt}(\gt)$:
	Of course, $p^\gt_0$ satisfies the requested conclusion.
	Hence we set $p^*_0 = p^\gt_0$ and the theorem is proved.
\item
	$\forall \gt < \len(\Es_\ga)$ $A^\gt \not\in \Es_{\gb^\gt}(\gt)$:
We claim that some shrinkage of $T^{p'_0}$ is enough to get
us into clause 2 of the theorem. Let us assume, by contradiction, 
that a small shrinkage
can not bring us to $p^*_0$.
This means that there is $\gt < \len(\Es_\ga)$ such that
\begin{align*}
\setof {\gns \in T^{p'_0}}
	{
		 \exists q_1 \append q_0  \leq^* p'_{0 \ordered{\gns}}\ 
			 q_1 \append q_0  \in D
	}
    \in E_\ga(\gt).
\end{align*}
Let $F'' \leq F^{p'_0}, g^\gt$ and $p''_0$ be $p'_0$ with $F''$ substituted
for $F^{p'_0}$. Due to openness of $D$ we still have
\begin{align*}
\setof {\gns \in T^{p'_0}}
	{
		 \exists q_1 \append q_0 \leq^* p''_{0 \ordered{\gns}}\ 
			 q_1 \append q_0  \in D
	}
    \in E_\ga(\gt).
\end{align*}
Hence, by the construction, we have
\begin{multline*}
\setof {\gns \in T^{p'_0}}
	{ 
	u_1^{\gns, h^\gt(\gns)} \union 
	\set{\ordered{
		\gas_1^{\gns, h^\gt(\gns)}, 
		T^{{\gns, h^\gt(\gns)}}_1, 
		f^{{\gns, h^\gt(\gns)}}_1, 
		F^{{\gns, h^\gt(\gns)}}_1
	}}
		\append
\\
	(u_{0 \ordered{\gp_{\ga, \ga_0^{\gns, h^\gt(\gns)}}(\gns)}}^
							{\gns, h^\gt(\gns)})_0
	\union \set{\ordered{
		\Es_{\ga_0^{\gns, h^\gt(\gns)}},
		 T^{\gns, h^\gt(\gns)}_0, 
		 f_0^{\gns, h^\gt(\gns)}, 
		F^{\gns, h^\gt(\gns)}_0
	}} \in D
	}
	    \in E_\ga(\gt).
\end{multline*}
Invoking $\gp^{-1}_{\gb^\gt, \ga}$ on the above set yields
\begin{multline*}
\setof {\gns \in T^{\gt}_0}
	{ 
	u_1^{\gp_{\gb^\gt, \ga}(\gns), h^\gt(\gp_{\gb^\gt, \ga}(\gns))} \union 
	\set{\ordered{
		\gas_1^{\gp_{\gb^\gt, \ga}(\gns), h^\gt(\gp_{\gb^\gt, \ga}(\gns))}, 
		T^{{\gp_{\gb^\gt, \ga}(\gns), h^\gt(\gp_{\gb^\gt, \ga}(\gns))}}_1, 
\\
		f^{{\gp_{\gb^\gt, \ga}(\gns), h^\gt(\gp_{\gb^\gt, \ga}(\gns))}}_1, 
		F^{{\gp_{\gb^\gt, \ga}(\gns), h^\gt(\gp_{\gb^\gt, \ga}(\gns))}}_1
	}}
		\append
\\
	(u_{0 \ordered{\gp_{\gb^\gt, \ga_0^{\gp_{\gb^\gt, \ga}(\gns), h^\gt(\gp_{\gb^\gt, \ga}(\gns))}}(\gns)}}^
							{\gp_{\gb^\gt,\ga}(\gns), h^\gt(\gp_{\gb^\gt, \ga}(\gns))})_0
	\union \set{\ordered{
		\Es_{\ga_0^{\gp_{\gb^\gt, \ga}(\gns), h^\gt(\gp_{\gb^\gt, \ga}(\gns))}},
		 T^{\gp_{\gb^\gt, \ga}(\gns), h^\gt(\gp_{\gb^\gt, \ga}(\gns))}_0, 
\\
		 f_0^{\gp_{\gb^\gt, \ga}(\gns), h^\gt(\gp_{\gb^\gt, \ga}(\gns))}, 
		F^{\gp_{\gb^\gt, \ga}(\gns), h^\gt(\gp_{\gb^\gt, \ga}(\gns))}_0
	}} \in D
	}
	    \in E_{\gb^\gt}(\gt).
\end{multline*}
Now, from the construction of $p^\gt$ we see that
\begin{multline*}
\setof {\gns \in T^{\gt}_0}
	{ 
	\exists t'_1 \leq^{**}
		(p^\gt_{\ordered{\gns}})_1\ 
	t'_1 \append (p^\gt_{\ordered{\gns}})_0 \leq^*
	u_1^{\gp_{\gb^\gt, \ga}(\gns), h^\gt(\gp_{\gb^\gt, \ga}(\gns))} \union 
\\
	\set{\ordered{
		\gas_1^{\gp_{\gb^\gt, \ga}(\gns), h^\gt(\gp_{\gb^\gt, \ga}(\gns))}, 
		T^{{\gp_{\gb^\gt, \ga}(\gns), h^\gt(\gp_{\gb^\gt, \ga}(\gns))}}_1, 
\\
		f^{{\gp_{\gb^\gt, \ga}(\gns), h^\gt(\gp_{\gb^\gt, \ga}(\gns))}}_1, 
		F^{{\gp_{\gb^\gt, \ga}(\gns), h^\gt(\gp_{\gb^\gt, \ga}(\gns))}}_1
	}}
		\append
\\
	(u_{0 \ordered{\gp_{\ga, \ga_0^{\gp_{\gb^\gt, \ga}(\gns), h^\gt(\gp_{\gb^\gt, \ga}(\gns))}}}}^
							{\gp_{\gb^\gt, \ga}(\gns), h^\gt(\gp_{\gb^\gt, \ga}(\gns))})_0
	\union 
\\
	\set{\ordered{
		\Es_{\ga_0^{\gp_{\gb^\gt, \ga}(\gns), h^\gt(\gp_{\gb^\gt, \ga}(\gns))}},
		 T^{\gp_{\gb^\gt, \ga}(\gns), h^\gt(\gp_{\gb^\gt, \ga}(\gns))}_0, 
\\
		 f_0^{\gp_{\gb^\gt, \ga}(\gns), h^\gt(\gp_{\gb^\gt, \ga}(\gns))}, 
		F^{\gp_{\gb^\gt, \ga}(\gns), h^\gt(\gp_{\gb^\gt, \ga}(\gns))}_0
	}}
	}
	    \in E_{\gb^\gt}(\gt).
\end{multline*}
Combining the above $2$ formulas and recalling that $D$ is open we get
\begin{align*}
\setof {\gns \in T^{\gt}_0}
	{ 
	\exists t'_1 \leq^{**}
		(p^\gt_{\ordered{\gns}})_1\ 
	t'_1 \append (p^\gt_{\ordered{\gns}})_0 \in D
	}
	    \in E_{\gb^\gt}(\gt).
\end{align*}
That is $A^\gt \in E_{\gb^\gt}(\gt)$. Contradiction.
So, we have shown that
\begin{align*}
& T^{p^*_0} =
	 \setof {\gns \in T^{p'_0}}
	{
		 \forall q_1 \append q_0  \leq^* p'_{0 \ordered{\gns}}\ 
			 q_1 \append q_0  \not \in D
	}
    \in \Es_\ga.
\end{align*}
By letting $p^*_0$ be $p'_0$ with $T^{p'_0}$ substituted by $T^{p^*_0}$
we get clause $2$.
\end{enumerate}
\end{proof}
\begin{lemma} \label{DenseHomogen-0}
Let $D$ be dense open in $\PE/\Pe$, $p = p_0 \in \PE/\Pe$.
Then there are $n<\gw$, $p^*_0 \leq^* p_0$ and
        $S \subseteq [T^{p^*_0}]^n$, an $\mc(p^*_0)$-fat tree,
such that 
\begin{multline*}
\forall \ordered{\gns_1,\dotsc,\gns_{n}} \in S\ 
\exists t'_n \append \dotsb \append t'_1 \leq^{**}
                (p^*_{0 \ordered{\gns_1,\dotsc,\gns_{n}}})_{n..1}
\\
        t'_n \append \dotsb \append t'_1 \append
		(p^*_{0 \ordered{\gns_1,\dotsc,\gns_{n}}})_0 \in D.
\end{multline*}
\end{lemma}
\begin{proof}
Construct, by repeat invocation of \ref{AlmostDenseHomogen-0} for each $n <\gw$,
a $\leq^*$-decreasing sequence $\ordof{p^n_0} {n <\gw}$.
Let $p^*_0 \leq^* p^n_0$ for all $n < \gw$.
Choose $q \in D$ such that $q \leq p^*_0$.
There is $\ordered{\gns_1, \dotsc, \gns_n} \in [T^{p^*_0}]^n$ such that
$q \leq^* p^*_{0 \ordered{\gns_1, \dotsc, \gns_n}}$.
By this we eliminated clause 2 of claim \ref{AlmostDenseHomogen-0} for $n$.
\end{proof}
\begin{claim} \label{DenseHomogenTriv-0}
Assume $\len(\ges) = 0$ and 
let $D$ be dense open in $\Pe/\Peii$, $p = p_0 \in \Pe/\Peii$.
Then there is $p^*_0 \leq^* p_0$ such that $p^*_0 \in D$.
\end{claim}
\begin{proof}
Of course this is completely trivial  as 
		$(\Pe/\Peii)/p_0 \simeq R(\gk(p^0_0), \gk^0(\ges))/p_0$.
In this case $\leq^*$ and $\leq$ are the same.
\end{proof}
\begin{lemma} \label{PeNoChange}
Assume $\len(\ges) = 0$ and let G be $\Pe/\Peii$ with $p = p_0 \in G$.
Then  $V_{\gk^0(\ges_2)+1} = V_{\gk^0(\ges_2)+1}^{V[G]}$ and
	$\Peii = P^{V[G]}_{\ges_2}$.
\end{lemma}
\begin{proof}
Let $\gn = \gk(p_0^0)$. As $\Pe/\Peii$ is $\gn^+$-closed we get
immediately that $V_{\gk^0(\ges_2)+1} = V_{\gk^0(\ges_2)+1}^{V[G]}$.

We have much more than that. Namely, $\gn^{++}, \dotsc, \gn^{+6}$ are
not collapsed. For this
we remind the reader that the $V$ we work with is a generic
extension of $V^*$ for a reverse Easton forcing.
Let $Q_1$ be the reverse Easton forcing up to $\gk^0(\ges_2)$ and
$H_1$ be its generic.
Let $Q_2$ be the forcing at stage $\gk^0(\ges_2)$
and $H_2$ be its generic over $V^*[H_1]$.
Let $Q_3$ be the rest of the reverse Easton forcing up to $\gk^0(\ges)$
and $H_3$ be its generic over $V^*[H_1][H_2]$.
Then we have
\begin{align*}
\Pset^{V[G]}(\gk^0(\ges)) = 
	 & \Pset^{V^*[H_1][H_2][H_3][G]}(\gk^0(\ges)).
\end{align*}
Note that $V^*[H_1][H_2][H_3][G]$ is a reflection of the
situation at \ref{prototype} and by this we see that 
$\gn^{++}, \dotsc, \gn^{+6}$ are not collapsed.

In fact we see that nothing has changed as far as the definition of
$\Peii$ in $V[G]$ is concerned. (We might have new anti-chains which
is no obstacle to us).
\end{proof}
\begin{lemma} \label{DenseHomogenTriv-l+1}
Assume $\len(\ges) = 0$ and let $p = p_{l..0} \in \Pe$.
Assume that \ref{FullHomogen} is true for $p_{l..1} \in \Peii$ and dense open
subsets of $\Peii$. Then it is true for $p_{l..0}$ and dense open subsets
of $\Pe$.
\end{lemma}
\begin{proof}
In order to avoid excess of indices we  give the proof of the
case $p = p_{1..0}$.

Let $G$ be $\Pe/\Peii$-generic with $p_0 \in G$
and let $D$ be dense open in $\Pe$.
Then $D_{\ges_2} = \setof {q \leq p_1} {q\append r \in D,\ r \in G} \in V[G]$
 is a dense
open subset of $\Peii$. By \ref{PeNoChange} and \ref{FullHomogen} for $p_1$
there are $p^*_1 \in V$ and $S^1 \in V$,
an $\mc(p^*_1)$-fat tree, such that $p^*_1 \leq^* p_1$ and
\begin{multline*}
\forall \ordered{\gns_1, \dotsc, \gns_n} \in S^1\ 
\exists t'_n \append \dotsb \append t'_1 \leq^{**}
                (p^*_{1 \ordered{\gns_1,\dotsc,\gns_{n}}})_{n..1}
\\
        t'_n \append \dotsb \append t'_1 \append
		(p^*_{1 \ordered{\gns_1,\dotsc,\gns_{n}}})_0 \in D_{\ges_2}.
\end{multline*}
Hence there is $p'_0 \leq^* p_0$ which forces the above.
That is for each $\ordered{\gns_1, \dotsc, \gns_n} \in S^1$
there is a maximal anti-chain, $A(\gns_1, \dotsc, \gns_n)$, of 
	$\Pe/\Peii$ below $p'_0$ such that
\begin{multline*}
\forall \ordered{\gns_1, \dotsc, \gns_n} \in S^1\ 
\exists t'_n \append \dotsb \append t'_1 \leq^{**}
                (p^*_{1 \ordered{\gns_1,\dotsc,\gns_{n}}})_{n..1}
\\
\forall q_0 \in A(\gns_1, \dotsc, \gns_n)
\\
        t'_n \append \dotsb \append t'_1 \append
		(p^*_{1 \ordered{\gns_1,\dotsc,\gns_{n}}})_0 
		\append
		q_0 \in D.
\end{multline*}
By noting that $\ordered{\Pe/\Peii, \leq^*}$ is $\gk^0(\ges_2)^+$-closed
and that $\power{S^1} = \gk^0(\ges_2)$ we see that there
is $p^*_0 \leq^* p'_0$ such that
\begin{multline*}
\forall \ordered{\gns_1, \dotsc, \gns_n} \in S^1\ 
\exists t'_n \append \dotsb \append t'_1 \leq^{**}
                (p^*_{1 \ordered{\gns_1,\dotsc,\gns_{n}}})_{n..1}
\\
        t'_n \append \dotsb \append t'_1 \append
		(p^*_{1 \ordered{\gns_1,\dotsc,\gns_{n}}})_0 
		\append
		p^*_0 \in D.
\end{multline*}
\end{proof}
Our main obstacle in proving \ref{FullHomogen} is that in general
$\PE/\Pe$ is $\gk^0(\ges)^+$-closed while $\Pe$ is $\gk^0(\ges)^{++}$-c.c.
However, when $\len(\ges)=0$ we have $\Pe$ is $\gk^0(\ges)^+$-c.c.
The following $2$ lemmas give us facts in this case which  help
us to overcome the obstacle.
\begin{lemma} \label{DenseAlmostHomogen-0l}
Let $D$ be dense open in $\PE$, $p = p_l \append \dotsb \append p_0 \in \PE$.
Then there are $n<\gw$, $p^*_0 \leq^* p_0$, $q \leq p_{l..1}$ and
        $S \subseteq [T^{p^*_0}]^n$, an $\mc(p^*_0)$-fat tree,
such that 
\begin{multline*}
\forall \ordered{\gns_1,\dotsc,\gns_{n}} \in S\ 
\exists t'_n \append \dotsb \append t'_1 \leq^{**}
                (p^*_{0 \ordered{\gns_1,\dotsc,\gns_{n}}})_{n..1}
\\
        q \append t'_n \append \dotsb \append t'_1 \append
		(p^*_{0 \ordered{\gns_1,\dotsc,\gns_{n}}})_0 \in D.
\end{multline*}
\end{lemma}
\begin{proof}
Let $D_\Es = \setof {r \leq p_0} {\exists s \leq p_{l..1 }\ 
			s \append r \in D}$.
Then $D_\Es$ is dense open in $\PE/\Pe$ below $p_0$ where $p_{l..1} \in \Pe$.
By \ref{DenseHomogen-0} there are $n < \gw$, $p^*_0 \leq^* p_0$,
$S \subseteq [T^{p^*_0}]^n$ such that
\begin{multline*}
\forall \ordered{\gns_1,\dotsc,\gns_{n}} \in S\ 
\exists t'_n \append \dotsb \append t'_1 \leq^{**}
                (p^*_{0 \ordered{\gns_1,\dotsc,\gns_{n}}})_{n..1}
\\
        t'_n \append \dotsb \append t'_1 \append
		(p^*_{0 \ordered{\gns_1,\dotsc,\gns_{n}}})_0 \in D_\Es.
\end{multline*}
By the definition of $D_\Es$ we see that there is a function
$q(\gns_1,\dotsc,\gns_{n})$  with domain $S$ such that
\begin{multline*}
\forall \ordered{\gns_1,\dotsc,\gns_{n}} \in S\ 
\exists t'_n \append \dotsb \append t'_1 \leq^{**}
                (p^*_{0 \ordered{\gns_1,\dotsc,\gns_{n}}})_{n..1}
\\
        q(\gns_1, \dotsc, \gns_n) \append t'_n \append \dotsb \append 
			t'_1 \append
		(p^*_{0 \ordered{\gns_1,\dotsc,\gns_{n}}})_0 \in D.
\end{multline*}
As $q(\gns_1,\dotsc,\gns_{n}) \leq p_{l..1}$ there is $q \leq p_{l..1}$
such that except on a measure $0$ set we have $q(\gns_1,\dotsc,\gns_{n}) = q$.
By removing this measure $0$ set from $S$ we get
\begin{multline*}
\forall \ordered{\gns_1,\dotsc,\gns_{n}} \in S\ 
\exists t'_n \append \dotsb \append t'_1 \leq^{**}
                (p^*_{0 \ordered{\gns_1,\dotsc,\gns_{n}}})_{n..1}
\\
        q \append t'_n \append \dotsb \append 
			t'_1 \append
		(p^*_{0 \ordered{\gns_1,\dotsc,\gns_{n}}})_0 \in D.
\end{multline*}
\end{proof}
\begin{lemma} \label{DenseAntiHomogen-0l}
Let $D$ be dense open in $\PE$, $p = p_{l..0} \in \PE$.
Let $\ges$ be such that $p_{l..1} \in \Pe$ where $\len(\ges) = 0$.
Then there are $p^*_0 \leq^* p_0$ and
	$\ordof{q^\gx} {\gx < \gk^0(\ges)}$ a maximal anti-chain below 
	$p_{l..1}$ such that
for each $\gx < \gk^0(\ges)$ there are $p'_0 \geq^* p^*_0$,
$n<\gw$ and
        $S \subseteq [T^{p'_0}]^n$, an $\mc(p'_0)$-fat tree,
such that 
\begin{multline*}
\forall \ordered{\gns_1,\dotsc,\gns_{n}} \in S\ 
\exists t'_n \append \dotsb \append t'_1 \leq^{**}
                (p'_{0 \ordered{\gns_1,\dotsc,\gns_{n}}})_{n..1}
\\
        q^\gx \append t'_n \append \dotsb \append t'_1 \append
		(p'_{0 \ordered{\gns_1,\dotsc,\gns_{n}}})_0 \in D
\end{multline*}
and
\begin{align*}
\setof{
        q^\gx \append t'_n \append \dotsb \append t'_1 \append
		(p'_{0 \ordered{\gns_1,\dotsc,\gns_{n}}})_0
	} 
	{
	\ordered{\gns_1,\dotsc,\gns_{n}} \in S
	}
\end{align*}
is pre-dense below $q^\gx \append p'_0$.

Of course an immediate corollary is for each $\gx < \gk^0(\ges)$ 
there are $n<\gw$ and
        $S \subseteq [T^{p^*_0}]^n$, an $\mc(p^*_0)$-fat tree,
such that 
\begin{multline*}
\forall \ordered{\gns_1,\dotsc,\gns_{n}} \in S\ 
\exists t'_n \append \dotsb \append t'_1 \leq^{**}
                (p^*_{0 \ordered{\gns_1,\dotsc,\gns_{n}}})_{n..1}
\\
        q^\gx \append t'_n \append \dotsb \append t'_1 \append
		(p^*_{0 \ordered{\gns_1,\dotsc,\gns_{n}}})_0 \in D.
\end{multline*}
\end{lemma}
\begin{proof}
Our first observation is that $\Pe$ is $\gk^0(\ges)^+$-c.c.
(As opposed to the usual $\gk^0(\ges)^{++}$-c.c. we have when $\len(\ges) > 0$).
We  construct, by induction,
the sequence $\ordof{q^\gx} {\gx < \gx_0}$ where $\gx_0 < \gk^0(\ges)^+$.
Together with it we  construct an auxiliary $\leq^*$-decreasing
sequence $\ordof{p^\gx_0} {\gx < \gx_0}$.
\begin{itemize}
\item
	$\gx_0 = 0$: By \ref{DenseAlmostHomogen-0l}, \ref{ToPreDense-n} 
		and openness of $D$ there are
	$q^0 \leq p_{l..1}$, $p^0_0 \leq^* p_0$, $S \subseteq [T^{p^0_0}]^n$
	such that
	\begin{multline*}
	\forall \ordered{\gns_1,\dotsc,\gns_{n}} \in S\ 
	\exists t'_n \append \dotsb \append t'_1 \leq^{**}
                (p^0_{0 \ordered{\gns_1,\dotsc,\gns_{n}}})_{n..1}
	\\
        	q^0 \append t'_n \append \dotsb \append t'_1 \append
		(p^0_{0 \ordered{\gns_1,\dotsc,\gns_{n}}})_0 \in D
	\end{multline*}
	and $\setof{
	        q^0 \append t'_n \append \dotsb \append t'_1 \append
			(p^0_{0 \ordered{\gns_1,\dotsc,\gns_{n}}})_0
		} 
		{
		\ordered{\gns_1,\dotsc,\gns_{n}} \in S
		}$
	is pre-dense below $q^0 \append p^0_0$.
\item
	$\gx_0 > 0$: If $\ordof {q^\gx} {\gx < \gx_0}$ is a maximal
	anti-chain below $p_{l..1}$ then the induction is finished.
	If it is not a maximal anti-chain we observe that
	$\gx_0 < \gk^0(\ges)^+$ as $\Pe$ is $\gk^0(\ges)^{+}$-c.c.
	Let $q' < p_{l..1}$ be such that
	$\forall \gx < \gx_0$ $q' \incompatible q^\gx$.
	As $\ordered{\PE, \leq^*}$ is $\gk^0(\ges)^+$-closed
	there is $p'_0$ such that $p'_0 \leq^* p^{\gx}_0$ for all $\gx < \gx_0$.
	By \ref{DenseAlmostHomogen-0l}, \ref{ToPreDense-n} starting from
	$q' \append p'_0$
	there are
	$q^{\gx_0} \leq q'$, $p^{\gx_0}_0 \leq^* p'_0$, 
		$S \subseteq [T^{p^{\gx_0}_0}]^n$
	such that
	\begin{multline*}
	\forall \ordered{\gns_1,\dotsc,\gns_{n}} \in S\ 
	\exists t'_n \append \dotsb \append t'_1 \leq^{**}
                (p^{\gx_0}_{0 \ordered{\gns_1,\dotsc,\gns_{n}}})_{n..1}
	\\
        	q^{\gx_0} \append t'_n \append \dotsb \append t'_1 \append
		(p^{\gx_0}_{0 \ordered{\gns_1,\dotsc,\gns_{n}}})_0 \in D
	\end{multline*}
	and $\setof{
	        q^{\gx_0} \append t'_n \append \dotsb \append t'_1 \append
			(p^{\gx_0}_{0 \ordered{\gns_1,\dotsc,\gns_{n}}})_0
		} 
		{
		\ordered{\gns_1,\dotsc,\gns_{n}} \in S
		}$
	is pre-dense below $q^{\gx_0} \append p^{\gx_0}_0$.
\end{itemize}
When the induction terminates we have a $\leq^*$-decreasing sequence
$\ordof {p^\gx_0} {\gx < \gx_0}$ where $\gx_0 < \gk^0(\ges)^+$.
By choosing $p^*_0 \leq^* p^\gx_0$ for all $\gx < \gx_0$ we finish 
the proof.
\end{proof}
\begin{lemma} \label{DenseHomogen-l+1}
Let $p = p_{l..0} \in \PE$.
Assume that \ref{FullHomogen} is true for $p_{l..1} \in \Pe$ and dense open
subsets of $\Pe$. Then it is true for $p_{l..0}$ and dense open subsets
of $\PE$.
\end{lemma}
\begin{proof}
In order to avoid too many indices we  prove the lemma for the
case $p = p_{1..0}$.

Choose $\gc$ large enough so that $H_\gc$  contains everything we are
interested in. Let $N \subelem H_\gc$ be such that $\power{N} = \gk$,
	$N \supset N^{\upto \gk}$, $p \in N$, $\PE \in N$.

Let $\Es_\gb = \mc(p_0)$.
Choose $\ga \in \dom \Es$ such that $\ga \Egt \gga$ for all 
$\gga \in \dom \Es \intersect N$.
Let $A = \setof {\gns \in \gp^{-1}_{\ga, \mc(p_0)} T^{p_0}} {\len(\gns) = 0}$. 
Note that $A \in E_\ga(0)$.
Let $\preceq$ be a well ordering of $A$ such that
$\forall \gns_1,\gns_2 \in A \ 
	\gns_1 \preceq \gns_2 \implies \gk^0(\gns_1) \leq \gk^0(\gns_2)$.
We shrink $A$ a bit so that the following is satisfied:
$\forall \gns \in A\ 
	\power{\setof{\gms \in A} {\gk^0(\gms) < \gk^0(\gns)}} 
		\leq \gk^0(\gns)$.
We start an induction on $\gns$ in which we build
\begin{align*}
\ordof{\ga^\gns, u^\gns_0, T^\gns_0, F^\gns_0}
      {\gns \in A}.
\end{align*}
Assume that we have constructed
$\ordof{\ga^\gns, u^\gns_0, F^\gns_0, T^\gns_0}
      {\gns \prec \gns_0}$.
We start working in $N$.
Set the following:
\begin{itemize}

\item $\gns_0$ is $\prec$-minimal:
    \begin{align*}
	& q' = p_0 \setminus \set{\ordered{\Es_\gb, T^{p_0}, f^{p_0},
					F^{p_0}}},
	\\
	& \ga' = \gb.
    \end{align*}

\item  $\gns_0$ is the immediate $\prec$-successor of $\gns$:
    \begin{align*}
	    & q' = u^\gns_0,
	    \\
	    & \ga' = \ga^\gns.
    \end{align*}
\item $\gns_0$ is $\prec$-limit: Choose $\ga' \in N$ such that
    $\forall \gns \prec \gns_0$ $\ga' \Egt \ga^\gns$ 
    and set
    \begin{align*}
	    q' = \bigunion_{\gns \prec \gns_0} u^\gns_0.
    \end{align*}
\end{itemize}
We make an induction on $i$ which  builds
$\ordof {\ga^{\gns_0, i}, u^{\gns_0, i}_0, T^{\gns_0, i}_0,
                 f^{\gns_0, i}_0,
                 F^{\gns_0, i}_0}
	       {i < \gk}$.
Assume  we have constructed
\begin{align*}
\ordof {\ga^{\gns_0, i}, u^{\gns_0, i}_0, T^{\gns_0, i}_0,
                 f^{\gns_0, i}_0, F^{\gns_0, i}_0}
	       {i < i_0},
\end{align*}
and we do step $i_0$.
\begin{itemize}

\item $i_0 = 0$:
    \begin{align*}
	& q'' = q',
	\\
	& \ga'' = \ga',
	\\
	& f'' = F^{p_0}(\gk(\gp_{\ga, \gb}(\gns_0)), \gp_{\ga'',\gb}(-)).
    \end{align*}

\item $i_0 = i + 1$: If 
	$\ordof{j_\Es(f^{\gns_0, i}_0)(\ga^{\gns_0, i})} {i < i_0}$
	is a maximal anti-chain below
		$j_\Es(F^{p_0})(\gk(\gp_{\ga, \gb}(\gns_0)), \gb)$
	we terminate the induction on $i$.
	Otherwise we choose $f'', \gb''$ such that
	$\forall i < i_0$ $j_\Es(f'')(\gb'') \incompatible 
		j_\Es(f^{\gns_0, i}_0)(\ga^{\gns_0, i})$
	and
	$j_\Es(f'')(\gb'') \leq
		j_\Es(F^{p_0})(\gk(\gp_{\ga, \gb}(\gns_0)), \gb)$.
 	We make sure to choose
	the $f''$ such that if $\gm < i_0$ is an inaccessible then
	$j_\Es(f'')(\gb'') \not\in j_\Es(R)(\gk^0(\gns_0), \gm)$.
	We set
    \begin{align*}
	    & q'' = u^{\gns_0, i}_0,
	    \\
	    & \ga'' \Ege \ga^{\gns_0, i}, \gb''.
    \end{align*}

\ifnum\article=0
\enlargethispage*{10pt}
\fi
\item $i_0$ is limit: If
	$\ordof{j_\Es(f^{\gns_0, i}_0)(\ga^{\gns_0, i})} {i < i_0}$
	is a maximal anti-chain below
		$j_\Es(F^{p_0})(\gk(\gp_{\ga, \gb}(\gns_0)), \gb)$
	we terminate the induction on $i$.
	Otherwise we choose $f'', \gb''$ such that
	$j_\Es(f'')(\gb'') \incompatible 
		j_\Es(f^{\gns_0, i}_0)(\ga^{\gns_0, i})$ 
	for all $i < i_0$
	and
	 $j_\Es(f'')(\gb'') \leq
		j_\Es(F^{p_0})(\gk(\gp_{\ga, \gb}(\gns_0)), \gb)$.
	We make sure to choose
	the $f''$ such that if $\gm < i_0$ is an inaccessible then
	$j_\Es(f'')(\gb'') \not\in j_\Es(R)(\gk^0(\gns_0), \gm)$.
	Choose $\ga'' \in N$ such that
	    $\forall i < i_0$ $\ga'' \Egt \ga^{\gns_0, i}, \gb''$.
	We set
    \begin{align*}
	    q'' = \bigunion_{i < i_0} u^{\gns_0, i}_0.
    \end{align*}
\end{itemize}
Set
    \begin{align*}
        & u''_1 = \big(q''_{\ordered{\gp_{\ga, \ga''}(\gns_0)}}\big)_1,
	\\
	& u''_0 = \big(q''_{\ordered{\gp_{\ga, \ga''}(\gns_0)}}\big)_0,
	\\
	&T''_0 =   \gp^{-1}_{\ga'',\gb} \,
			T^{p_0} \setminus \gp_{\ga, \ga''}(\gns_0),
	\\
	&F''_{0} =
	   F^{p_0}  \circ \gp_{\ga'', \gb},
	\\
	&T''_1 = \emptyset,
	\\
	& f''_1 = f^{p_0} \circ \gp_{\ga, \gb} (\gk(\gns_0)),
	\\
	&F''_1 = \emptyset.
\end{align*}
Using the corollary of \ref{DenseAntiHomogen-0l} construct 
$q'''_0 \leq^* u''_0 \union 
                \set{\ordered{
			\Es_{\ga''}, T''_0, 
			f'' \circ \gp_{\ga'', \gb''},
			F''_0
		}}$ 
	and 
		$B^{\gns_0, i_0}$ a maximal anti-chain below
	$p_{l..1} \append u''_1 \union \set{f''_1}$.
So for each $b \in B^{\gns_0, i_0}$ there is
$S \subseteq [T^{q'''_0}]^n$, an $\mc(q'''_0)$-tree, such that
	\begin{multline*}
	\forall \ordered{\gns_1,\dotsc,\gns_{n}} \in S\ 
	\exists t'_n \append \dotsb \append t'_1 \leq^{**}
                (q'''_{0 \ordered{\gns_1,\dotsc,\gns_{n}}})_{n..1}
	\\
        	b \append t'_n \append \dotsb \append t'_1 \append
		(q'''_{0 \ordered{\gns_1,\dotsc,\gns_{n}}})_0 \in D.
	\end{multline*}
We set
\begin{align*}
& \ga^{\gns_0, i_0} = \gk(\mc (q'''_0)),
\\
& u^{\gns_0, i_0}_0 = q'' \union
	\setof{\ordered{\Es_\gga, q_0^{\prime\prime\prime \Es_\gga}}}
		{\Es_\gga \in \supp q'''_0 \setminus \supp q''},
\\
& T^{\gns_0, i_0}_0 = T^{q'''_0},
\\
& f^{\gns_0, i_0}_0 = f^{q'''_0},
\\
& F^{\gns_0, i_0}_0 = F^{q'''_0}.
\end{align*}
When the induction on $i$ terminates we have
\begin{align*}
\ordof {\ga^{\gns_0, i}, u^{\gns_0, i}_0, T^{\gns_0, i}_0,
                 f_0^{\gns_0, i}, F_0^{\gns_0, i}}
	       {i < \gk}.
\end{align*}
We complete step $\gns_0$ by setting
\begin{align*}
&  \ga^{\gns_0} \in N,
\\
& \forall i < \gk \  \ga^{\gns_0} \Egt \ga^{\gns_0, i},
\\
& u^{\gns_0}_0 = \bigunion_{i < \gk}
		u^{\gns_0, i}_0,
\\
& T^{\gns_0}_0 = \dsintersect_{i < \gk}
		\gp^{-1}_{\ga^{\gns_0}, \ga^{\gns_0, i}} \, T^{\gns_0, i}_0,
\\
& \forall i < \gk\ 
	F_0^{\gns_0} \leq F_0^{\gns_0, i}
	\circ
		\gp_{\ga^{\gns_0}, \ga^{\gns_0, i}}.
\end{align*}
When the induction on $\gns$ terminates 
we return to work in $V$ and
we have
\begin{align*}
\ordof{\ga^\gns, u^\gns_0, T^\gns_0, F^\gns_0}
      {\gns \in A}.
\end{align*}
We define the following function with domain $A$:
\begin{align*}
g(\gns)= \ordof{j_\Es(f^{\gns, i}_0)(\ga^{\gns, i})} 
			{i < \gk}.
\end{align*}
By the construction, $g(\gns)$ is a maximal anti-chain below 
$j_\Es(F^{p_0})(\gk(\gp_{\ga, \gb}(\gns)), \gb)$.
We note that $\ordof{j_0(f^{\gns, i}_0)(\ga^{\gns, i})} 
			{i < \gk} \in M_0$ as $M_0$ is closed under
$\gk$-sequences.
Hence $g(\gns) \in \ME$ as  $g(\gns)= i_{0, \Es}
	(\ordof{j_0(f^{\gns, i}_0)(\ga^{\gns, i})} 
			{i < \gk})$.
So  $j_\Es(g)(\ga) \in M^2_\Es$ is a maximal anti-chain below
$j^2_\Es(F^{p_0})(\gb, j_\Es(\gb))$.
As $I(\Es)$ is $j^2_\Es(R)(\gk, j_\Es(\gk))$-generic over $M^2_\Es$, there is $f_0$
such that
$j^2_\Es(f_0)(\ga, j_\Es(\ga)) \in I(\Es)$ and
$j^2_\Es(f_0)(\ga, j_\Es(\ga))$ is stronger than a condition in
$j_\Es(g)(\ga)$.
Note that we can use $\ga$ here because the generic was build through
the normal measure. If we would not have had this property we
would have enlarged $\ga$ to accommodate the intersection.
We combine everything into one condition, $p^{\prime*}_0$, as follows:
\begin{align*}
& p^{\prime*}_0 = \bigunion_{\gns \in A} u^{\gns}_0,
\\
& T^{p^{\prime*}_0} = \dsintersect_{\gns \in A} \gp^{-1}_{\ga, \ga^\gns} \, T^{\gns}_0,
\\
& f^{p^{\prime*}_0}(\gn_1) = f^{p_0} \circ \gp_{\ga, \gb} (\gn_1),
\\
& \forall \gns \in A \ F^{p^{\prime*}_0} \leq F_0^\gns \circ 
			\gp_{\ga, \ga^\gns}, f_0.
\end{align*}
We write what we have gained so far: For each $\gns \in A$
	there is $B_\gns$, a maximal anti-chain below
	$p_1 \append (p^{\prime*}_{0 \ordered{\gns}})_1$, such that
for each $b \in B_\gns$ there is
$S \subseteq [T^{p^{\prime*}_0}]^n$, an $\mc(p^{\prime*}_0)$-tree,
such that
\begin{multline*}
	\forall \ordered{\gns_1,\dotsc,\gns_{n}} \in S\ 
	\exists t'_n \append \dotsb \append t'_1 \leq^{**}
                (p^{\prime*}_{0 \ordered{\gns_1,\dotsc,\gns_{n}}})_{n..1}
	\\
        	b \append t'_n \append \dotsb \append t'_1 \append
		(p^{\prime*}_{0 \ordered{\gns_1,\dotsc,\gns_{n}}})_0 \in D.
\end{multline*}
We set $D_\gns = \setof {r \in P_\gns} {r\leq q,\ q \in B_\gns}$.
Then $D_\gns$ is a dense open subset of $P_\gns$. 
By invoking \ref{DenseHomogenTriv-l+1} for $D_\gns$, $p_1 \append (p^{\prime*}_{0\ordered{\gns}})_1$
we find $p_1(\gns) \append h(\gns) \leq^* p_1 \append (p^{\prime*}_{0\ordered{\gns}})_1$,
$S^1(\gns)$
such that
\begin{multline*}
	\forall \ordered{\gns_1,\dotsc,\gns_{n}} \in S^1(\gns)\ 
	\exists t'_n \append \dotsb \append t'_1 \leq^{**}
                (p_1(\gns)_{\ordered{\gns_1,\dotsc,\gns_{n}}})_{n..1}
	\\
        	t'_n \append \dotsb \append t'_1 \append
		(p_1(\gns)_{\ordered{\gns_1,\dotsc,\gns_{n}}})_0  \append
		h(\gns)
		\in D_\gns.
\end{multline*}
Immediately we see that there are $p^*_1$, $S^1$ such
that by removing a measure $0$ set from $A$ we get
$\forall \gns \in A$ $p^*_1 = p_1(\gns)$, $S^1 = S^1(\gns)$.
So after the shrinkage of $A$ we have for each $\gns \in A$
\begin{multline*}
	\forall \ordered{\gns_1,\dotsc,\gns_{n}} \in S^1\ 
	\exists t'_n \append \dotsb \append t'_1 \leq^{**}
                (p^*_{1\ordered{\gns_1,\dotsc,\gns_{n}}})_{n..1}
	\\
        	t'_n \append \dotsb \append t'_1 \append
		(p^*_{1\ordered{\gns_1,\dotsc,\gns_{n}}})_0  \append
		h(\gns)
		\in D_\gns.
\end{multline*}
We gather the additional information we have by setting 
$f^{p^*_0}(\gns) = h(\gns)$ and letting the condition $p^*_0$ be
$p^{\prime*}_0$ with $f^{p^*_0}$ substituted for 
$f^{p^{\prime*}_0}$.
So at this point we have the following
\begin{multline*}
	\forall \ordered{\gns_{1,1},\dotsc,\gns_{1,n_1}} \in S^1\ 
	\exists t'_{1,n_1} \append \dotsb \append t'_{1,1} \leq^{**}
                (p^*_{1\ordered{\gns_{1,1},\dotsc,\gns_{1,n_1}}})_{n_1..1}
	\\
	\forall \gns \in A\ 
	\exists S^0 \subseteq [T^{p^*_0}]^{n_0}\ 
	\\
	\forall \ordered{\gns_{0, 1},\dotsc,\gns_{0, n_0}} \in S^0\ 
	\exists t'_{0,n_0+1} \append \dotsb \append t'_{0,1} \leq^{**}
                (p^*_{0\ordered{\gns,\gns_{0,1},\dotsc,\gns_{0,n_0}}})_{n_0+1..1}
	\\
        	t'_{1,n_1} \append \dotsb \append t'_{1,1} \append
		(p^*_{1\ordered{\gns_{1,1},\dotsc,\gns_{1, n_1}}})_0  \append
	\\
        	t'_{0,n_0+1} \append \dotsb \append t'_{0,1} \append
		(p^*_{0\ordered{\gns, \gns_{0,1},\dotsc,\gns_{0, n_0}}})_0
		\in D.
\end{multline*}
Of course as $A \in E_\ga(0)$ the above is just a convoluted form of
\begin{multline*}
	\forall \ordered{\gns_{1,1},\dotsc,\gns_{1,n_1}} \in S^1\ 
	\exists t'_{1,n_1} \append \dotsb \append t'_{1,1} \leq^{**}
                (p^*_{1\ordered{\gns_{1,1},\dotsc,\gns_{1,n_1}}})_{n_1..1}
	\\
	\exists S^0 \subseteq [T^{p^*_0}]^{n_0}\ 
	\\
	\forall \ordered{\gns_{0, 1},\dotsc,\gns_{0, n_0}} \in S^0\ 
	\exists t'_{0,n_0} \append \dotsb \append t'_{0,1} \leq^{**}
                (p^*_{0\ordered{\gns_{0,1},\dotsc,\gns_{0,n_0}}})_{n_0..1}
	\\
        	t'_{1,n_1} \append \dotsb \append t'_{1,1} \append
		(p^*_{1\ordered{\gns_{1,1},\dotsc,\gns_{1, n_1}}})_0  \append
	\\
        	t'_{0,n_0} \append \dotsb \append t'_{0,1} \append
		(p^*_{0\ordered{\gns_{0,1},\dotsc,\gns_{0, n_0}}})_0
		\in D.
\end{multline*}
\end{proof}
\begin{proof}[proof of \ref{FullHomogen}]
The proof is done by induction on $l$.

The case $l=0$ is done in \ref{DenseHomogenTriv-0} for $\len(\Es)=0$ and
	\ref{DenseHomogen-0} for $\len(\Es) > 0$.

The case $l+1$ is done in \ref{DenseHomogenTriv-l+1} for $\len(\Es)=0$ and
	\ref{DenseHomogen-l+1} for $\len(\Es) > 0$.
\end{proof}
%
%
%
%
%
%
\ifnum\article=0
\newpage
\fi
\section{Prikry's condition} \label{Prikry'sCondition}
\ifnum\article=0
\enlargethispage*{10pt}
\fi
\begin{lemma} \label{WeakPrikry} \label{ToNoTree}
Let $\gs$ be a formula in the forcing language, $q \append p_k \append r 
	\in \PE$ and $S \subseteq [T^{p_k}]^m$ an $\mc(p_k)$-fat tree
such that
\begin{multline*}
\forall \ordered{\gns_1,\dotsc,\gns_{m}} \in S\ 
	\exists t'_m \append \dotsb \append t'_1 \leq^{**}
                (p_{k \ordered{\gns_1,\dotsc,\gns_{m}}})_{m..1}
	\\
                q \append t'_m \append \dotsb \append t'_1 \append
		(p_{k \ordered{\gns_1,\dotsc,\gns_{m}}})_0 \append r\decides \gs.
\end{multline*}
Then there is $p^*_k \leq^{**} p_k$ such that $q \append p^*_k \append r\decides \gs$.
\end{lemma}
\begin{proof}
Let
\begin{multline*}
 A^1_{\gns_1, \dotsc, \gns_{m-1}} = \setof 
		{\gns_m \in \Suc_{S}({\gns_1, \dotsc, \gns_{m-1}})}
		 {
\\
               q \append t'_m \append \dotsb \append t'_1 \append
		(p_{k \ordered{\gns_1,\dotsc,\gns_{m}}})_0 \append 
				r \forces \gs
		},
\end{multline*}
\begin{multline*}
 A^2_{\gns_1, \dotsc, \gns_{m-1}} = \setof 
		{\gns_m \in \Suc_{S}({\gns_1, \dotsc, \gns_{m-1}})}
		 {
\\
                q \append t'_m \append \dotsb \append t'_1 \append
		(p_{k \ordered{\gns_1,\dotsc,\gns_{m}}})_0 \append
					r \forces \lnot\gs
		}.
\end{multline*}
Then
\begin{align*}
& \Suc_{S}({\gns_1, \dotsc, \gns_{m-1}}) = A^1_{\gns_1, \dotsc, \gns_{m-1}} \union A^2_{\gns_1, \dotsc, \gns_{m-1}},
\\
& A^1_{\gns_1, \dotsc, \gns_{m-1}} \intersect A^2_{\gns_1, \dotsc, \gns_{m-1}} = 
			\emptyset.
\end{align*}
Let $A_{\gns_1, \dotsc, \gns_{m-1}}$ be the
$A^i_{\gns_1, \dotsc, \gns_{m-1}}$, $i \in \set{1,2}$ such that
\begin{align*}
& \exists \gx<\len(\Es)\ 
	A_{\gns_1, \dotsc, \gns_{m-1}} \in 
		{\mc(p_k)}(\gx).
\end{align*}
We choose $t'_m(\gns_1) \append \dotsb \append t'_1(\gns_1, \dotsc, \gns_m)$
such that for all $\ordered{\gns_1, \dotsc, \gns_{m-1}}$,
for all $\ordered{\gns_m} \in A_{\gns_1, \dotsc, \gns_{m-1}}$
\begin{align*}
	q \append t'_m(\gns_1) \append \dotsb \append t'_1(\gns_1, \dotsc, \gns_m) 
		\append
		(p_{k \ordered{\gns_1,\dotsc,\gns_{m}}})_0 \append r
			\forces \gs^i
\end{align*}
where $\gs^i \in \set{\gs, \lnot \gs}$ according to the selection of $i$.
By \ref{ToPreDense}, 
for each $\ordered{\gns_1, \dotsc, \gns_{m-1}} \in S$
there are $B_{\gns_1, \dotsc, \gns_{m-1}} \subseteq 
			A_{\gns_1, \dotsc, \gns_{m-1}}$, 
	$p^{m-1}_k(\gns_1, \dotsc, \gns_{m-1}) \leq^{**} 
		(p_{k \ordered{\gns_1, \dotsc, \gns_{m-1}}})_0$
such that 
below
 $q \append t'_{m}(\gns_1) \append \dotsb \append t'_2(\gns_1, \dotsc, \gns_{m-1}) 
		\append
		p^{m-1}_k(\gns_1,\dotsc, \linebreak[0] \gns_{m-1}) \append r$
\begin{multline*}
\setof{
	q \append t'_m(\gns_1) \append \dotsb \append 
		t'_1(\gns_1, \dotsc, \gns_m) 
		\append
		(p^{m-1}_k(\gns_1,\dotsc,\gns_{m-1})_{\ordered{\gns_m}})_0
		\append
		r
	}
	{
	\\
	 \gns_m \in B_{\gns_1, \dotsc, \gns_{m-1}} 
	}
\end{multline*}
is pre-dense.
Hence
\begin{align*}
q \append t'_{m}(\gns_1) \append \dotsb \append 
		t'_2(\gns_1, \dotsc, \gns_{m-1}) 
		\append
		p^{m-1}_k(\gns_1,\dotsc,\gns_{m-1}) \append r
			\decides \gs.
\end{align*}
Let $T^{m-1} = T^{p_k} \intersect \dsintersect_{\gns_1, \dotsc, \gns_{m-1}} 
				T^{p^{m-1}_k(\gns_1, \dotsc, \gns_{m-1})}$.
Let $p^{m-1}_k$ be the condition $p_k$ with its measure $1$ set substituted
by $T^{m-1}$. We get that for all 
	$\ordered{\gns_1, \dotsc, \gns_{m-1}} \in S$
\begin{align*}
q \append t'_{m}(\gns_1) \append \dotsb \append 
	t'_2(\gns_1, \dotsc, \gns_{m-1}) 
		\append
		(p^{m-1}_{k \ordered{\gns_1,\dotsc,\gns_{m-1}}})_0
		\append
		r
			\decides \gs.
\end{align*}
Letting $S^{m-1}$ be $S$ restricted to $m-1$ levels bring us
to the beginning of the proof but with $m-1$ instead of $m$.

Hence, repeating another $m-1$ steps as the above  build 
$p^0_k \leq^{**} \dotsb \leq^{**} p^{m-1}_k \leq^{**} p_k$ and 
$q \append p^0_k \append r \decides \gs$.
\end{proof}
\begin{theorem} \label{PrikryCondition}
Let $p \in \PE$, $\gs$ a formula in the forcing language. Then there
is $p^* \leq^* p$ such that $p^* \decides \gs$.
\end{theorem}
\begin{proof}
As usual we give the proof for the case $p = p_1 \append p_0$.

Let $D = \setof {r \leq p} {r \decides \gs}$. $D$ is a dense open subset
in $\PE$. By \ref{FullHomogen} there is $p' = p'_{1..0} \leq^* p$ such that
\begin{multline*}
\exists S^1\ \forall \ordered{\gns_{1, 1}, \dotsc,\gns_{1, n_1}} \in S^1\ 
\exists t'_{1, n_1} \append \dotsb \append t'_{1, 1} \leq^{**}
                (p'_{1 \ordered{\gns_{1, 1}, \dotsc, \gns_{1, n_1}}})_{n_1..1}
\\
\exists S^0\ \forall \ordered{\gns_{0, 1}, \dotsc,\gns_{0, n_0}} \in S^0\ 
\exists t'_{0, n_0} \append \dotsb \append t'_{0, 1} \leq^{**}
                (p'_{0 \ordered{\gns_{0, 1}, \dotsc, \gns_{0, n_0}}})_{n_0..1}
\\
t'_{1, n_1} \append \dotsb \append t'_{1, 1} \append
                (p'_{1 \ordered{\gns_{1, 1}, \dotsc, \gns_{1, n_1}}})_{0}
		\append
\\
t'_{0, n_0} \append \dotsb \append t'_{0, 1} \append
                (p'_{0 \ordered{\gns_{0, 1}, \dotsc, \gns_{0, n_0}}})_{0}
		\decides \gs.
\end{multline*}
We use the above formula to fix $S^1$. Then for each 
$\ordered{\gns_{1,1}, \dotsc, \gns_{1, n_1}} \in S^1$ we fix $t'_1(\gns_{1, 1})$,
$\dotsc$, $t'_{1, n_1}(\gns_{1, 1}, \dotsc, \gns_{1, n_1})$, 
$S^0(\gns_{1, 1}, \dotsc, \gns_{1, n_1})$.
In the same way we fix $t'_{0,1}$, $\dotsc$, $t'_{0, n_0}$ for each 
$\ordered{\gns_{0,1}, \dotsc, \gns_{0, n_0}} \in 
	S^0(\gns_{1, 1}, \dotsc, \gns_{1, n_1})$.

By \ref{WeakPrikry} for each 
	$\ordered{\gns_{1,1}, \dotsc, \gns_{1, n_1}} \in S^1$
there is $p'_0(\gns_{1,1}, \dotsc, \gns_{1, n_1}) \leq^{**} p'_0$ such that
\begin{align*}
t'_{1, n_1} \append \dotsb \append t'_{1, 1} \append
                (p'_{1 \ordered{\gns_{1, 1}, \dotsc, \gns_{1, n_1}}})_{0}
		\append
                p'_0(\gns_{1,1}, \dotsc, \gns_{1, n_1})
		\decides \gs.
\end{align*}
We choose $p^*_0 \leq^{**} p'_0(\gns_{1,1}, \dotsc, \gns_{1, n_1})$
for all $\ordered{\gns_{1,1}, \dotsc, \gns_{1, n_1}} \in S^1$.
Hence we get
\begin{align*}
\forall \ordered{\gns_{1,1}, \dotsc, \gns_{1, n_1}} \in S^1\ 
t'_{1, n_1} \append \dotsb \append t'_{1, 1} \append
                (p'_{1 \ordered{\gns_{1, 1}, \dotsc, \gns_{1, n_1}}})_{0}
		\append
                p^*_0
		\decides \gs.
\end{align*}
Invoking \ref{WeakPrikry} again we get $p^*_1 \leq^{**} p'_1$ such that
$p^*_1 \append p^*_0 \decides \gs$.
\end{proof}
With Prikry condition at our hand and $\ordered{\PE/\Pe, \leq^*}$ being
$\gk^0(\ges)^+$-closed
we get
\begin{theorem} \label{NoNewSubsets}
Let	$G$ be $\PE/\Pe$-generic.
Then $\Pset^{V[G]}(\gk^0(\ges)) = \Pset(\gk^0(\ges))$.
\end{theorem}
%
%
%
%
%
%
%
%
%
\ifnum\article=0
\newpage
\fi
\section{Properness} \label{Properness}
\ifnum\article=0
\enlargethispage*{10pt}
\fi
The following definitions (which are \emph{not} used in this work)
are due to Saharon Shelah \cite{ProperForcing}.
\begin{definition*}
Let $\gc$ be large enough, $N \subelem H_\gc$,
	$\power{N} = \ha_0$, 	$\PE \in N$.
$p \in \PE$ is called $\ordered{N, \PE}$-generic if
\begin{align*}
p \forces_{\PE} \formula{\CN{G} \text{ is } \VN{P}_\Es \text{-generic over } \VN{N}}.
\end{align*}
\end{definition*}
\begin{definition*}
The forcing $\PE$ is called proper if given
$N \subelem H_\gc$,
	$\power{N} = \ha_0$, 
	$\PE \in N$, $p \in \PE \intersect N$
there is $q \leq p$ which is $\ordered{N, \PE}$-generic.
\end{definition*}
We adapt the above definitions to handle elementary submodels of size $\gk$.
We keep the names from the original definitions.
\begin{definition}
Let $\gc$ be large enough, $N \subelem H_\gc$,
	$\power{N} = \gk$, $N \supset \gk$, $N \supset N^{\upto\gk}$,
	$\PE \in N$.
$p \in \PE$ is called $\ordered{N, \PE}$-generic if
\begin{align*}
p \forces_{\PE} \formula{\CN{G} \text{ is } \VN{P}_\Es \text{-generic over } \VN{N}}.
\end{align*}
\end{definition}
\begin{definition}
The forcing $\PE$ is called proper if given
$N \subelem H_\gc$,
	$\power{N} = \gk$, $N \supset \gk$, $N \supset N^{\upto\gk}$,
	$\PE \in N$, $p \in \PE \intersect N$
there is $q \leq p$ which is $\ordered{N, \PE}$-generic.
\end{definition}
\begin{theorem} \label{Proper}
$\PE$ is proper.
\end{theorem}
\begin{proof}
Let $\gc$ be large enough, $N \subelem H_\gc$,
	$\power{N} = \gk$, $N \supset \gk$, $N \supseteq N^{\upto\gk}$,
	$\PE \in N$,
	$p = p_{l..0} \in \PE \intersect N$.
We  find $p^*_0 \leq^* p_0$ such that $p_{l..1} \append p^*_0$ is
$\ordered{N, \PE}$-generic.

Let $\setof{D_\gx} {\gx < \gk}$ be enumeration of all dense open subsets
of $\PE$ appearing in $N$. Note that for $\gx_0 < \gk$ we have
	$\setof{D_\gx} {\gx < \gx_0} \in N$.

Let $\Es_\gb = \mc(p_0)$.
Choose $\ga \in \dom \Es$ such that $\ga \Egt \gga$ for all 
$\gga \in \dom \Es \intersect N$.
Let $A' = T^{p_0}$.
We shrink $A'$ a bit so that the following is satisfied:
$\forall \gns \in A'\ 
	\power{\setof{\gms \in A} {\gk^0(\gms) < \gk^0(\gns)}} 
		\leq \gk^0(\gns)$.
Let $A = \setof {\ordered{\gns_1, \dotsc, \gns_n} \in [A']^{\upto \gw} }
			{\len(\gns_n) = 0,\ 
		\gk^0(\gns_1) < \dotsb < \gk^0(\gns_n)}$.
Elements of $A$ are written in the form $\gnv$.
That is $\gnv = \ordered{\gns_1, \dotsc, \gns_n}$.
By $\max^0 \gnv$ we mean $\gns^0_n$.
Let $\preceq$ be well ordering of $A$ such that
$\forall \gnv,\gmv \in A \ 
	\gnv \preceq \gmv \implies \max^0 \gnv \leq \max^0 \gmv$.
We start an induction on $\gnv$ in which we build
\begin{align*}
\ordof{\ga^\gnv, u^\gnv_0, T^\gnv_0, F^\gnv_0}
      {\gnv \in A}.
\end{align*}
Assume that we have constructed
$\ordof{\ga^\gnv, u^\gnv_0, F^\gnv_0, T^\gnv_0}
      {\gnv \prec \gnv_*}$.
Recall our convention: $\gnv_* = \ordered{\gns_{*1}, \dotsc, \gns_{*k}}$.
We start working in $N$.

Set the following:
\begin{itemize}
\item 
        $\gnv_*$ is $\prec$-minimal:
    \begin{align*}
	& q' = p_0 \setminus \set{\ordered{\Es_\gb, T^{p_0}, f^{p_0},
					F^{p_0}}},
	\\
	& \ga' = \gb.
    \end{align*}
\item  
        $\gnv_*$ is the immediate $\prec$-successor of $\gnv$:
    \begin{align*}
	    & q' = u^\gnv_0,
	    \\
	    & \ga' = \ga^\gnv.
    \end{align*}
\item 
        $\gnv_*$ is $\prec$-limit: Choose $\ga' \in N$ such that
    $\forall \gnv \prec \gnv_*$ $\ga' \Egt \ga^\gnv$ 
    and set
    \begin{align*}
	    q' = \bigunion_{\gnv \prec \gnv_*} u^\gnv_0.
    \end{align*}
\end{itemize}
We make an induction on $i$ which builds
$\ordof {\ga^{\gnv_*, i}, u^{\gnv_*, i}_0, T^{\gnv_*, i}_0,
                 f^{\gnv_*, i}_0,
                 F^{\gnv_*, i}_0}
	       {i < \gk}$.
Assume we have constructed
\begin{align*}
\ordof {\ga^{\gnv_*, i}, u^{\gnv_*, i}_0, T^{\gnv_*, i}_0,
                 f^{\gnv_*, i}_0, F^{\gnv_*, i}_0}
	       {i < i_*},
\end{align*}
and we do step $i_*$.
\ifnum\article=0
\enlargethispage*{10pt}
\fi
\begin{itemize}
\item 
        $i_* = 0$:
    \begin{align*}
	& q'' = q',
	\\
	& \ga'' = \ga',
	\\
	& f'' = F^{p_0}(\gk(\gp_{\ga, \gb}(\gns_{*k})), \gp_{\ga'',\gb}(-)).
    \end{align*}
\item 
        $i_* = i + 1$: If 
	$\ordof{j_\Es(f^{\gnv_*, i}_0)(\ga^{\gnv_*, i})} {i < i_*}$
	is a maximal anti-chain below
		$j_\Es(F^{p_0})(\gk(\gp_{\ga, \gb}(\gns_{*k})), \gb)$
	we terminate the induction on $i$.
	Otherwise we pick $f'', \gb''$ such that
	$\forall i < i_*$ $j_\Es(f'')(\gb'') \incompatible 
		j_\Es(f^{\gnv_{*}, i}_0)(\ga^{\gnv_*, i})$
	and
	$j_\Es(f'')(\gb'') \leq
		j_\Es(F^{p_0})(\gk(\gp_{\ga, \gb}(\gns_{*k})), \gb)$.
 	We make sure to choose
	the $f''$ such that if $\gm < i_*$ is an inaccessible then
	$j_\Es(f'')(\gb'') \not\in j_\Es(R)(\gk^0(\gns_{*k}), \gm)$.
	We set
    \begin{align*}
	    & q'' = u^{\gnv_*, i}_0,
	    \\
	    & \ga'' \Ege \ga^{\gnv_*, i}, \gb''.
    \end{align*}
\item 
        $i_*$ is limit: If
	$\ordof{j_\Es(f^{\gnv_*, i}_0)(\ga^{\gnv_*, i})} {i < i_*}$
	is a maximal anti-chain below
		$j_\Es(F^{p_0})(\gk(\gp_{\ga, \gb}(\gns_{*k})), \gb)$
	we terminate the induction on $i$.
	Otherwise we pick $f'', \gb''$ such that
	$\forall i < i_*$ $j_\Es(f'')(\gb'') \incompatible 
		j_\Es(f^{\gnv_*, i}_0)(\ga^{\gnv_*, i})$ 
	and
	 $j_\Es(f'')(\gb'') \leq
		j_\Es(F^{p_0})(\gk(\gp_{\ga, \gb}(\gns_{*k})), \gb)$.
	We make sure to choose
	the $f''$ such that if $\gm < i_*$ is an inaccessible then
	$j_\Es(f'')(\gb'') \not\in j_\Es(R)(\gk^0(\gns_{*k}), \gm)$.
	Choose $\ga'' \in N$ such that
	    $\forall i < i_*$ $\ga'' \Egt \ga^{\gnv_*, i}, \gb''$ .
	We set
    \begin{align*}
	    q'' = \bigunion_{i < i_*} u^{\gnv_*, i}_0.
    \end{align*}
\end{itemize}
Set
    \begin{align*}
	& u''_{k+1..0} = (q'' \union \set{ \Es_{\ga''}, 
				\gp^{-1}_{\ga'', \gb} T^{p_0},
				f'' \circ \gp_{\ga'', \gb''},
			        F^{p_0}  \circ \gp_{\ga'', \gb}
			})_
			{\ordered{\gp_{\ga, \ga''}(\gnv_*)}}.
\end{align*}
We construct 	$q'''_0 \leq^* u''_0$
by invoking \ref{DenseAntiHomogen-0l} repeatedly for each $\setof {D_\gx} {\gx < \max^0 \gnv}$
starting from $p_{l..1} \append u''$.
We write explicitly what we have here: For each $\gx < \max^0 \gnv$
there is $B^{\gnv_*, i_*, \gx}$, a maximal anti-chain below
$p_{l..1} \append u''_{k+1..1}$ such that for each $b \in B^{\gnv_*, i_*, \gx}$
there are $r_0 \geq^* q'''_0$, $S \subseteq [T^{r_0}]^n$, an $\mc(r_0)$-fat
tree such that
\begin{multline*}
\forall \ordered{\gms_1,\dotsc,\gms_{n}} \in S\ 
\exists t'_n \append \dotsb \append t'_1 \leq^{**}
                (r_{0 \ordered{\gms_1,\dotsc,\gms_{n}}})_{n..1}
\\
        b \append t'_n \append \dotsb \append t'_1 \append
		(r_{0 \ordered{\gms_1,\dotsc,\gms_{n}}})_0 \in D,
\end{multline*}
and
\begin{align*}
\setof{
        b \append t'_n \append \dotsb \append t'_1 \append
		(r_{0 \ordered{\gms_1,\dotsc,\gms_{n}}})_0
	} 
	{
	\ordered{\gms_1,\dotsc,\gms_{n}} \in S
	}
\end{align*}
is pre-dense below $b \append r_0$.
We set
\begin{align*}
& \ga^{\gnv_*, i_*} = \gk(\mc (q'''_0)),
\\
& u^{\gnv_*, i_*}_0 = q'' \union
	\setof{\ordered{\Es_\gga, q_0^{\prime\prime\prime \Es_\gga}}}
		{\Es_\gga \in \supp q'''_0 \setminus \supp q''},
\\
& T^{\gnv_*, i_*}_0 = T^{q'''_0},
\\
& f^{\gnv_*, i_*}_0 = f^{q'''_0},
\\
& F^{\gnv_*, i_*}_0 = F^{q'''_0}.
\end{align*}
When the induction on $i$ terminates we have
\begin{align*}
\ordof {\ga^{\gnv_*, i}, u^{\gnv_*, i}_0, T^{\gnv_*, i}_0,
                 f_0^{\gnv_*, i}, F_0^{\gnv_*, i}}
	       {i < \gk}.
\end{align*}
We complete step $\gnv_*$ by setting
\begin{align*}
&  \ga^{\gnv_*} \in N,
\\
& \forall i < \gk \  \ga^{\gnv_*} \Egt \ga^{\gnv_*, i},
\\
& u^{\gnv_*}_0 = \bigunion_{i < \gk}
		u^{\gnv_*, i}_0,
\\
& T^{\gnv_*}_0 = \dsintersect_{i < \gk}
		\gp^{-1}_{\ga^{\gnv_*}, \ga^{\gnv_*, i}} \, T^{\gnv_*, i}_0,
\\
& \forall i < \gk\ 
	F_0^{\gnv_*} \leq F_0^{\gnv_*, i}
	\circ
		\gp_{\ga^{\gnv_*}, \ga^{\gnv_*, i}}.
\end{align*}
When the induction on $\gnv$ terminates 
we return to work in $V$ and
we have
\begin{align*}
\ordof{\ga^\gnv, u^\gnv_0, T^\gnv_0, F^\gnv_0}
      {\gnv \in A}.
\end{align*}
For each $\ordered{\gns_1, \dotsc, \gns_k} \in [T^{p^*_0}]^{\upto \gw}$
we define the following function with domain 
	$\setof {\gms \in A'} {\len(\gms) = 0}$
\begin{align*}
g^{\gns_1, \dotsc, \gns_k}(\gms)= \ordof{j_\Es(f^{\ordered{\gns_1, \dotsc, \gns_k, \gms}, i}_0)
			(\ga^{\ordered{\gns_1, \dotsc, \gns_k, \gms}, i})} 
			{i < \gk}.
\end{align*}
So  $g^{\gns_1, \dotsc, \gns_k}(\gms)$ is a maximal 
anti-chain below 
$j_\Es(F^{p_0})(\gk(\gp_{\ga, \gb}(\gms)), \gb)$.

We note that 
	$\ordof{j_0(f^{\ordered{\gns_1, \dotsc, \gns_k,\gms}, i}_0)
		(\ga^{\ordered{\gns_1, \dotsc, \gns_k, \gms}, i})} 
			{i < \gk} \in M_0$
as $M_0$ is closed under $\gk$-sequence.
As
 $g^{\gns_1, \dotsc, \gns_k}(\gms)= i_{0, \Es}
	(\ordof{j_0(f^{\ordered{\gns_1, \dotsc, \gns_k,\gms}, i}_0)
		(\ga^{\ordered{\gns_1, \dotsc, \gns_k, \gms}, i})} 
			{i < \gk})$
we get
 $g^{\gns_1, \dotsc, \gns_k}(\gms) \in \ME$.
So  $j_\Es(g^{\gns_1, \dotsc, \gns_k})(\ga) \in M^2_\Es$ is a maximal 
anti-chain below $j^2_\Es(F^{p_0})(\gb, j_\Es(\gb))$.
As $I(\Es)$ is $j^2_\Es(R)(\gk, j_\Es(\gk))$-generic over $M^2_\Es$, there is 
$f^{\gns_1, \dotsc, \gns_k}$
such that
$j^2_\Es(f^{\gns_1, \dotsc, \gns_k})(\ga, j_\Es(\ga)) \in I(\Es)$ and
$j^2_\Es(f^{\gns_1, \dotsc, \gns_k})(\ga, j_\Es(\ga))$ is stronger than 
a condition in
$j_\Es(g^{\gns_1, \dotsc, \gns_k})(\ga)$.
Note that we can use $\ga$ here because the generic was built through
the normal measure. If we would not have had this property we
would have enlarged $\ga$ to accommodate the intersection.
We combine everything into one condition, $p^{*}_0$, as follows:
\begin{align*}
& p^*_0 = \bigunion_{\gnv \in A} u^{\gnv}_0,
\\
& T^{p^*_0} = \dsintersect_{\gnv \in A} \gp^{-1}_{\ga, \ga^\gnv} \, T^{\gnv}_0,
\\
& f^{p^*_0}(\gn_1) = f^{p_0} \circ \gp_{\ga, \gb} (\gn_1),
\\
& \forall \ordered{\gns_1, \dotsc, \gns_k, \gms} \in A \ 
		F^{p^*_0} \leq F_0^{\ordered{\gns_1, \dotsc, \gns_k, \gms}} \circ 
			\gp_{\ga, \ga^{\ordered{\gns_1, \dotsc, \gns_k, \gms}}},
		 f^{\gns_1, \dotsc, \gns_k}.
\end{align*}
We claim that $p_{l..1} \append p^*_0$ is $\ordered{N, \PE}$-generic.

So, let $G$ be $\PE$-generic with $p_{l..1} \append p^*_0 \in G$.
Let $D \in N$ be a dense open subset of $\PE$.
There is $\gx < \gk$ such that $D = D_\gx$.
Let $\gnv = \ordered{\gns_1, \dotsc, \gns_k, \gms} \in A$ be such that
$p_{l..1} \append p^*_{0 \ordered{\gnv}} \in G$ and $\gk^0(\gms) > \gx$.
For convenience let us set $u_{k+2..0} = p^*_{0 \ordered{\gnv}}$.
By the construction there is $B \in N$, a maximal anti-chain below
$p_{l..1} \append u_{k+2..1}$, such that for each $b \in B$ there are
$r_0 \geq^* u_0$, $r_0 \in N$ and $S \in N$, an $\mc(r_0)$-fat tree, such that
\begin{multline*}
\forall \ordered{\gms_1,\dotsc,\gms_{n}} \in S\ 
\exists t'_n \append \dotsb \append t'_1 \leq^{**}
                (r_{0 \ordered{\gms_1,\dotsc,\gms_{n}}})_{n..1}
\\
        b \append t'_n \append \dotsb \append t'_1 \append
		(r_{0 \ordered{\gms_1,\dotsc,\gms_{n}}})_0 \in D_\gx
\end{multline*}
and
\begin{align*}
\setof{
        b \append t'_n \append \dotsb \append t'_1 \append
		(r_{0 \ordered{\gms_1,\dotsc,\gms_{n}}})_0
	} 
	{
	\ordered{\gms_1,\dotsc,\gms_{n}} \in S
	}
\end{align*}
is pre-dense below $b \append r_0$. Moreover, this pre-dense set is
contained in $N$.

We do the natural factoring $G/p_{l..1} \append u_{k+2..0} = G_\gms/p_{l..1}\append u_{k+2..1} \times
	G_\Es/u_0$. 
By genericity there is $b \in B \intersect G_\gms/p_{l..1} \append u_{k+2..1}$.
Necessarily $b \append u_0 \in G$. Hence $b \append r_0 \in G$.
So there is 	$\ordered{\gms_1,\dotsc,\gms_{n}} \in S$
such that        $b \append t'_n \append \dotsb \append t'_1 \append
		(r_{0 \ordered{\gms_1,\dotsc,\gms_{n}}})_0 \in G$.

So we finally got 
$b \append t'_n \append \dotsb \append t'_1 \append
		(r_{0 \ordered{\gms_1,\dotsc,\gms_{n}}})_0 \in 
		D \intersect G \intersect N$.
\end{proof}
We remind the reader of our convention that when $\gt_1 < \gt_2$ we have
$P_{\Es\restricted \gt_2} \subseteq P_{\Es \restricted \gt_1}$.
For the sake of completeness we mention the following rather obvious
propositions.
\begin{proposition} \label{RemainPreDense}
Assume that we have $p = p_l \append \dotsb \append p_0 \in \PE$, 
$S \subseteq [T^{p_0}]^n$ an $\mc(p_0)$-fat tree and
$t'_n(\gns_1) \append \dotsb \append t'_1(\gns_1, \dotsc, \gns_n) 
	\leq^{**} (p_{0 \ordered{\gns_1, \dotsc, \gns_n}})_{n..1}$
 such that
$D = \setof {p_{l..1} \append t'_n(\gns_1) \append \dotsb \append 
	t'_1(\gns_1, \dotsc, \gns_n) \append
	(p_{0 \ordered{\gns_1, \dotsc, \gns_n}})_0} 
	{\ordered{\gns_1, \dotsc, \gns_n} \in S}$
is pre-dense  below $p$. Let $\gt \leq \len(\Es)$ be such that $S$ is
$\mc(p_0) \restricted \gt$-fat tree, 
then $D$ is pre-dense in $P_{\Es \restricted \gt}$ below $p$.
\end{proposition}
\begin{proposition}
Let $\gc$ be large enough, $N \subelem H_\gc$,
	$\power{N} = \gk$, $N \supset \gk$, $N \supset N^{\upto\gk}$,
	$\PE \in N$,
	$p = p_{l..0} \in \PE \intersect N$.
Assume $S \in N$, $S \subseteq [T^{p_0}]^n$ is an $\mc(p_0)$-fat tree.
Then $S$ is $mc(p_0) \restricted \gt$-fat tree for each
$\sup \setof {\gt'+1} {\gt' \in \len(\Es) \intersect N} \leq \gt \leq \len(\Es)$.
\end{proposition}
With these propositions in mind we see that the properness proof
actually gave us more:
\ifnum\article=0
\enlargethispage*{10pt}
\fi
\begin{theorem} \label{PEgeneric}
Let $\gc$ be large enough, $N \subelem H_\gc$,
	$\power{N} = \gk$, $N \supset \gk$, $N \supset N^{\upto\gk}$,
	$\PE \in N$,
	$p = p_{l..0} \in \PE \intersect N$.
Then there is $p^*_0 \leq^* p_0$ such that
for all
$\sup \setof {\gt'+1} {\gt' \in \len(\Es) \intersect N} \leq 
					\gt \leq \len(\Es)$
\begin{align*}
p_{l..1} \append p^*_0 \forces_{P_{\Es \restricted \gt}}
	\formula{
		\CN{G} \text{ is } \VN{P}_\Es\text{-generic over } \VN{N}
	}.
\end{align*}
\end{theorem}
The following is a method to get a generic over elementary submodel from 
a `small' forcing.
It is given here, even though it does not belong to this section, 
as it has the same proof as \ref{Proper}.
\begin{definition}
Let $s \subseteq \Es$. We define
\begin{align*}
P_{\Es} \restricted s = \setof 
	{ p \in \PE}
	       { \supp p_0 \union \mc(p_0) \subseteq s}
\end{align*}
with $\leq$, $\leq^*$ inherited from $\PE$.
\end{definition}
\begin{definition}
Let $s \subseteq \Es$. If $G$ is $\PE$-generic then
\begin{align*}
G  \restricted s = G \intersect \PE \restricted s.
\end{align*}
\end{definition}
When $\power{s} \leq \gk$ this forcing is a somewhat convoluted Radin forcing.
Simple analysis reveals that
$\forces_{P_\Es \restricted s} \formula{2^{\gk} = \gk^{+}}$.
We also point out that $P_{\Es} \restricted s$ is completely embedded in 
$\PE$.
Hence, if $G$ is $\PE$-generic then  $G \restricted s$
is $P_\Es \restricted s$-generic.
\begin{theorem} \label{GenericAlpha}
Let $\gc$ be large enough, $N \subelem H_\gc$, $\power{N} = \gk$,
$N \supset \gk$, $N \supset N^{\upto \gk}$,
$\PE \in N$, $p = p_{l..0} \in N$.
Then there is $p^*_0 \leq^* p_0$ such that
\begin{align*}
p_{l..1} \append p^*_0 \forces_{\PE \restricted s} \formula {
		\CN{G} \text{ is } \VN{P}_\Es \text{-generic over } \VN{N}
	}
\end{align*}
where $s = \supp p^*_0 \union \mc(p^*_0)$.

Moreover, if $G$ is $\PE$-generic with $p_{l..1} \append p^*_0 \in G$
then $H = G \intersect 
P_{\Es} \restricted s$ is $\PE \restricted s$-generic and, obviously,  
$L[G] = L[H]$.
\end{theorem}
\ifnum\article=0
\newpage
\fi
\section{Cardinals in $V^{\PE}$} \label{CardinalStructure}
The following claim is just an exercise in properness.
\begin{claim} \label{NoCollapseSpecial}
$\forces_{\PE} \formula{ \VN{\gk^+} \text{ is cardinal} }$.
\end{claim}
\begin{proof}
If $\len(\Es) = 0$ then the claim is trivial. Hence we assume that
$\len(\Es) > 0$.
Let
$p \forces \formula{\GN{f} \func \VN{\gk} \to \VN{\gk^+}}$.
Choose $\gc$ large enough so that $H_\gc$ contains everything we are
interested in. 
By \ref{Proper} there are $p^* \leq^* p$, $N \subelem H_\gc$
such that
\begin{enumerate}
\item $p, \PE, \GN{f} \in N$,
\item $\power{N}= \gk$,
\item $N \supset \gk$,
\item $N \supset N^{\upto \gk}$,
\item $p^*$ is $\ordered{N, \PE}$-generic.
\end{enumerate}
Let us set
$\gl = N \intersect \gk^+$.
Note that $\gl < \gk^+$. 

Let $G$ be $\PE$-generic with $p^* \in G$.
The $\ordered{N,\PE}$-genericity ensures us that for all $\gx < \gk$ 
        $\GN{f}(\gx)^{N[G]} \in N$,
        $\GN{f}(\gx)^{N[G]}= \GN{f}(\gx)^{V[G]}$.
Hence $\ran \GN{f}^{V[G]} \subseteq \gl$. That is
$p^* \forces \formula{\GN{f} \text{ is bounded in } \gk^+}$.
\end{proof}
\begin{claim} \label{NoCollapseAbove}
No cardinals $> \gk$ are collapsed by $\PE$.
\end{claim}
\begin{proof}
$\gk^+$ is not collapsed by \ref{NoCollapseSpecial}.
No cardinals $\geq \gk^{++}$ are collapsed as 
$\PE$ satisfies $\gk^{++}$-c.c.
\end{proof}
\begin{claim} \label{LimitPower}
Assume $\len(\Es) > 0$.
$\forces_{\PE} \formula {2^{\gk}= \gk^{+3}}$.
\end{claim}
\begin{proof}
Let $G$ be $\PE$-generic.
For each $\ga \in \dom \Es$ define
$M^\ga = \bigunion \setof {p_0^{\Es_\ga}} {p \in G}$.
It is routine to check that for $\ga \not= \gb$ we have
$M^\ga \not= M^\gb$. Hence
	$\forces_{\PE} \formula {2^{\gk^0(\Es)} \geq \gk^0(\Es)^{+3}}$.

For the other direction let 
	$\forces_{\PE} \formula {\GN{A} \subseteq \VN{\gk}}$.
By \ref{GenericAlpha} there are
$\gc$ large enough, $N \subelem H_\gc$, $\power{N} = \gk$,
$\PE, \GN{A} \in N$ and $p \in G$ such that
$G \restricted s$ is $\PE$-generic over N
where $s = \supp p_0 \union
\mc(p_0)$.
Hence, $\GN{A}[G] \in V[G \restricted s]$. 
That is $\Pset^{V[G]}(\gk) = \bigunion_{s \in ([\Es]^\gk)_V} 
	\Pset^{V[G \restricted s]}(\gk)$. As $V[G] \satisfies \formula{
		\power{\Pset^{V[G \restricted s]}(\gk)}	= \gk^{+}
		}$
when $\power{s} \leq \gk$
and $\power{([\Es]^\gk)_V} = \gk^{+3}$
the proof is completed.
\end{proof}
%
%
Simple counting of anti-chains shows that we did not destroy
the behavior of the power function on $\gk^+$, $\gk^{++}$, $\gk^{+3}$:
\begin{claim}
Assume $\len(\Es) > 0$.
$\forces_{\PE} \formula {2^{(\gk^+)}= \gk^{+4},\ 
		2^{(\gk^{++})} = \gk^{+5},\ 
		2^{(\gk^{+3})} = \gk^{+6}
	}$.
\end{claim}
%
\begin{lemma}
\label{SmallForcing}
Let $p = p_l \append \dotsb \append p_k \append \dotsb \append p_0 \in G$ 
and $\ges$ be such that
$p_{l..k} \in \Pe$.
Let $G/p = G_\ges \times G_\Es$ be the obvious factoring.
Then
$\Pset(\gk^0(\ges))^{V[G]} = \Pset(\gk^0(\ges))^{V[G_\ges]}$.
\end{lemma}
\begin{proof}
This is immediate due to 
	$\Pset(\gk^0(\ges))^{V[G]} = \Pset(\gk^0(\ges))^{V[G/p]}$ and
\ref{NoNewSubsets}.
\end{proof}
\begin{claim} \label{NoCollapseSuccessor} \label{SuccessorStatus}
Let $G$ be $\PE$-generic with 
$p = p_l \append \dotsb \append p_k \append \dotsb \append p_0 \in G$ 
and $\ges$ be such that
$p_{l..k} \in \Pe$ and $\len(\ges) = 0$. Let $\gn = \gk(p_k^0)$.
Then, in V[G],
$\gn^+, \dotsc, \gn^{+6}$ remain cardinals,
all cardinals in $[\gn^{+7}, \gk^0(\ges)]$
are collapsed and 
	$2^{\gn^{+}} = \gn^{+4}$,
	$2^{\gn^{++}} = \gn^{+5}$,
	$2^{\gn^{+3}} = \gn^{+6}$,
	$2^{\gn^{+4}} = \gk^0(\ges)^+$,
	$2^{\gn^{+5}} = \gk^0(\ges)^{++}$,
	$2^{\gn^{+6}} = \gk^0(\ges)^{+3}$.
\end{claim}
\begin{proof}
Let
$G/p = G_{\ges} \times G_\Es$ be the natural factoring.
By \ref{SmallForcing} the fate of the cardinals in question is
decided by $P_{\ges}$. We note that
	$P_{\ges} = P_{\ges_2} \times R(\gn, \gk^0(\ges))$
where $p_{l..k+1} \in P_{\ges_2}$. So we factor $G_\ges = G_{\ges_2} 
		\times G_{\gn}$.
As it stands $R(\gn ,\gk^0(\ges))$ is $\gn^+$-closed. So in order
to prove the claim we  make some finer analysis.

We remind the reader that the $V$ we work with is a generic
extension of $V^*$ for a reverse Easton forcing.
Let $Q_1$ be the reverse Easton forcing up to $\gk^0(\ges_2)$ and
$H_1$ be its generic.
Let $Q_2$ be the forcing at stage $\gk^0(\ges_2)$
and $H_2$ be its generic over $V^*[H_1]$.
Let $Q_3$ be the rest of the reverse Easton forcing up to $\gk^0(\ges)$
and $H_3$ be its generic over $V^*[H_1][H_2]$.
Then we have
\begin{align*}
\Pset^{V[G]}(\gk^0(\ges)) = 
	 & \Pset^{V^*[H_1][H_2][H_3][G_\gn][G_{\ges_2}]}(\gk^0(\ges))
\end{align*}
Comparison of the forcings used to construct 
	$M^*_\Es[G_{\upto \gk}][G^\Es_\gk][G^\Es_{\downto \gk}][I_\Es]$ 
		(section \ref{PreparationForcing}) 
		and
	$V^*[H_1][H_2][H_3][G_\gn]$ 
shows that
the cardinal structure and power function
of the model $V^*[H_1][H_2][H_3][G_\gn]$ 
in the range $[\gn^+,\gk^0(\ges)^{+3}]$ behave in the same way as the 
cardinal structure and
power function of the model 
	$M^*_\Es[G_{\upto \gk}][G^\Es_\gk][G^\Es_{\downto \gk}][I_\Es]$=
        $M_\Es[I_\Es]$
in the range $[\gk^+, j_\Es(\gk)^{+3}]$.
From \ref{prototype} we see that in $V^*[H_1][H_2][H_3][G_\gn]$:
there are no cardinals in $[\gn^{+7}, \gk^0(\ges)]$,
$\gn^+, \dotsc ,\gn^{+6}$ are cardinals,
	$2^{\gn^+} = \gn^{+4}$,
	$2^{\gn^{++}} = \gn^{+5}$,
	$2^{\gn^{+3}} = \gn^{+6}$,
	$2^{\gn^{+4}} = \gk^0(\ges)^{+}$,
	$2^{\gn^{+5}} = \gk^0(\ges)^{++}$,
	$2^{\gn^{+6}} = \gk^0(\ges)^{+3}$.

Forcing with $P_{\ges_2}$ does not change the power function
and does not collapse cardinals
above $\gk^0(\ges_2)$ by the previous claims adapted to the current context.
\end{proof}
\begin{claim} \label{LimitStays}
Assume $\len(\Es) > 0$.
$\forces_{\PE} \formula{ \gk \text{ is a cardinal}}$.
\end{claim}
\begin{proof}
$\gk$ is limit ordinal
and by \ref{NoCollapseSuccessor}, there are unbounded number of cardinals
below $\gk$ which are preserved. Hence $\gk$ is preserved.
\end{proof}
With \ref{PEgeneric} at our disposal we can give a direct proof
of the following theorem.
It is the same one given
in \cite{GeneralRadinForcing} for proving the theorem
in  Radin forcing context.
\begin{theorem} \label{RemainsRegular}
If $\cf\len(\Es) > \gk$ then $\forces_{\PE} \formula{ \VN{\gk}
\text{ is regular} }$.
\end{theorem}
\begin{proof}
Let $\gl < \gk$, $p \in \PE$.
Let $\gc$ be large enough.
By \ref{PEgeneric} we have $N \subelem H_\gc$, 
$\power{N} = \gk$,
	$N \supset \gk$, $N \supset N^{\upto \gk}$,
	$p, \,\PE\in N$,
$p^* \leq^* p$ such that $p^* \forces_{P_{\Es \restricted \gt}}
	\formula{
	\CN{G} \text{ is } \VN{P}_\Es\text{-generic over } \VN{N}
	}$ 
for each
$\sup \setof{\gt'+1} {\gt' \in \len(\Es) \intersect N} \leq \gt 
				\leq \len(\Es)$.

Choose $\gt$ such that 
$\sup \setof{\gt'+1} {\gt' \in \len(\Es) \intersect N} \leq \gt 
				\leq \len(\Es)$,
$\gl < \cf\gt \leq \gk$.
This is possible because  $\cf \len(\Es) > \gk$, $\power{N} = \gk$.
Choose $q \leq_{P_{\Es \restricted \gt}} p^*$ such that
$q \forces_{P_{\Es \restricted \gt}} \formula{
	\VN{\gl} < \cf \VN{\gk} < \VN{\gk} }$.
Let $G$ be $P_{\Es \restricted \gt}$-generic with $q \in G$.
Of course, $p^* \in G$ also. Hence, $G$ is $\PE$-generic over $N$.

As $N[G] \in V[G]$ we have $\cf^{N[G]} \gk \geq \cf^{V[G]} \gk > \gl$.
So there is $r \in \PE \intersect G$, $r \leq p^*$ such that
$r \forces_\PE \formula { \cf \VN{\gk} > \VN{\gl} }$.
\end{proof}
\ifnum\article=0
\newpage
\fi
\section{Consistency theorem} \label{ConsistencyTheorem}
We state the consistency theorem we worked so much for.
\begin{theorem}
If there is $\Es$ such that $\power{\Es} = \gk^{+3}$,
$\cf \len(\Es) > \gk$ then it is consistent to have
the power function $2^\gm = \gm^{+3}$ for all cardinals $\gm$.
\end{theorem}
\begin{proof}
Let $p^* \in \PES$ such that $\gk(p^{*0})$ is inaccessible
and  $G$ be $\PE$-generic with $p^* \in G$.
(Forcing below an element of $\PES$ eliminates a finite number
of exceptions which we might otherwise have. That is if 
$p_1 \append p_0 \in G$ and $\gk^0(\mc(p_1)) < \gk(p^0_0)$ then
the interval $[\gk^0(\mc(p_1)), \gk(p^0_0)]$ is untouched by the forcing).
We set
\begin{align*}
& M = \bigunion \setof {p_0^{\Es_\gk}} {p \in G},
\\
& C = \bigunion \setof {\gk(p_0^{\Es_\gk})} {p \in G}.
\end{align*}
Note that $M$ is a Radin generic sequence for the extender sequence
$\Es_\gk$.
Hence $C \subset \gk$ is a club. The first ordinal in this 
club is $\gl = \gk(p^{*0})$.
We investigate the range $(\gl, \gk)$ in $V[G]$.
We note that, by \ref{SmallForcing}, for $\ges \in M$ it is enough to 
use $\Pe$ in order to understand $V_{\gk^0(\ges)}^{V[G]}$.
So let $\gm \in C$, $\gm > \gl$.
\begin{itemize}
\item
	$\gm \in \lim C$: Then there is $\ges \in M$ such that
		$\len(\ges) > 0$ and $\gk(\ges) = \gm$.
	By \ref{LimitStays}, $\gm$ remains a cardinal
	and by \ref{LimitPower}, $2^\gm = \gm^{+3}$.
\item
	$\gm \in C \setminus \lim C$: Then there is $\ges \in M$ such that
		$\len(\ges) = 0$ and $\gk(\ges) = \gm$.
	Let $\gm_2 \in C$ be the $C$-immediate predecessor of $\gm$.
	By \ref{SuccessorStatus} we have: $\gm_2^+, \dotsc, \gm_2^{+6}$
	are cardinals, there are no cardinals in $[\gm_2^{+7}, \gm]$,
	$2^{\gms_2^+} = \gm_2^{+4}$,
	$2^{\gms_2^{++}} = \gm_2^{+5}$,
	$2^{\gms_2^{+3}} = \gm_2^{+6}$,
	$2^{\gms_2^{+4}} = \gm^{+}$,
	$2^{\gms_2^{+5}} = \gm^{++}$,
	$2^{\gms_2^{+6}} = \gm^{+3}$.
\end{itemize}
In fact due to all the cardinals collapsed we have
$\setof {\gm \text{ is a cardinal }}{\gl < \gm < \gk}  = 
	\lim C \union \setof {\gm^+,\dotsc,\gm^{+6}} {\gm \in C}$.
Hence if $\gm \in (\gl, \gk)$ is a cardinal then
$2^{\gm} = \gm^{+3}$.
By \ref{RemainsRegular}, $\gk$ is an inaccessible cardinal.
Let $H$ be $\Col(\ha_0, \gl^+)_{V[G]}$-generic over $V[G]$.
In $V[G][H]$ $\gk$ remains inaccessible and
 $\forall \gm < \gk$ $2^\gm = \gm^{+3}$.
So $V_\gk^{V[G][H]}$ is a model of ZFC satisfying 
$\forall \gm$ $2^\gm = \gm^{+3}$.
\end{proof}
\ifnum\article=0
\newpage
\fi
\section{Concluding remarks} \label{ConcludingRemarks}
\subsection{Regarding The Power Function in Our Model}
Our forcing divides the cardinals into
$3$ categories. The first category contains the cardinals appearing in the 
club, $C$, generated by the normal Radin sequence. 
The second category contains the successors of cardinals in $C$ which are
below the length of the extender we  use. The third category contains
the cardinals above the length of the extender.
The gap on cardinals in each of these categories can be different.
We  give several examples to clarify this point. In all of them
we assume that $\cf \len (\Es) > \gk$ and $V_\gk$ of the generic extension
is the model of ZFC we are interested in.

Example:
	By just doing the extender based Radin forcing (that is,
	without the extra cardinal collapsing and Cohen forcings)
	starting from $j_\Es \func V \to \ME \supset V_{\gk+n}$
	we get that there is a generic extension with a club 
	$C \subset \gk$	and a power function
\ifnum\article=1
\enlargethispage*{10pt}
\fi
	\begin{align*}
	& 2^\gm = \begin{cases}
		\gl^{+n}	& \gl \in \lim C,\ \gm = \gl^{+k},\ 0 \leq k < n 
		\\
		\gm^{+}		& \text{otherwise}
		\end{cases}.
	\end{align*}

Example:
	By adding to the previous example the  collapse
	$\Col(\gl_1^{+n+1}, \gl_2)$ for each $\gl_1, \gl_2 \in C$ 
	successive points, we get the same power function.
	Our gain here is that the cardinals of the new model
	are `close' to $C$. Namely 
	the cardinals are $\bigunion \setof {\gl^+, \dotsc, \gl^{+n+1}} 
					{\gl \in C} \union \lim C$.

Example:
	The collapse we chose in the previous example is the lowest
	possible. We can use others if the need arises.
	Let $n = 3$.
	By doing a reverse Easton preparation on the 
	inaccessibles of $\C(\gl^+, \gl^{+5}) \times
	\C(\gl^{++}, \gl^{+7}) \times
	\C(\gl^{+3}, \gl^{+10})$
	and then invoking $\PE$ with the collapse
		$\Col(\gl_1^{+ 10}, \gl_2)$ we  get that
	the cardinals are $\bigunion \setof {\gl^+, \dotsc, \gl^{+10}} {\gl \in C}
	\union \lim C$ with 
	power function
	\begin{align*}
	& 2^\gm = \begin{cases}
		\gm^{+3}	& \gm \in \lim C
		\\
		\gm^{+4}	& \gm = \gl^+,\ \gl \in C
		\\
		\gm^{+5}	& \gm = \gl^{++},\ \gl \in C
		\\
		\gm^{+7}	& \gm = \gl^{+3},\ \gl \in C
		\\
		\gm^{+6}	& \gm = \gl^{+4},\ \gl \in C
		\\
		\gm^{+5}	& \gm = \gl^{+5},\ \gl \in C
		\\
		\gm^{+4}	& \gm = \gl^{+6},\ \gl \in C
		\\
		\gm^{+3}	& \gm = \gl^{+7},\ \gl \in C
		\\
		\gm^{+2}	& \gm = \gl^{+8},\ \gl \in C
		\\
		\gm^{+}		& \text{otherwise}
		\end{cases}.
	\end{align*}

Example:
	If we do the reverse Easton forcing as in the previous
	example and then invoke $\PE$ with the forcing
	$\Col(\gl_1^{ + 10}, \gl_2) \times \C(\gl_1^{+4}, \gl_2^{+}) \times
	\C(\gl_1^{+6}, \gl_2^{+5})$ we  get the same cardinals and
	the power function
	\begin{align*}
	& 2^\gm = \begin{cases}
		\gm^{+3}	& \gm \in \lim C
		\\
		\gm^{+4}	& \gm = \gl^+,\ \gl \in C
		\\
		\gm^{+5}	& \gm = \gl^{++},\ \gl \in C
		\\
		\gm^{+7}	& \gm = \gl^{+3},\ \gl \in C
		\\
		\gm^{+7}	& \gm = \gl^{+4},\ \gl \in C
		\\
		\gm^{+6}	& \gm = \gl^{+5},\ \gl \in C
		\\
		\gm^{+9}	& \gm = \gl^{+6},\ \gl \in C
		\\
		\gm^{+8}	& \gm = \gl^{+7},\ \gl \in C
		\\
		\gm^{+7}	& \gm = \gl^{+8},\ \gl \in C
		\\
		\gm^{+6}	& \gm = \gl^{+9},\ \gl \in C
		\\
		\gm^{+5}	& \gm = \gl^{+10},\ \gl \in C
		\\
		\gm^{+}		& \text{otherwise}
		\end{cases}.
	\end{align*}
	Note the following limitation. If our reverse Easton preparation
	would have contained $\C(\gl^+, \gl^{+4})$ then the case
	$\gm = \gl^+$ in the power function would have been like this
	\begin{align*}
	& 2^\gm = \begin{cases}
		\gm^{+3}	& \gm = \gl^+,\ \gl \in \lim C
		\\
		\gm^{+4}	& \gm = \gl^+,\ \gl \in C \setminus \lim C
		\end{cases}.
	\end{align*}

As can be seen we have quite a lot of freedom in setting the power
of the successors in these models. However, we do have a major limitation.
We get the same behavior over and over again. This is inherent to
our forcing. Another point is that this freedom is on a somewhat
insignificant set. It is a non-stationary set. And a very thin
non-stationary. It contains no limit cardinals.

A generalization of the second example above is as follows.
Assume that we have $f \func \gk \to \gk$ such that
$\cf \gl^{+f(\gl)} > \gl$ on a measure $1$ set
and
$j_\Es \func V \to \ME \supset V_{\gk + j(f)(\gk)}$.
If we force with $\PE$ adding the collapse $\Col(\gl_1^{+f(\gl_1)+1}, \gl_2)$
for each $\gl_1, \gl_2 \in C$ successive points then the cardinals
in the new model are $\bigunion \setof{\gm \text{ is cardinal}} 
	{\gl^+ \leq \gm \leq \gl^{+f(\gl)+1},\ \gl \in C} \union \lim C$
with power function
	\begin{align*}
	& 2^\gm = \begin{cases}
		\gl^{+f(\gl)}	& \gm \in \lim C
		\\
		\gl^{+f(\gl)}	& \gl < \gm < \gl^{+f(\gl)},\ 
					\gl \in \lim C,\ 
					\cf \gl^{+f(\gl)} > \gm
		\\
		\gl^{+f(\gl)+1}	& \gl < \gm < \gl^{+f(\gl)},\ 
					\gl \in \lim C,\ 
					\cf \gl^{+f(\gl)} \leq \gm
		\\
		\gm^{+}		&  \text{otherwise}
		\end{cases}.
	\end{align*}
We can, of course, do a preparation forcing and add Cohen forcings
along the normal Radin sequence as before. However, if $\gm \not \in C$
is a singular cardinal we have SCH on it.
A different method is needed in order to generate a gap on such cardinal.

We suggest the following attack and we stress that it is a
\emph{suggestion}. Unlike the previous examples which are immediate
consequences of our forcing notion, this attack require a deeper
research.
So, we assume that $\ME$ thinks that there is a cardinal $\gn$
between	$\gk$ and $j_\Es(\gk)$ carrying a $\gn+3$-strong extender, $F$.
Let $\Fs$ be the extender sequence of length 1 built from $F$.
Hence we can define $Q_\Fs$ in $\ME$, the forcing for adding
$\gn+3$ Prikry sequences to $\gn$.
Our idea is to add along $C$ reflections
of $Q_\Fs$. Hence, if $\gl_1, \gl_2 \in C$ are successive points
then we force with $Q_\ges$, the forcing
notion for adding $\gn'$+3 Prikry sequences to $\gn'$ for $\gn'$
lying between $\gl_1$ and $\gl_2$.
Recall that in order to have a Prikry like condition we need to have
a $Q_\Fs$-generic filter over $\ME$. Alas, we do not have one.
However, we do have a $\ordered{Q_\Fs, \leq^*}$-generic filter over $\ME$.
(Construct a generic filter over the normal ultrapower and then
 send it to $\ME$).
We think that with some modifications our proofs  go through using
this weaker generic filter.
\subsection{Regarding The Power Function}
Let $\text{Reg}$ be the class of regular cardinals.
We recall Easton's theorem
\begin{theorem*}
Assume $\text{GCH}$ and 
let $F\func \text{Reg} \to \Card$ be a class function such that
\begin{enumerate}
\item
	$\gl_1 < \gl_2 \implies F(\gl_1) \leq F(\gl_2)$,
\item
	$\cf F(\gl) > \gl$,
\end{enumerate}
then there is a cofinalities preserving  generic extension 
satisfying $\forall \gl \in \text{Reg}\ 2^\gl = F(\gl)$ and
$\text{SCH}$.
\end{theorem*}
Nowadays view on the power function is that we look for $\text{ZFC}$ theorems.
Forcing is used in order to show that some theorem is not possible or
to gain intuition on what is possible. Looking on Easton's theorem
from this point of view, it says: The only theorems we have regarding
the power function on the regular
cardinals are monotonicity and K\"{o}nig's lemma.

The question which, still, stands before us is how to include the singular
cardinals in an Easton like theorem. Even the formulation of such a theorem
is troublesome. For example, let $F\func \gw+1 \to \On$ be defined as:
$F(n) = n+1$, $F(\gw) = \gw_5$. Can we have a generic extension
in which $\forall \ga < \gw+1$ $2^{\ha_\ga} = \ha_{F(\ga)}$?
The answer to this question, as posed, is positive.
Neither our lack of knowledge in blowing $2^{\ha_\gw}$ above $\ha_{\gw_1}$
nor
Shelah's bound $2^{\ha_\gw} < \ha_{\gw_4}$ come into play
here due to the non-absoluteness of $\ha_{\gw_5}$.
By picking a successor $\ga < \gw_1$ and starting with
$\Col(\ga, \gw_5)$ we are in a position to invoke 
Gitik-Magidor forcing realizing the required power function.

Let us rephrase the question. For this we set $\gy(\gl, \gm) = \formula{
        (\gl < \ha_\gw \implies \gm = \gl^+)  \land
        (\gl = \ha_{\gw+1} \implies \gm = \ha_{\gw_5}) }$.
This time we ask: Does $\text{ZFC}$ + $\formula{ \forall \gl,\gm$ $\gy(\gl, \gm) \implies
        2^{\gl} = \gm}$ consistent?
And the answer is, of course, negative due to Shelah's bound.

So a possible attempt at including the singulars is:
For what formulas $\gy(\gl, \gm)$, satisfying 
        $\forall \gl \ \exists! \gm$
                $\gy(\gl, \gm)$,
the theory
   $\text{ZFC}$ + $\formula{\forall \gl,\gm\ \gy(\gl, \gm) \implies 2^\gl = \gm}$
is consistent?

Of course this theory should satisfy the $2$ `trivialities':
\begin{enumerate}
\item
	(Monotonicity)
	$\forall \gl_1, \gm_1, \gl_2, \gm_2$ 
	$\gl_1 < \gl_2 \land \gy(\gl_1, \gm_1) \land \gy(\gl_2, \gm_2)
		\implies
		\gm_1 \leq \gm_2$,
\item
	(K\"{o}nig's lemma)
	$\forall \gl, \gm$ $\gy(\gl, \gm) \implies \cf \gm > \gl$.
\end{enumerate}
Let us assume the theory satisfies:
\begin{enumerate}
\item
        (Galvin-Hajnal) If $\gy(\gl, \gm)$ and $\gl = \ha_\gh$ is a singular 
        strong limit of uncountable cofinality
                               then $\gm < \ha_{\gx^+}$ where
        $\gy(\power{\gh}, \gx)$,
\item
        (Silver) $\gw < \cf \gl < \gl$ $\land$
			$\setof {\gk < \gl} {\gy(\gk, \gk^+)}$ is stationary
			$\implies$
		$\gy(\gl, \gl^+)$,
\item
        (Shelah) $\gl$ is strong limit $\land$ $\gl = \ha_{\gx + \gz}$
		$\land$ $\gx < \ha_\gz$ 
		$\land$ $\gy(\gl, \gm)$ 
		$\implies$ 
			$\gm < \ha_{\gx + \power{\gz}^{+4}}$.
\end{enumerate}
Are these restrictions enough to ensure consistency of the theory?

Such a general theorem is beyond our knowledge at this time.
Note, this $\gy$ does not preclude infinite gaps
and we  miss a lot of information for such gaps.
Already for the first singular, $\ha_\gw$, assuming it is strong
limit,
we are lacking the technology
to blow $2^{\ha_\gw}$ above $\ha_{\gw_1}$ while 
the best known bound is $2^{\ha_\gw} < \ha_{\gw_4}$.

So let us restrict ourselves to finite gaps.
The forcing presented in this work showed the consistency of
the theory $\text{ZFC}$ $+$ $\forall \gl\ 2^\gl = \gl^{+n}$.
In the previous subsection
several generalizations and the principle limitations of it
were shown.
The main point was the appearance of a club with a fixed gap on it.
And the question is: Can we do without such a club?

For example, can the cardinals be partitioned into 2 stationary classes
such that on one of them we have gaps of $2$ and on the other $3$?
We do not know the answer (also) to this question.

In fact we do not know if it is possible to realize
a similar situation even below $\ha_{\gw_1}$.
That is, can we have stationary subsets of $\gw_1$, $S_1$ and $S_2$,
such that $S_1 \union S_2 = \gw_1$ and
$\ga \in S_1 \implies 2^{\ha_\ga} = \ha_{\ga + 2}$,
$\ga \in S_2 \implies 2^{\ha_\ga} = \ha_{\ga + 3}$?
Note, by Silver's theorem, we must have 
$2^{\ha_{\gw_1}} \leq \ha_{\gw_1 + 2}$.

It is interesting to note it looks as if the situation
$\ga \in S_1 \implies 2^{\ha_\ga} = \ha_{\ga + 1}$,
$\ga \in S_2 \implies 2^{\ha_\ga} = \ha_{\ga + 2}$
(in which case $2^{\ha_{\gw_1}} = \ha_{\gw_1 + 1}$)
is simpler to attack than the previous situation.

On the other hand, 
for a simple enough $\gy$
the club appearance is a must due to the following.
Let $\gy$ be an \emph{absolute} formula such that
$\forall \ga$ $\exists n< \gw$ $\gy(\ga, \ga+n)$
and assume the theory we work with is
$\text{ZFC}$ $+$ $\forall \ga,\gb\ \gy(\ga, \gb) \implies 2^{\ha_\ga} = \ha_\gb$.
(Note the change in $\gy$'s parameters: from cardinals
to ordinals).
We set $C_n = \setof{\ga} {\gy(\ga, \ga+n) }$.
Each $C_n$ is an $L$-class
and $\On = \bigunion \setof {C_n} {n < \gw}$.
Hence there is $n<\gw$
such that $C_n$ contains one of the $L$-indiscernibles hence all of them.
So $C_n$ contains a club.
\ifnum\article=0
\enlargethispage*{10pt}
\fi
\subsection{Regarding Our Forcing Notion}
We showed here only that $\gk$ is regular if $\cf \len(\Es) > \gk$.
	We have some preliminary work showing that if we have a repeat
	point, in the sense that $\PE = P_{\Es \restricted \gt}$,
	then $\gk$ remains measurable. Unlike \cite{MeselfPublishedPhdI},
	this is our only requirement.

\vspace{\baselineskip}
Let $G^*$ be $\ordered{\PE, \leq^*}$-generic. We think 
there is $H$ in
	$V[G^*]$ which is $j_{\Es}^\gx(\PE)$-generic
	over $M_{\Es}^\gx$ for a properly chosen $\gx$.
	So far, this is the closest we come
	to getting a generic by iteration.

	Moti Gitik had pointed that $\ordered{\PE, \leq^*}$ collapses
	$\gk^+$.

\vspace{\baselineskip}
We think it is of interest to find
	the connection between $M_{\Es}^\gx[H]$ and
	$\bigintersect_{\gx'<\gx} M_{\Es}^{\gx'}$.
\nocite{Cantor1}
\nocite{CofinalityAndNS}
\ifnum\article=0
\newpage
\fi
\bibliographystyle{plain}
\bibliography{carmi}
\end{document}